\documentclass[11pt]{article}
\usepackage{epsfig}
\usepackage{amssymb}
\usepackage{graphicx}
\usepackage{amsmath}
\usepackage{amsthm}
\usepackage{amscd}
\usepackage{arydshln}
\usepackage{mathtools}
\usepackage{units}
\usepackage{hyperref}
\usepackage[all,cmtip]{xy}
\usepackage{diagxy}
\usepackage{xymatrix}
\usepackage[all]{xy}
\usepackage{fancyhdr}
\usepackage[nottoc,numbib]{tocbibind}
\usepackage{graphicx}
\usepackage{mathrsfs}
\usepackage[perpage]{footmisc}
\usepackage{titlesec}
\usepackage{color}
\usepackage{geometry}
\usepackage{multicol}
\usepackage{slashed}

\newcommand{\tens}[1]{\mathbin{\mathop{\otimes}\limits_{#1}}}

\newtheoremstyle{basic}{11pt}{11pt}{}{}{\bfseries}{.}{0.5em}{}
\newtheoremstyle{proof}{11pt}{11pt}{}{}{\scshape}{:}{0.5em}{}

\newtheorem{prop}{Proposition}
\newtheorem{cor}{Corollary}
\newtheorem{lem}{Lemma}
\newtheorem{thm}{Theorem}

\theoremstyle{basic}
\newtheorem{defn}{Definition}

\theoremstyle{basic}
\newtheorem{exa}{Example}

\theoremstyle{basic}
\newtheorem{rem}{Remark}

\theoremstyle{proof}
\newtheorem*{prf}{Proof}

\makeatletter
\newcommand{\colim@}[2]{%
  \vtop{\m@th\ialign{##\cr
    \hfil$#1\operator@font colim$\hfil\cr
    \noalign{\nointerlineskip\kern1.5\ex@}#2\cr
    \noalign{\nointerlineskip\kern-\ex@}\cr}}%
}
\newcommand{\colim}{%
  \mathop{\mathpalette\colim@{\rightarrowfill@\textstyle}}\nmlimits@
}
\makeatother

\titleformat*{\section}{\large\bfseries}
\titleformat*{\subsection}{\normalsize\bfseries}
\newcommand{\harpoon}{\overset{\rightharpoonup}}

\begin{document}

\title{Non-commutative covering spaces}

\author{Clarisson Rizzie Canlubo\\
University of Copenhagen\\
\textit{clarisson@math.ku.dk}}

\maketitle

%\pagestyle{empty}{\tableofcontents}

%\newpage

\pagestyle{fancy}
\headheight=17pt
\textwidth=500pt
\fancyhead{}
\fancyhead[L]{\footnotesize{\leftmark}}
\fancyhead[R]{\thepage}
\renewcommand{\headrulewidth}{1pt}
\renewcommand{\footrulewidth}{1pt}
\fancyfoot{}
\fancyfoot[R]{}
\fancyfoot[L]{\scriptsize{Non-commutative covering spaces}}
\fancyfoot[C]{}

%%%%%%%%%%%%%%%%%%%%%%%%%%%%%%%%%%%%%%%%%%%%Text

\begin{abstract}
In this article, we will define non-commutative covering spaces using Hopf-Galois theory. We will look at basic properties of covering spaces that still hold for these non-commutative analogues. We will describe examples including coverings of commutative spaces and coverings of non-commutative tori.

\


\

\end{abstract}

\tableofcontents

\normalsize

\section{Introduction}\label{S1.0}

The fundamental group of a topological space $X$ is a very important and well-used invariant in classical geometry. It is defined as the group of homotopy classes of loops in $X$ based at some fixed point. This does not readily generalize to noncommutative spaces since there are no \textit{spaces} to work with let alone have a good notion of homotopy. A more subtle problem arise in algebraic geometry where spaces are too rigid to have a good notion of paths and homotopies. One can naively define the fundamental group of a scheme as the one we expect by simply considering the underlying topological space of that scheme. Explicitly, a loop is a Zariski-continuous map $\gamma:I\longrightarrow X$ where $I$ is the unit interval and $X$ is the scheme under consideration equipped with the Zariski topology. To stay in the realm of the algebraic category, we want to impose an algebraicity condition on $\gamma$ and eventually on $I$ but $I$ is far from being algebraic. We can relax this condition and settle for the usual euclidean topology on $I$. However, the fundamental group we will get is a rough one. To be precise, it cannot distinguish among affine schemes defined by integral domains which is more commonly known as affine varieties. Undeniably, they constitute an important class of schemes. More specifically, the fundamental group of the spectrum of an integral domain we will get by this naive definition is trivial. We will give a different formulation of the fundamental group analogous to Grothendieck's formulation in algebraic geometry \cite{b002}. He considered the category of finite Galois covering spaces of a scheme and defined the \'etale fundamental group as the inverse limit of the associated groups of deck transformations. To this end, we will develop in this paper the noncommutative analogue of covering spaces.

In the rest of this section, we will recall the necessary aspects of classical coverings spaces that we need. In section \ref{S2.0}, we will develop the necessary exposition for Hopf algebroids that we will use in the development of noncommutative covering spaces. Hopf algebroids should be taken as the noncommutative analogue of groupoids. We will enumerate examples of Hopf algebroids that will play a crucial part in the rest of the article. We will develop the necessary representation theoretic and Galois theoretic properties of Hopf algebroids in the remainder of that section. In section \ref{S3.0}, we will give our formulation of a noncommutative covering space and look at the appropriate notion of their equivalences. Section \ref{S4.0} deals with the structure of noncommutative coverings of commutative spaces. We will give a characterization of coverings of a point. We will give a reconstruction theorem that let us recover classical covering spaces when the algebraic objects involved are commutative. We will show that with centrality assumptions, noncommutative coverings of commutative spaces are bundles of coverings of a point. The remainder of that section tackles the special case of coverings with semisimple and cleft fibers. Section \ref{S5.0} deals with the noncommutative coverings of the noncommutative torus.

Let $X$ be a connnected and locally path connected space. An \textit{(unramified) covering} of $X$ is a space $Y$ together with a continuous surjection $Y\stackrel{p}{\longrightarrow} X$ such that any point $x\in X$ has an open neighborhood $U$ whose preimage is a disjoint union of homeomorphic copies of $U$, i.e. $p^{-1}(U)=\coprod\limits_{\alpha\in I}V_{\alpha}$ where each $V_{\alpha}$ are homeomorphic via $p$ to $U$. A \textit{ramified covering} of $X$ is a space $Y$ together with a continuous surjection $Y\stackrel{p}{\longrightarrow}X$ such that outside a nowhere dense set in $X$, $p$ is a unramified. The smallest such nowhere dense set is called the ramification locus of $p$. We will briefly refer to unramified coverings as coverings. The collection of all coverings of a given space $X$ forms a category $Cov(X)$. A morphism from a covering $Y\stackrel{p}{\longrightarrow}X$ to a covering $Z\stackrel{q}{\longrightarrow}X$ is a continuous map $Y\stackrel{r}{\longrightarrow}Z$ such that $p=q\circ r$. It is obvious that $r$ itself is a covering map. Given a covering $Y\stackrel{p}{\longrightarrow}X$, we can associate a group $Aut_{X}(Y)$. This group is called the group of \textit{deck transformations} of the covering $Y\stackrel{p}{\longrightarrow}X$. We say that $Y\stackrel{p}{\longrightarrow}X$ is \textit{Galois} if this group acts free and transitively on the fibers.

There is another useful description of covering spaces. The category $Cov(X)$ is equivalent to the functor category on the fundamental groupoid of $X$ with values in the category of sets. The latter category is easily seen to be complete and cocomplete.

Given a space $X$, let us denote by $\tilde{X}$ its universal cover and by $\pi_{1}(X,a)$ its fundamental group based at $a\in X$ (we will just write $\pi_{1}(X)$ if the group is independent of the base point, the case when for example $X$ is path-connected). We say that a (pointed) covering $(Y,b)\stackrel{p}{\rightarrow} (X,a)$ is \textit{intermediate} to the covering $(Z,c)\stackrel{q}{\rightarrow} (X,a)$ if there is a (pointed) map $(Z,c)\stackrel{\varphi}{\rightarrow} (Y,b)$ such that $p\circ\varphi=q$. This induces a partial order on the set of coverings of $X$ and incidentally gives a notion of equivalence. The group of autoequivalences of $(Y,b)\stackrel{p}{\rightarrow} (X,a)$ is precisely the group of deck transformations. We will be mostly interested in the case of connected covers $Y$. If $Aut_{(X,a)}(Y,b)$ acts transitively on the fibers of $Y\stackrel{p}{\rightarrow} X$, we call such covering \textsl{normal}. The covering map $p$ induces a monomorphism $p_{*}$ between fundamental groups. By the classification theorem for coverings (cf \cite{b001}), for every subgroup $G\leqslant \pi_{1}(X)$ there is a connected covering $(Y,b)\stackrel{p}{\rightarrow} (X,a)$ such that $p_{*}(\pi_{1}(Y))=G$. If $G$ is normal in $\pi_{1}(X)$ then $Aut_{X}Y=\pi_{1}(X)/G$. In this case, $Aut_{X}Y$ acts transitively on the fibers of $(Y,b)\stackrel{p}{\rightarrow}(X,a)$ and hence a normal covering. In general, $Aut_{X}Y=Nor(G)/G$ where $Nor(G)$ stands for the normalizer of $G$ in $\pi_{1}(X)$. Two coverings $(Y,b)\stackrel{p}{\rightarrow}(X,a)$ and $(Z,c)\stackrel{q}{\rightarrow}(X,a)$ are equivalent if the images of the fundamental groups of $Y$ and $Z$ coincides in $\pi_{1}(X)$. More generally, a covering $(Y,b)\stackrel{p}{\rightarrow}(X,a)$ associated to the subgroup $G_{Y}$ is intermediate to the covering $(Z,c)\stackrel{q}{\rightarrow}(X,a)$ associated to the subgroup $G_{Z}$ if $G_{Z}\subseteq G_{Y}$. The pointed coverings associated to conjugate subgroups are equivalent as coverings (rather than pointed coverings).

The above discussion will be briefly referred to as the Galois theory for coverings. In analogy with the Galois theory for fields, normal coverings correspond to Galois extensions, intermediate coverings correspond to intermediate extensions, and deck transformation groups correspond to Galois groups. Note that in classical Galois theory, a Galois extension is an algebraic extension which is both normal and separable. Since we are dealing with unramified coverings (which topologically corresponds to separable extensions), normal coverings automatically corresponds to Galois extensions. For a detailed exposition on this correspondence, one may consult Khovanskii \cite{b001}.

Let $\Sigma$ denote the category whose objects are $\pi_{1}(X)$-sets and whose morphisms are $\pi_{1}(X)$-equivariant maps. We will show that $Cov(X)$ and $\Sigma$ are equivalent categories. Given a covering $Y\stackrel{p}{\rightarrow} X$, there is an induced action of $\pi_{1}(X)$ on $p^{-1}(a)$. This defines a functor from $Cov(X)$ to $\Sigma$. Now, let $S$ be a $\pi_{1}(X)$-set. Let $S=\coprod_{\alpha\in I}S_{\alpha}$ be its decomposition into $\pi_{1}(X)$-orbits. Given a representative $s_{\alpha}$ of $S_{\alpha}$, we get a bijection between $S_{\alpha}$ and $\pi_{1}(X)/stab(s_{\alpha})$ by the orbit-stabilizer theorem. Then $stab(s_{\alpha})$ acts on $\tilde{X}$ and turns $\tilde{X}/stab(s_{\alpha})$ into a covering of $X$. Thus, we get $Y=\coprod_{\alpha\in I}\tilde{X}/stab(s_{\alpha})$ as a covering of $X$. This defines a functor inverse to the previous one. Note under this equivalence, the connected coverings are precisely the ones corresponding to homogenous $\pi_{1}(X)$-sets. A natural question to ask is whether $\pi_{1}(X)$ is completely determined by $Cov(X)$. The answer turns out to be affirmative using the following result:

\begin{thm}\label{thm1}
The group of natural automorphisms of the forgetful functor from $\pi_{1}(X)$-Sets to Sets is isomorphic to $\pi_{1}(X)$.
\end{thm}

By an automorphism $\alpha$ of the forgetful functor $\mathcal{F}$ we mean a family of automorphism $\mathcal{F}(S)\stackrel{\alpha_{S}}{\rightarrow}\mathcal{F}(S)$ such that for any morphism of $\pi_{1}(X)$-sets $S\stackrel{\sigma}{\rightarrow} T$, the following commutes

$$\xymatrix{
\mathcal{F}(S) \ar[d(1.6)]_-{\alpha_{S}} \ar[r(1.7)]^-{\mathcal{F}(\sigma)}
& & \mathcal{F}(T) \ar[d(1.6)]^-{\alpha_{T}}\\
& & \\
\mathcal{F}(S) \ar[r(1.6)]_-{\mathcal{F}(\sigma)}
& & \mathcal{F}(T).}$$

The theorem above is a very important theorem. Since the categories $Cov(X)$ and $\Sigma$ are equivalent, a problem concerning $Cov(X)$ is equally difficult in $\Sigma$. However, one can \textit{approximate} the answer by considering nice full subcategories of $\Sigma$ and the automorphism group of the forgetful functors for those subcategories. For example, if one considers the full subcategory of finite dimensional $\pi_{1}(X)$-representations, one gets the algebraic hull of $\pi_{1}(X)$. If one considers the full subcategory of finite $\pi_{1}(X)$-sets, the automorphism of the forgetful functor to sets is the profinite completion of $\pi_{1}(X)$. In a way, the above theorem serves as our guide in formulating the notion of a fundamental group and fundamental groupoid for a noncomutative space.

%%%%%%%%%%%%%%%%%%%%%%%%%%%%%%%%%%%%%%%%%%%%%%%%%%%%%%%%%%%%%%%%%%%%%%%%%%%%%%%%%%%%%%%%%%

\textbf{Acknowledgement.} I would like to thank my PhD supervisor Ryszard Nest for guiding me through my studies in noncommutative geometry and Ehud Meir and Olivier Gabriel for the valuable discussions that help me write this article. I would also like to thank DSF Grant, UP Diliman and the Center for Symmetry and Deformation, KU for the generous support.

%%%%%%%%%%%%%%%%%%%%%%%%%%%%%%%%%%%%%%%%%%%%%%%%%%%%%%%%%%%%%%%%%%%%%%%%%%%%%%%%%%%%%%%%%%

\section{Hopf algebroids}\label{S2.0}

\subsection{Definitions}\label{S2.1}

It has been a general consensus in noncommutative geometry that the analogue of groups are certain class of Hopf algebras called quantum groups. A \textit{Hopf algebra} $H$ (over a field $k$, or over a commutative unital ring in general) is an associative unital algebra $(H,m,1)$ together with algebra maps $H\stackrel{\Delta}{\longrightarrow} H\otimes H$ (\textit{coproduct}), $H\stackrel{\varepsilon}{\longrightarrow}k$ (\textit{counit}) and a linear map $H\stackrel{S}{\longrightarrow}H$ (\textit{antipode}) making the following diagrams commute.

\[ \xymatrix{
H \ar[d(.85)]_-{\Delta} \ar[r(1.7)]^-{\Delta}
& & H\otimes H \ar[d(.85)]^-{\Delta\otimes id}\\
H\otimes H \ar[r(1.6)]_-{id \otimes\Delta}
& & H\otimes H\otimes H \\
H \ar[rrd]^-{\Delta} \ar@{-}[r(1.8)]^-{\cong} \ar@{-}[d(.85)]_-{\cong} & & H\otimes k \\
k\otimes H & & H \otimes H \ar[l(1.5)]^-{\varepsilon \otimes id} \ar[u(.85)]_-{id \otimes \varepsilon} } \hspace{.25in}
\xymatrix@R=7mm{
& & H\otimes H \ar[r(1.6)]^-{S\otimes id} & & H\otimes H \ar[rd(1.6)]^-{m} & & \\
& & & & & & \\
H \ar[ru(1.8)]^-{\Delta} \ar[rd(1.8)]_-{\Delta} \ar[r(3.5)]^-{\varepsilon} & & & k \ar[r(2.5)]^-{1} & & & H \\
& & & & & & \\
& & H\otimes H \ar[r(1.6)]_-{id\otimes S} & & H\otimes H \ar[ru(1.6)]_-{m} & & \\
} \]

\noindent The two leftmost diagrams expresses the \textit{coassociativity} of $\Delta$ and its \textit{counitality} with respect to $\varepsilon$. With $\Delta$ and $\varepsilon$, $End(H)$ becomes a unital ring under convolution

\[ f\star g : H\stackrel{\Delta}{\longrightarrow} H\otimes H \stackrel{f\otimes g}{\longrightarrow} H\otimes H \stackrel{m}{\longrightarrow} H \]

\noindent with $H\stackrel{\varepsilon}{\longrightarrow} k \stackrel{1}{\longrightarrow} H$ as the unit. The diagram above involving $S$ expresses the fact that $S$ is the convolution inverse of $id$. From this, we immediately see that given a bialgebra $H$ (i.e. an algebra $H$ with coproduct and a counit which are algebra maps), there is at most one antipode which makes it into a Hopf algebra. We call a Hopf algebra a \textit{quantum group} if it has a bijective antipode. We will use Sweedler notation and Einstein summation convention all through out this paper. Explicitly, for any $h\in H$, instead of writing $\Delta(h)=\sum\limits_{i=1}^{n}(h_{1})_{i}\otimes (h_{2})_{i}$, we will write it as $\Delta(h)=h_{(1)}\otimes h_{(2)}$.

Recently, there has been great interest in Hopf-like structures in which the base ring is not necessarily commutative. Originally, we sought to develop the theory of noncommutative covering spaces using only Hopf algebras but there has been a great need to use a more general structure, one in which the base ring is possibly noncommutative. We will describe one which suits our purpose called a Hopf algebroid. A detailed discussion about Hopf algebroids, related structures and references can be found on B\"ohm \cite{bohm}. Through the remainder of this section, $k$ will be an associative, commutative unital ring and $R$ and $L$ will be associative unital $k$-algebras. A Hopf algebroid resembles a Hopf algebra$-$ it will have bialgebra-like structures defined over $R$ and $L$ and an antipode that relates them. Since we are mainly interested in the situation where Hopf algebroids are seen as further generalization of quantum groups, we will assume all throughout that Hopf algebroids have bijective antipodes. As it turns out, $R$ and $L$ will be anti-isomorphic $k$-algebras. However, for notational convenience it will be better to denote them accordingly, where we will use $R$ and $L$ to denote right and left structures, respectively. In addition, whenever we have a Hopf-like structure we will use Sweedler notation and Einstein summation convention to write down coproduct and coaction images. For a ring $R$, we will denote by $_{R}\mathcal{M}$ and by $\mathcal{M}_{R}$ the categories of left and right $R$-modules.

Before giving the definition of a Hopf algebroid, let us define first several intermediate structures. An $R$-\textit{ring} is a monoid object in the category of $R$-bimodules. Explicitly, an $R$-ring is a triple $(A,\mu,\eta)$ where $A\otimes_{R}A\stackrel{\mu}{\longrightarrow}A$ and $R\stackrel{\eta}{\longrightarrow}A$ are $R$-bimodule maps satisfying the associativity and unit axioms similar for algebras over commutative rings. A morphism of $R$-rings is a monoid morphism in category of $R$-bimodules. It is important to note that there is a bijection between $R$-rings $(A,\mu,\eta)$ and $k$-algebra morphisms $R\stackrel{\eta}{\longrightarrow}A$. Similar to the case of algebras over commutative rings, we can define modules over $R$-rings. For an $R$-ring $(A,\mu,\eta)$, a \textit{right} (resp. \textit{left}) $(A,\mu,\eta)$-\textit{module} is an algebra for the monad $-\otimes_{R}A$ (resp. $A\otimes_{R}-$) on the category $\mathcal{M}_{R}$ (resp. ${}_{R}\mathcal{M}$) of right (resp. left) modules over $R$.

We can dualize all the objects we have defined in the previous paragraph. An $R$-\textit{coring} is a comonoid in the category of $R$-bimodules, i.e a triple $(C,\Delta,\epsilon)$ where $C\stackrel{\Delta}{\longrightarrow}C\otimes_{R}C$ and $C\stackrel{\epsilon}{\longrightarrow}R$ are $R$-bimodule maps satisfying the coassociativity and counit axioms dual to those axioms satisfied by the structure maps of an $R$-ring. A morphism of $R$-corings is a morphism of comonoids. Given an $R$-coring $(C,\Delta,\epsilon)$, similar to coalgebras over commutative rings, we define a \textit{right} (resp. \textit{left}) $(C,\Delta,\epsilon)$-\textit{comodule} as a coalgebra for the comonad $-\otimes_{R}C$ (resp. $C\otimes_{R}-$) on the category $\mathcal{M}_{R}$ (resp. ${}_{R}\mathcal{M}$).

A \textit{right} (resp. \textit{left}) $R$-\textit{bialgebroid} $B$ is an $R\otimes_{k}R^{op}$-ring $(B,s,t)$ and an $R$-coring $(B,\Delta,\epsilon)$ satisfying:

\begin{enumerate}
\item [(a)] $R\stackrel{s}{\longrightarrow}B$ and $R^{op}\stackrel{t}{\longrightarrow}B$ are $k$-algebra maps with commuting images defining the $R\otimes_{k}R^{op}$-ring structure on $B$ which is compatible to the $R$-bimodule structure as an $R$-coring thru the following relation:
\[ r\cdot b \cdot r':=bs(r')t(r), \ \ (\text{resp.} \ r\cdot b \cdot r':=s(r)t(r')b,) \hspace{.15in} \forall r,r'\in R, b\in B. \]

\item [(b)] With the above $R$-bimodule structure on $B$ one can form $B\otimes_{R}B$. The coproduct $\Delta$ is required to corestrict to a $k$-algebra map to
\[ B\times_{R}B:=\left\{\sum\limits_{i}b_{i}\otimes_{R}b_{i}'\left|\sum\limits_{i}s(r)b_{i}\otimes_{R}b_{i}'=\sum\limits_{i}b_{i}\otimes_{R}t(r)b_{i}',\forall r\in R\right.\right\} \]
\noindent respectively,
\[ B \prescript{}{R}{\times} \ B:=\left\{\sum\limits_{i}b_{i}\otimes_{R}b_{i}'\left|\sum\limits_{i}b_{i}t(r)\otimes_{R}b_{i}'=\sum\limits_{i}b_{i}\otimes_{R}b_{i}'s(r),\forall r\in R\right.\right\}. \]

\item [(c)] The counit $B\stackrel{\epsilon}{\longrightarrow}R$ extends the right (resp. left) regular $R$-module structure on $R$ to a right (resp. left) $(B,s)$-module.
\end{enumerate}

\noindent A \textit{morphism} of $R$-bialgebroids is a morphism of $R\otimes R^{op}$-rings and $R$-corings.

\begin{rem}
\begin{enumerate}
\item[]
\item[(i)] The $k$-algebra maps $s$ and $t$ define a $k$-algebra map $\eta=s\otimes_{k}t$. As we have noted, such $k$-algebra uniquely determines an $R\otimes_{k}R^{op}$-ring structure on $B$. The maps $s$ and $t$ are called the \textit{source} and \textit{target} maps, respectively.

\item[(ii)] The $k$-submodule $B\times_{R}B$ (resp. $B\prescript{}{R}{\times} \ B$) of $B\otimes_{R}B$ is a $k$-algebra with factorwise multiplication. This is called the \textit{Takeuchi product}. The map $R\otimes_{k}R^{op}\longrightarrow B\times_{R}B$, $r\otimes_{k}r'\mapsto t(r')\otimes_{R}s(r)$ is easily seen to be a $k$-algebra morphism and hence, $B\times_{R}B$ is an $R\otimes_{k}R^{op}$-ring. The corestriction of $\Delta$ is an $R\otimes_{k}R^{op}$-bimodule map. Hence, $\Delta$ is an $R\otimes R^{op}$-ring map. The same is true for $B\prescript{}{R}{\times} \ B$.

\item[(iii)] The source map $s$ is a k-algebra map and so it defines a unique $R$-ring structure on $B$. The right version of condition (c) explicitly means that $r\cdot b:=\epsilon(s(r)b)$, $\forall r\in R, b\in B $ defines a right $(B,s)$-action on $R$.
\end{enumerate}
\end{rem}

We now have the necessary ingredients to define what a Hopf algebroid is.

\begin{defn}\label{D1}
Let $k$ be a commutative, associative unital ring and let $L$ and $R$ be associative $k$-algebras. A \textit{Hopf algebroid} $\mathcal{H}$ is a triple $\mathcal{H}=(\mathcal{H}_{L},\mathcal{H}_{R},S)$. $\mathcal{H}_{L}$ and $\mathcal{H}_{R}$ are bialgebroids having the same underlying $k$-algebra $H$. Specifically, $\mathcal{H}_{L}$ is a left $L$-bialgebroid with $(H,s_{L},t_{L})$ and $(H,\Delta_{L},\epsilon_{L})$ as its underlying $L\otimes_{k}L^{op}$-ring and $L$-coring structures. Similarly, $\mathcal{H}_{R}$ is a right $R$-bialgebroid with $(H,s_{R},t_{R})$ and $(H,\Delta_{R},\epsilon_{R})$ as its underlying $R\otimes_{k}R^{op}$-ring and $R$-coring structures. Let us denote by $\mu_{L}$ (resp. $\mu_{R}$) the multiplication on $(H,s_{L})$ (resp. $(H,s_{R})$). $S$ is a (bijective) $k$-module map $H\stackrel{S}{\longrightarrow}H$, called the \textit{antipode}. The compatibility conditions of these structures are as follows.

\begin{enumerate}
\item[(a)] the sources $s_{R},s_{L}$, targets $t_{R},t_{L}$ and counits $\epsilon_{R},\epsilon_{L}$ fit in commutative diagrams

\[ \xymatrix{
& & R^{op} \ar[lld]_-{t_{R}} \ar[rrd]^-{t_{R}} & & \\
H & & & & H \ar[lld]^-{\epsilon_{L}} \\
& & L \ar[llu]^-{s_{L}} \ar[lld]_-{s_{L}} \ar[rrd]^-{s_{L}} & & \\
H & & & & H \ar[lld]^-{\epsilon_{R}} \\
& & R^{op} \ar[llu]^-{t_{R}} & & \\
} \hspace{.25in}
\xymatrix{
& & R \ar[lld]_-{s_{R}} \ar[rrd]^-{s_{R}} & & \\
H & & & & H \ar[lld]^-{\epsilon_{L}} \\
& & L^{op} \ar[llu]^-{t_{L}} \ar[lld]_-{t_{L}} \ar[rrd]^-{t_{L}} & & \\
H & & & & H \ar[lld]^-{\epsilon_{R}} \\
& & R \ar[llu]^-{s_{R}} & & }\]

\item[(b)] the left- and right-regular comodule structures commute, i.e.

\[ \xymatrix{
H \ar[d(1.7)]_-{\Delta_{L}} \ar[r(1.7)]^-{\Delta_{R}}
& & H\tens{R} H \ar[d(1.6)]^-{\Delta_{L}\tens{R} id}\\
& & \\
H\tens{L} H \ar[r(1.6)]_-{id \tens{L}\Delta_{R}}
& & H\tens{L} H\tens{R} H } \hspace{0.75in}
\xymatrix{
H \ar[d(1.7)]_-{\Delta_{R}} \ar[r(1.7)]^-{\Delta_{L}}
& & H\tens{L} H \ar[d(1.6)]^-{\Delta_{R}\tens{L} id}\\
& & \\
H\tens{R} H \ar[r(1.6)]_-{id \tens{R}\Delta_{L}}
& & H\tens{R} H\tens{L} H } \]

\item[(c)] for all $l\in L, r\in R$ and for all $h\in H$ we have $S(t_{L}(l)ht_{R}(r))=s_{R}(r)S(h)s_{L}(l)$.

\item[(d)] $S$ is the convolution inverse of the identity map i.e., the following diagram commute

\[ \xymatrix{
& & H\tens{L} H \ar[rrrr]^-{S\tens{L} id} & & & & H\tens{L} H \ar[rrd]^-{\mu_{L}} & & \\
H \ar[rru]^-{\Delta_{L}}  \ar[rrrr]^-{\epsilon_{R}} & & & & R \ar[rrrr]^-{s_{R}} & & & & H \\
H \ar[rrd]_-{\Delta_{R}}  \ar[rrrr]_-{\epsilon_{L}} & & & & L \ar[rrrr]_-{s_{L}} & & & & H \\
& & H\tens{R} H \ar[rrrr]_-{id\tens{R} S} & & & & H\tens{R} H \ar[rru]_-{\mu_{R}} & & \\
} \]

\end{enumerate}
\end{defn}

\begin{rem}
\begin{enumerate}
\item[]
\item[(i)] In the constituent bialgebroids $\mathcal{H}_{R}$ and $\mathcal{H}_{L}$, the counits $\epsilon_{R}$ and $\epsilon_{L}$ extend the regular module structures on the base rings $R$ and $L$ to the $R$-ring $(H,s_{R})$ and to the $L$-ring $(H,s_{L})$, respectively. Equivalently, the counits extend the regular module structures on the base rings $R$ and $L$ to the $R^{op}$-ring $(H,t_{R})$ and to the $L^{op}$-ring $(H,t_{L})$. This particularly implies that the maps $s_{L}\circ \epsilon_{L}$, $t_{L}\circ \epsilon_{L}$, $s_{R}\circ \epsilon_{R}$ and $t_{R}\circ \epsilon_{R}$ are idempotents. This means that the images of $s_{R}$ and $t_{L}$ coincides in $H$. Same is true for the images of $s_{L}$ and $t_{R}$.

\item[(ii)] This implies that $\Delta_{L}$, apart from being an $L$-bimodule map, is also an $R$-bimodule map. Similarly, $\Delta_{R}$ is an $L$-bimodule map and so the diagrams in condition (b) make sense.

\item[(iii)] We can equip $H$ with two $(R,L)$-bimodule structures one using $t_{R}$ and $t_{L}$ and the other using $s_{R}$ and $s_{L}$. Condition (c) relates these two $(R,L)$-bimodules structures via the antipode $S$ which in turn makes the diagram in condition (d) defined.

\item[(iv)] The convolution structure condition (d) refers to a convolution structure one can define analogous to the one for linear maps from a coalgebra to an algebra. See \ref{S2.4} for this convolution structure.

\item[(v)] Let us note that condition (c) in the definition of a bialgebroid implies that $\epsilon_{L}\circ s_{L}:L\longrightarrow L$ is the identity. Similarly, $\epsilon_{R}\circ s_{R}:R\longrightarrow R$ is also the identity. Using condition (a) in the definition of a Hopf algebroid, we see that the following compositions define pairs of inverse $k$-algebra maps.

\[ \xymatrix{
L \ar[rr]^-{\epsilon_{R}\circ s_{L}} & & R^{op} \ar[rr]^-{\epsilon_{L}\circ t_{R}} & & L} \hspace{.5in}
\xymatrix{
R \ar[rr]^-{\epsilon_{L}\circ s_{R}} & & L^{op} \ar[rr]^-{\epsilon_{R}\circ t_{L}} & & R}\]

\noindent This is particular implies that $R$ and $L$ are anti-isomorphic $k$-algebras.

\item[(vi)] Since there are two coproducts involved in a Hopf algebroid, namely $\Delta_{L}$ and $\Delta_{R}$, we will use different Sweedler notations for their corresponding components. We will write $\Delta_{L}(h)=h_{[1]}\otimes_{L}h_{[2]}$ and $\Delta_{R}(h)=h^{[1]}\otimes_{R}h^{[2]}$ for $h\in H$.

\item[(vii)] With a fixed bijective antipode $S$, the constituent left- and right-bialgebroids of a Hopf algebroid determine each other, see for example \cite{bohm-szlachanyi}. In view of this and the fact that $L$ and $R$ are anti-isomorphic, in the sequel where we will be mainly interested with Hopf algebroids with bijective antipodes we will simply call $\mathcal{H}$ a Hopf algebroid \textit{over} $R$ instead of explicitly mentioning $L$.

\end{enumerate}
\end{rem}

\noindent Let $(\mathcal{H}_{L},\mathcal{H}_{R},S)$ and $(\mathcal{H}_{L}^{'},\mathcal{H}_{R}^{'},S^{'})$ be Hopf algebroids over $R$. An \textit{algebraic morphism}

\[ (\mathcal{H}_{L},\mathcal{H}_{R},S)\longrightarrow(\mathcal{H}_{L}^{'},\mathcal{H}_{R}^{'},S^{'}) \]

\noindent of Hopf algebroids is a pair $(\varphi_{L},\varphi_{R})$ of a left-bialgebroid morphism $\varphi_{L}$ and a right-bialgebroid morphism $\varphi_{R}$ for which the following diagrams commute

\[ \xymatrix{
\mathcal{H}_{L} \ar[dd]_-{\varphi_{L}} \ar[rr]^-{S}
& & \mathcal{H}_{R} \ar[dd]^-{\varphi_{R}}\\
& & \\
\mathcal{H}_{L}^{'} \ar[rr]_-{S^{'}}
& & \mathcal{H}_{R}^{'} } \hspace{0.75in}
\xymatrix{
\mathcal{H}_{R} \ar[dd]_-{\varphi_{R}} \ar[rr]^-{S}
& & \mathcal{H}_{L} \ar[dd]^-{\varphi_{L}}\\
& & \\
\mathcal{H}_{R}^{'} \ar[rr]_-{S^{'}}
& & \mathcal{H}_{L}^{'} } \]

\noindent and composition of such a pair is componentwise.

Let $R$ and $R^{'}$ be $k$-algebras and $(\mathcal{H}_{L},\mathcal{H}_{R},S)$ and $(\mathcal{K}_{L^{'}},\mathcal{K}_{R^{'}},S^{'})$ be Hopf algebroids over $R$ and $R^{'}$, respectively. In view of remark (vii) above, denote by $L=R^{op}$ and $L^{'}=(R^{'})^{op}$. A \textit{geometric morphism} $(\mathcal{H}_{L},\mathcal{H}_{R},S)\longrightarrow(\mathcal{K}_{L^{'}},\mathcal{K}_{R^{'}},S^{'})$ of Hopf algebroids is a pair $(f,\phi)$ of $k$-algerba maps $R\stackrel{f}{\longrightarrow}R^{'}$ and $H\stackrel{\phi}{\longrightarrow} K$, where $H,K$ denote the underlying $k$-algebra structures of the Hopf algebroids under consideration. These two maps satisfy the following compatibility conditions.

\begin{enumerate}

\item[(a)] $f$ and $\phi$ intertwines the source, target and counit maps of the left-bialgebroid structures of $\mathcal{H}$ and $\mathcal{K}$, i.e.

\[ \xymatrix{
H \ar[dd]_-{\phi} \ar[rr]^-{\epsilon^{H}_{L}}
& & L \ar[dd]^-{f}\\
& & \\
K \ar[rr]_-{\epsilon^{K}_{L}}
& & L^{'} } \hspace{0.15in}
\xymatrix{
L \ar[dd]_-{f} \ar[rr]^-{t^{H}_{L}}
& & H \ar[dd]^-{\phi}\\
& & \\
L^{'} \ar[rr]_-{t^{K}_{L}}
& & K } \hspace{0.15in}
\xymatrix{
L \ar[dd]_-{f} \ar[rr]^-{s^{H}_{L}}
& & H \ar[dd]^-{\phi}\\
& & \\
L^{'} \ar[rr]_-{s^{K}_{L}}
& & K. } \]

\noindent Same goes for the source, target and counit maps of the right-bialgebroid structures.

\item[(b)] In view of condition $(a)$, the $k$-bimodule map $\phi\otimes_{k}\phi$ defines $k$-bimodule maps

\[  \xymatrix{ H \prescript{}{L}{\otimes} \ H \ar[rr]^-{\phi\prescript{}{f}{\otimes} \ \phi} && K \prescript{}{L^{'}}{\otimes} \ K, } \hspace{.25in} \xymatrix{ H \otimes_{R} \ H \ar[rr]^-{\phi\ \otimes_{f} \phi} && K \otimes_{R^{'}} \ K. } \]

\noindent We then require that the following diagrams commute

\[  \xymatrix{ H \prescript{}{L}{\otimes} \ H \ar[rr]^-{\phi\prescript{}{f}{\otimes} \ \phi} \ar[dd]_-{\mu^{H}_{L}} && K \prescript{}{L^{'}}{\otimes} \ K \ar[dd]^-{\mu^{K}_{L}} \\
&& \\
H \ar[rr]_-{\phi} && K } \hspace{.25in}
\xymatrix{ H \otimes_{R} \ H \ar[rr]^-{\phi\ \otimes_{f} \phi} \ar[dd]_-{\mu^{H}_{R}} && K \otimes_{R^{'}} \ K \ar[dd]^-{\mu^{K}_{R}} \\
&& \\
H \ar[rr]_-{\phi} && K} \]

\item[(c)] Also by of condition $(a)$, the $k$-bimodule maps $\phi\prescript{}{f}{\otimes} \ \phi$ and $\phi\ \otimes_{f} \phi$ of condition $(b)$ further define $k$-bimodule maps

\[  \xymatrix{ H \prescript{}{L}{\times} \ H \ar[rr]^-{\phi\prescript{}{f}{\times} \ \phi} && K \prescript{}{L^{'}}{\times} \ K, } \hspace{.25in} \xymatrix{ H \times_{R} \ H \ar[rr]^-{\phi\ \times_{f} \phi} && K \times_{R^{'}} \ K. } \]

\noindent We then require that the following diagrams commute.

\[  \xymatrix{ H \ar[rr]^-{\phi} \ar[dd]_-{\Delta^{H}_{L}} && K \ar[dd]^-{\Delta^{K}_{L}} \\
&& \\
H \prescript{}{L}{\times} \ H \ar[rr]_-{\phi\prescript{}{f}{\times} \ \phi} && K \prescript{}{L^{'}}{\times} \ K } \hspace{.25in}
\xymatrix{ H \ar[rr]^-{\phi} \ar[dd]_-{\Delta^{H}_{R}} && K \ar[dd]^-{\Delta^{K}_{R}} \\
&& \\
H \times_{R} \ H \ar[rr]_-{\phi\ \times_{f} \phi} && K \times_{R^{'}} \ K } \]

\item[(d)] $\phi$ intertwines the antipodes of $\mathcal{H}$ and $\mathcal{K}$, i.e. $\phi\circ S_{H}=S_{K}\circ\phi$.

\end{enumerate}

\begin{rem}
\begin{enumerate}
\item[]

\item[(i)] For a $k$-algebra $R$, let us denote by $HALG^{alg}(R)$ the category whose objects are Hopf algebroids over $R$ and morphisms are algebraic morphisms. For a fixed $k$, let us denote by $HALG^{geom}(k)$ the category whose objects are Hopf algebroids over $k$-algebras and morphisms are geometric morphisms. The existence of these two naturally defined categories reflect the fact that Hopf algebroids are generalization of both Hopf algebras and groupoids.

\item[(ii)] Equip $R^{e}$ with the Hopf algebroid structure defined in example 5 of the next section. Let $(\mathcal{H}_{L},\mathcal{H}_{R},S)$ be a Hopf algebroid over $R$. Then the unit maps $\eta_{L},\eta_{R}$ together with the identity map on $R$ define geometric morphisms $(id,\eta_{L}):R^{e}\longrightarrow \mathcal{H}$ and $(id,\eta_{R}):R^{e}\longrightarrow \mathcal{H}$.
\end{enumerate}
\end{rem}

\subsection{Examples and properties}\label{S2.2}

In this section, we will enumerate examples of Hopf algebroids that will play a crucial role in the rest of the article.

\begin{exa}\label{exa1} \textbf{Hopf algebras}. A Hopf algebra $H$ over the commutative unital ring $k$ gives an example of a Hopf algebroid. Here, we take $R=L=k$ as $k$-algebras, take $s_{L}=t_{L}=s_{R}=t_{R}=\eta$ to be the source and target maps, set $\epsilon_{L}=\epsilon_{R}=\epsilon$ to be the counits, and $\Delta_{L}=\Delta_{R}=\Delta$ to be the coproducts.
\end{exa}

\

\begin{exa}\label{exa2} \textbf{Coupled Hopf algebras}. It might be tempting to think that Hopf algebroids for which $R=L=k$ must be Hopf algebras. This is not entirely the case. We will give a general set of examples for which this is not true. Two Hopf algebra structures $H_{1}=(H,m_{1},\eta_{1},\Delta_{1},\epsilon_{1},S_{1})$ and $H_{2}=(H,m_{2},\eta_{2},\Delta_{2},\epsilon_{2},S_{2})$ over the same $k$-module $H$ are said to be \textit{coupled} if

\begin{enumerate}
\item[(i)] there exists a $k$-module map $C:H_{1}\longrightarrow H_{2}$, called the \textit{coupling map} such that

\[ \xymatrix{
& & H\otimes H \ar[r(3)]^-{C\otimes id} & & & & H\otimes H \ar[rd(1.6)]^-{m_{1}} & & \\
& & & & & & & & \\
H \ar[ru(1.8)]^-{\Delta_{1}} \ar[rd(1.8)]_-{\Delta_{2}} \ar@<1ex>[rrrr]^-{\epsilon_{2}} \ar@<-1ex>[rrrr]_-{\epsilon_{1}} & & & & k \ar[rrrr]_-{\eta} & & & & H \\
& & & & & & & & \\
& & H\otimes H \ar[r(3)]_-{id\otimes C} & & & & H\otimes H \ar[ru(1.6)]_-{m_{2}} & & \\
} \]

\noindent commutes, and

\item[(ii)] the coproducts $\Delta_{1}$ and $\Delta_{2}$ commute.
\end{enumerate}

\noindent Coupled Hopf algebras give rise to Hopf algebroids over $k$. The left $k$-bialgebroid is the underlying bialgerba of $H_{1}$ while the right $k$-bialgebroid is the underlying bialgebra of $H_{2}$. The coupling map plays the role of the antipode.

Let us give examples of coupled Hopf algebras. Connes and Moscovici constructed \textit{twisted} antipodes in \cite{cm}. Let us show that such a twisted antipode is a coupling map for some coupled Hopf algebras. Let $H=(H,m,1,\Delta,\epsilon,S)$ be a Hopf algebra. Take $H_{1}=H$ as Hopf algebras. Let $\sigma:H\longrightarrow k$ be a character. Define $\Delta_{2}:H\longrightarrow H\otimes H$ by $h\mapsto h_{(1)}\otimes\sigma(S(h_{(2)}))h_{(3)}$. Take $\epsilon_{2}=\sigma$. Define $S_{2}:H\longrightarrow H$ by $h\mapsto \sigma(h_{(1)})S(h_{(2)})\sigma(h_{(3)})$. Note the Sweedler-legs of $h$ appearing in the definition of $S_{2}$ is the one provided by $\Delta$ and not by $\Delta_{2}$. Then, $H_{2}=(H,m,1,\Delta_{2},\epsilon_{2},S_{2})$ is a Hopf algebra coupled with $H_{1}$ by the coupling map $S^{\sigma}:H\longrightarrow H$ defined by $h\mapsto \sigma(h_{(1)})S(h_{(2)})$.
\end{exa}

\

\begin{exa}\label{exa3} \textbf{Groupoid algebras}. Given a small groupoid $\mathcal{G}$ with finitely many objects and a commutative unital ring $k$, we can construct what is called the groupoid algebra of $\mathcal{G}$ over $k$, denoted by $k\mathcal{G}$. For such a groupoid $\mathcal{G}$, let us denote by $\mathcal{G}^{(0)}$ its set of objects, $\mathcal{G}^{(1)}$ its set of morphisms, $s,t:\mathcal{G}^{(1)}\longrightarrow \mathcal{G}^{(0)}$ the source and target maps, $\iota:\mathcal{G}^{(0)}\longrightarrow \mathcal{G}^{(1)}$ the unit map, $\nu:\mathcal{G}^{(1)}\longrightarrow \mathcal{G}^{(1)}$ the inversion map, $\mathcal{G}^{(2)}=\mathcal{G}^{(1)}\prescript{}{t}\times_{s} \ \mathcal{G}^{(1)}$ the set of composable pairs of morphisms, and $m:\mathcal{G}^{(2)}\longrightarrow \mathcal{G}^{(1)}$ the partial composition. The groupoid algebra $k\mathcal{G}$ is the $k$-algebra generated by $\mathcal{G}^{(1)}$ subject to the relation

\[ ff^{'}= 
\begin{cases}
    f\circ f^{'},& \text{if } f,f^{'} \ \text{are composable}\\
		& \\
    0,              & \text{otherwise}
\end{cases}
\]

for $f,f^{'}\in \mathcal{G}^{(1)}$. The groupoid algebra $k\mathcal{G}$ is a Hopf algebroid as folows. The base algebras $R$ and $L$ are both equal to $k\mathcal{G}^{(0)}$ and the two bialgebroids $H_{R}$ and $H_{L}$ are isomorphic as bialgebroids with underlying $k$-module $k\mathcal{G}^{(1)}$. The partial groupoid composition $m$ dualizes and extends to a multiplication $m:k\mathcal{G}^{(1)}\otimes k\mathcal{G}^{(1)}\longrightarrow k\mathcal{G}^{(1)}$ which then factors through the canonical surjection $k\mathcal{G}^{(1)}\otimes k\mathcal{G}^{(1)}\longrightarrow k\mathcal{G}^{(1)}\otimes_{k\mathcal{G}^{(0)}} k\mathcal{G}^{(1)}$ to give the product $k\mathcal{G}^{(1)}\otimes_{k\mathcal{G}^{(0)}} k\mathcal{G}^{(1)}\longrightarrow k\mathcal{G}^{(1)}$. The source and target maps $s,t$ of the groupoid give the source and target maps $s,t:k\mathcal{G}^{(0)}\longrightarrow k\mathcal{G}^{(1)}$, respectively. The unit map gives the counit map $\epsilon:k\mathcal{G}^{(1)}\longrightarrow k\mathcal{G}^{(0)}$. Finally, the inversion map gives the antipode map $S:k\mathcal{G}^{(1)}\longrightarrow k\mathcal{G}^{(1)}$.

With this example, we immediately see that if the groupoid is a group, the construction above gives a Hopf algebra over $k$. This justifies the name Hopf algebroid. Just like in the case for groups, there is a dual construction to the one we presented here. We will present that in the beginning of section \ref{S3.1}.
\end{exa}

\

\begin{exa}\label{exa4} \textbf{Weak Hopf algebras}. Another structure that generalize Hopf algebras, called weak Hopf algebras, also are Hopf algebroids. Explicitly, a weak Hopf algebra $H$ over a commutative unital ring $k$ is a unitary associative algebra together with $k$-linear maps $\Delta:H\longrightarrow H\otimes H$ (weak coproduct), $\epsilon:H\longrightarrow k$ (weak counit) and $S:H\longrightarrow H$ (weak antipode) satisfying the following axioms:

\begin{enumerate}
\item[(i)] $\Delta$ is multiplicative, coassociative, and weak-unital, i.e.
\[(\Delta(1)\otimes 1)(1\otimes \Delta(1))=\Delta^{(2)}(1)=(1\otimes\Delta(1))(\Delta(1)\otimes 1),\]

\item[(iii)] $\epsilon$ is counital, and weak-multiplicative, i.e. for any $x,y,z\in H$
\[ \epsilon(xy_{(1)})\epsilon(y_{(2)}z)=\epsilon(xyz)=\epsilon(xy_{(2)})\epsilon(y_{(1)}z),\]

\item[(v)] for any $h\in H$, $S(h_{(1)})h_{(2)}S(h_{(3)})=S(h)$ and
\[ h_{(1)}S(h_{(2)})=\epsilon(1_{(1)}h)1_{(2)}, \hspace{.5in} S(h_{(1)})h_{(2)}=1_{(1)}\epsilon(h1_{(2)}) \]
\end{enumerate}

Let us sketch a proof why a weak Hopf algebra $H$ is a Hopf algebroid. Consider the maps $p_{R}:H\longrightarrow H$, $h\mapsto 1_{(1)}\epsilon(h1_{(2)})$ and $p_{L}:H\longrightarrow H$, $h\mapsto \epsilon(1_{(1)}h)1_{(2)}$. By $k$-linearity and weak-multiplicativity of $\epsilon$, $p_{R}$ and $p_{L}$ are idempotents.

Multiplicativity and coassiociativity of $\Delta$ and counitality of $\epsilon$ implies that for any $h\in H$,

\[ h_{(1)}\otimes p_{L}(h_{(2)})=1_{(1)}h\otimes 1_{(2)} \hspace{.5in} p_{R}(h_{(1)})\otimes h_{(2)}=1_{(1)}\otimes h1_{(2)}.\]

\noindent Now, using these relations and coassiociativity of $\Delta$ we get

\[ 1_{(1)}1_{(1')}\otimes 1_{(2)} \otimes 1_{(2')} = 1_{(1')(1)}\otimes p_{L}(1_{(1')(2)})\otimes 1_{(2')}=1_{(1)}\otimes p_{L}(1_{(2)})\otimes 1_{(3)}\]

\[ 1_{(1)}\otimes 1_{(1')}\otimes 1_{(2)}1_{(2')} = 1_{(1)}\otimes p_{L}(1_{(2)(1)})\otimes 1_{(2)(2)}=1_{(1)(1)}\otimes p_{L}(1_{(1)(2)})\otimes 1_{(2)}\]

\noindent Thus, the first tensor factor of the left-hand side of the first equation above is in the image of $p_{R}$. Similarly, the last tensor factor of the left-hand side of the second equation above is in the image of $p_{L}$. Clearly, $p_{R}(1)=p_{L}(1)=1$. Hence, the images of $p_{R}$ and $p_{L}$ are unitary subalgebras of $H$. Denote these subalgebras by $R$ and $L$, respectively. By the weak-unitality of $\Delta$ we see that these subalgebras are commuting subalgebras of $H$.

Taking the source map $s$ as the inclusion $R\longrightarrow H$ and the target map as $t:R^{op}\longrightarrow H$, $r\mapsto\epsilon(r1_{(1)})1_{(2)}$ equips $H$ with an $R\otimes_{k}R^{op}$-ring structure. Taking $\epsilon_{R}=p_{R}$ and $\Delta_{R}$ as the composition

\[ \xymatrix{H \ar[rr]^-{\Delta}  & & H\otimes_{k}H \ar@{->>}[rr] & & H\otimes_{R}H }\]

\noindent equips $H$ with an $R$-coring structure $(H,\Delta_{R},\epsilon_{R})$. The ring and coring structures just constructed gives $H$ a structure of right $R$-bialgebroid $H_{R}$.

Using $R^{op}$ in place of $R$ in the above construction, we get a left $R^{op}$-bialgebroid $H_{R^{op}}$. Together with the right $R$-bialgebroid constructed and the existing weak antipode $S$, we get a Hopf algebroid $(H_{R^{op}},H_{R},S)$.
\end{exa}

\

Weak Hopf algebras also has a well-understood representation theory. Given a weak Hopf algebra $H$ over a field $k$, the category ${}_{H}\mathcal{M}$ of finitely-generated left modules over $H$ is a fusion category. A \textit{fusion category} $\mathcal{C}$ over $k$ is a $k$-linear rigid semisimple category with finitely-many inequivalent simple objects such that the hom-spaces are finite-dimensional and the endomorphism algebra of the unit object $\mathbb{1}_{\mathcal{C}}$ is $k$. By Tannaka duality, any fusion category is equivalent to a module category of a weak Hopf algebra. This phenomenon has a nice symmetry. Similar to Hopf algebras, the dual $H^{*}$ of a finitely generated weak Hopf algebra $H=(H,m,1,\Delta,\varepsilon,S)$ has a natural weak Hopf algebra structure. Using this idea, one can show that the category $\mathcal{M}^{H}$ of finitely-generated right comodules over $H$ is a fusion category as well.

\

\begin{exa}\label{exa5} \textbf{Group algebras over noncommutative rings.} One of the most studied yet mysterious class of a Hopf algebras are group algebras over commutative rings. In this section, we will show a similar construction of a group algebra over a noncommutative base ring and see that such is a Hopf algebroid. This further justifies the banner of Hopf algebroids being a generalization of Hopf algebras over noncommutative rings.

Let $A$ be an associative unital algebra over a commutative ring $k$. Denote by $A^{e}=A\otimes A^{op}$ its universal enveloping algebra. Consider a finite group $G$ acting on $A$ via $G\stackrel{\alpha}{\longrightarrow}Aut(A)$. This action extends to a $kG$-module structure on $A^{e}$ via the usual coproduct on $kG$. Consider the smash product algebra $A^{e}\#kG$. The underlying $k$-module of this algebra is $A^{e}\otimes kG$. The multiplication is defined as

\[ \left(\sum \left(a^{1}\otimes a^{2}\right)\# g\right)\left(\sum \left(b^{1}\otimes b^{2}\right)\# h\right)=\sum \left(a^{1}\otimes a^{2}\right)\alpha_{g}\left( b^{1}\otimes b^{2}\right)\# gh \]

\noindent Note that this construction generalize to the case of a bialgebra $H$ in place of $kG$ where the two appearance of $g$'s in the defining relation for the multiplication is played by the legs of coproduct applied to the appropriate tensor factor. If the action of $G$ is trivial, we get the algebra $A^{e}G$ which we call the group algebra of $G$ over $A^{e}$. Let us show that $A^{e}G$ is a Hopf algebroid over $A$. The right $A$-bialgebroid structure consists of $A^{e}G$ as the underlying $k$-module. The right source $s_{R}$, target $t_{R}$ and counit maps $\epsilon_{R}$ are

\[ \xymatrix @R=2mm {A \ar[r]^-{s_{R}} & A\otimes A^{op} \# kG \\
a \ar@{|->}[r] & \left(a \otimes 1\right) \# e} \hspace{.5in} \xymatrix @R=2mm {A \ar[r]^-{t_{R}} & A\otimes A^{op} \# kG \\
a \ar@{|->}[r] & \left(1 \otimes a\right) \# e} \hspace{.5in} \xymatrix@R=2mm{A\otimes A^{op} \# kG \ar[r]^-{\epsilon_{R}} & A. \\
\left(a \otimes a^{'}\right) \# g \ar@{|->}[r] & aa^{'}}\]

\noindent where $e$ stands for the identity element of $G$. The right coproduct $\Delta_{R}$ is the following map.

\[ \xymatrix@R=3mm{ A\otimes A^{op} \# kG \ar[rr]^-{\Delta_{R}} & & \left(A\otimes A^{op} \# kG\right) \tens{A} \left(A\otimes A^{op} \# kG\right) \\
\left(a\otimes a^{'}\right) \# g \ar@{|->}[rr] & & \left(1\otimes a^{'}\right) \# g \tens{A} \left(a\otimes 1\right) \# g } \]

The left $A$-bialgebroid is the opposite co-opposite of the right $A$-bialgebroid we just constructed. The map

\[ \xymatrix@R=2mm{ A\otimes A^{op} \# kG \ar[rr]^-{S} & & A^{op}\otimes A \# kG \\
\left(a \otimes a^{'}\right) \# g \ar@{|->}[rr] & & \left(a^{'} \otimes a\right) \# g^{-1}  } \]

\noindent is the antipode. In particular, taking $G$ to be the trivial group makes $A^{e}$ a Hopf algebroid over $A$. Any of the underlying coring structures of $A^{e}$ is what is commonly known in the literature as the canonical coring associated to $A$. With this, we call $A^{e}$ the \textit{canonical Hopf algebroid} over $A$.
\end{exa}

\subsection{Representation theory of Hopf algebroids and their descent}\label{S2.3}

In this section, we will look at representations of Hopf algebroids. Towards the end of the section, we will look at the descent theoretic aspect of a special class of modules over Hopf algebroids, the so called relative Hopf modules. Let $\mathcal{H}=(\mathcal{H}_{L},\mathcal{H}_{R},S)$ be a Hopf algebroid with underlying $k$-module $H$. $H$ carries both a left $L$-module sctructure and a left $R$-module structure via the maps $s_{L}$ and $t_{R}$, respectively. A \textit{right} $\mathcal{H}$-\textit{comodule} $M$ is a right $L$-module and a right $R$-module together with a right $\mathcal{H}_{R}$-coaction $\rho_{R}:M\longrightarrow M\otimes _{R}H$ and a right $\mathcal{H}_{L}$-coaction $\rho_{L}:M\longrightarrow M\otimes_{L}H$ such that $\rho_{R}$ is an $\mathcal{H}_{L}$-comodule map and $\rho_{L}$ is an $\mathcal{H}_{R}$-comodule map.

For the coaction $\rho_{R}$, let us use the following Sweedler notation:

\[ \rho_{R}(m) = m^{[0]}\tens{R} m^{[1]} \]

\noindent and for the coaction $\rho_{L}$, let us use the following Sweedler notation:

\[ \rho_{L}(m) = m_{[0]} \tens{L} m_{[1]}. \]

\noindent With these notations, the conditions above explicitly means that for all $m\in M$, $l\in L$ and $r\in R$ we have

\[ (m\cdot l)^{[0]}\tens{R}(m\cdot l)^{[1]}=\rho_{R}(m\cdot l)=m^{[0]}\tens{R} t_{L}(l)m^{[1]} \]

\[ (m\cdot r)_{[0]}\tens{L}(m\cdot r)_{[1]}=\rho_{L}(m\cdot r)=m_{[0]}\tens{L} m_{[1]}s_{R}(r). \]

\noindent We further require that the two coactions satify the following commutative diagrams

\[ \xymatrix{
M \ar[d(1.7)]_-{\rho_{R}} \ar[r(1.7)]^-{\rho_{L}}
& & M\tens{L} H \ar[d(1.6)]^-{\rho_{R}\tens{L} id}\\
& & \\
M\tens{R} H \ar[r(1.6)]_-{id \tens{R}\Delta_{L}}
& & M\tens{R} H\tens{L} H } \hspace{0.75in}
\xymatrix{
M \ar[d(1.7)]_-{\rho_{L}} \ar[r(1.7)]^-{\rho_{R}}
& & M\tens{R} H \ar[d(1.6)]^-{\rho_{L}\tens{R} id}\\
& & \\
M\tens{L} H \ar[r(1.6)]_-{id \tens{L}\Delta_{R}}
& & M\tens{L} H\tens{R} H } \]

We will denote by $\mathcal{M}^{\mathcal{H}}$ the category of right $\mathcal{H}$-comodules. Symmetrically, we can define left $\mathcal{H}$-comodules and we denote the category of a such by $^{\mathcal{H}}\mathcal{M}$.

Comodules over Hopf algebroids are comodules over the constituent bialgebroids. Thus, one can speak of two different coinvariants, one for each bialgebroid. For a given right $\mathcal{H}$-comodule $M$, they are defined as follows:

\[ M^{co \ \mathcal{H}_{R}} = \left\{m\in M\left| \ \rho_{R}(m)=m\tens{R} 1\right.\right\}, \]

\[ M^{co \ \mathcal{H}_{L}} = \left\{m\in M\left| \ \rho_{L}(m)=m\tens{L} 1\right.\right\}. \]

\noindent In the general case, we have $M^{co \ \mathcal{H}_{R}}\subseteq M^{co \ \mathcal{H}_{L}}$. But in our case, where we assume $S$ is bijective these two spaces coincide. This will be important in the formulation of Galois theory for Hopf algebroids. To see that these coinvariants coincide, consider the following map

\[ \Phi_{M}:M\tens{R}H\longrightarrow M\tens{L}H \]
\[ m\tens{R} h \mapsto \rho_{L}(m)\cdot S(h) \]

\noindent Here, $H$ acts on the right of $M\otimes_{L}H$ through the second factor. If $m\in M^{co \ \mathcal{H}_{R}}$, then we have

\begin{eqnarray}
\nonumber\label{} \rho_{L}(m)&=&\rho_{L}(m)\cdot S(h) = \Phi_{M}(m\tens{R}1) = \Phi_{M}(\rho_{R}(m))\\
\nonumber &=& \Phi_{M}(m^{[0]}\tens{R} m^{[1]}) = \rho_{L}(m^{[0]})\cdot S(m^{[1]})\\
\nonumber &=& (m^{[0]}_{[0]}\tens{L} m^{[0]}_{[1]})\cdot S(m^{[1]}) = m^{[0]}_{[0]}\tens{L} m^{[0]}_{[1]}S(m^{[1]})\\
\nonumber &=& m_{[0]}\tens{L} m^{[0]}_{[1]}S(m^{[1]}_{[1]})= m_{[0]}\tens{L} s_{L}(\epsilon_{L}(m_{[1]}))\\
\nonumber &=& m_{[0]}s_{L}(\epsilon_{L}(m_{[1]}))\tens{L} 1 = m\tens{L} 1\\
\nonumber
\end{eqnarray}

\noindent This shows the inclusion $M^{co \ \mathcal{H}_{R}}\subseteq M^{co \ \mathcal{H}_{L}}$. To show the other inclusion, one can run the same computation but using the inverse of $\Phi_{M}$ which is the following map

\[ \Phi_{M}^{-1}:M\tens{L}H\longrightarrow M\tens{R}H \]
\[ m\tens{L}h\mapsto S^{-1}(h)\cdot\rho_{R}(m). \]

\noindent In this case, we can simply write $M^{co \ \mathcal{H}}$ for $M^{co \ \mathcal{H}_{R}}=M^{co \ \mathcal{H}_{L}}$ and refer to it as the $\mathcal{H}$-coinvariants of $M$ instead of distinguishing the $\mathcal{H}_{R}$- from the $\mathcal{H}_{L}$-coinvariants, unless it is necessary to do so.

Let us now discuss monoid objects in $\mathcal{M}^{\mathcal{H}}$. They are called $\mathcal{H}$-comodule algebras. A right $\mathcal{H}$-\textit{comodule algebra} is an $R$-ring $(M,\mu,\eta)$ such that $M$ is a right $\mathcal{H}$-comodule and $\eta:R\longrightarrow M$ and $\mu:M\otimes_{r}M\longrightarrow M$ are $\mathcal{H}$-comodule maps. Using Sweedler notation for coactions, this explicitly means that for any $m,n\in M$ we have

\[ (mn)^{[0]}\tens{R}(mn)^{[1]}=\rho_{R}(mn)=m^{[0]}n^{[0]}\tens{R} m^{[1]}n^{[1]}, \]

\[ (mn)_{[0]}\tens{L}(mn)_{[1]}=\rho_{L}(mn)=m_{[0]}n_{[0]}\tens{L} m_{[1]}n_{[1]}, \]

\[ 1_{M}^{[0]}\tens{R}1_{M}^{[1]}=\rho_{R}(1_{M})=1_{M}\tens{R} 1_{H}, \]

\[ (1_{M})_{[0]}\tens{L}(1_{M})_{[1]}=\rho_{L}(1_{M})=1_{M}\tens{L} 1_{H}. \]

Let $M$ be a right $\mathcal{H}$-comodule algebra. A \textit{right-right relative} $(M,\mathcal{H})$-\textit{Hopf module} $W$ is a right module of the $R$-ring $M$ such that the module structure $(\cdot):W\otimes_{R}M\longrightarrow W$ is a right $\mathcal{H}$-comodule map, i.e.

\[ (w\cdot m)^{[0]} \tens{R} (w\cdot m)^{[1]} = w^{[0]}\cdot m^{[0]} \tens{R} w^{[1]}m^{[1]} \]

\[ (w\cdot m)_{[0]} \tens{L} (w\cdot m)_{[1]} = w_{[0]}\cdot m_{[0]} \tens{L} w_{[1]}m_{[1]} \]

\noindent for any $w\in W$ and $m\in M$. We denote by $\mathcal{M}^{\mathcal{H}}_{M}$ the category of right-right relative $(M,\mathcal{H})$-Hopf modules. One can symmetrically define left-right, left-left and right-left relative Hopf modules, whose categories will be denoted by $_{M}\mathcal{M}^{\mathcal{H}}$, $_{M}^{\mathcal{H}}\mathcal{M}$ and $^{\mathcal{H}}\mathcal{M}_{M}$, respectively.

With the previous set-up, where $M$ is a right $\mathcal{H}$-comodule algebra, let us denote by $N=M^{co \ \mathcal{H}_{R}}$. Then we have the following adjunction

\[ \xymatrix{
\mathcal{M}_{N} \ar@<1ex>[rrr]^{-\otimes_{N}M}
& & & \mathcal{M}_{M}^{\mathcal{H}} \ar@<1ex>[lll]^{(-)^{co\ \mathcal{H}_{R}}}
} \]

\noindent The unit of the adjunction is 

\[ V\longrightarrow (V\tens{N}M)^{co\ \mathcal{H}_{R}}\]
\[ v\mapsto v\tens{N} 1\]

\noindent while the counit is

\[ W^{co \ \mathcal{H}_{R}}\tens{N}M\longrightarrow W\]
\[ w\tens{N}m\mapsto w\cdot m. \]

The Hopf algebroid $\mathcal{H}$ is itself a right $\mathcal{H}$-comodule algebra whose $\mathcal{H}_{R}$-coinvariants is the image of $t_{R}$, or equivalently the image of $L\stackrel{s_{L}}{\longrightarrow}H$. The associated induction functor $-\otimes_{L}H:\mathcal{M}_{L}\longrightarrow \mathcal{M}^{H}_{H}$ is an adjoint equivalence.

\subsection{Galois theory of Hopf algebroids}\label{S2.4}

Let $\mathcal{H}=(\mathcal{H}_{L},\mathcal{H}_{R},S)$ be a Hopf algebroid with underlying $k$-module $H$. A $k$-algebra extension $A\subseteq B$ is said to be (\textit{right}) $\mathcal{H}_{R}$-\textit{Galois} if $B$ is a right $\mathcal{H}_{R}$-comodule algebra with $B^{co \ \mathcal{H}_{R}}=A$ and the map

\[ \xymatrix{B\tens{A}B \ar[rr]^-{\mathfrak{gal}_{R}} & & B\tens{R}H } \]
\[ a\tens{A}b\longmapsto ab^{[0]}\tens{R}b^{[1]}\]

\noindent is a bijection. The map $\mathfrak{gal}_{R}$ is called the Galois map associated to the bialgebroid extension $A\subseteq B$. Symmetrically, the extension $A\subseteq B$ is (\textit{right}) $\mathcal{H}_{L}$-\textit{Galois} if $B$ is a right $\mathcal{H}_{L}$-comodule algebra with $B^{co \ \mathcal{H}_{L}}=A$ and the map

\[ \xymatrix{B\tens{A}B \ar[rr]^-{\mathfrak{gal}_{L}} & & B\tens{L}H } \]
\[ a\tens{A}b\longmapsto a_{[0]}b\tens{L}a_{[1]}\]

\noindent is a bijection. We say that a $k$-algebra extension $A\subseteq B$ is $\mathcal{H}$-\textit{Galois} if it is both $\mathcal{H}_{R}$-Galois and $\mathcal{H}_{L}$-Galois. It is not known in general if the bijectivity of $\mathfrak{gal}_{R}$ and $\mathfrak{gal}_{L}$ are equivalent. However, if the antipode $S$ is bijective (which is part of our standing assumption) then $\mathfrak{gal}_{R}$ is bijective if and only if $\mathfrak{gal}_{L}$. To see this, note that $\mathfrak{gal}_{L}=\Phi_{B} \circ \mathfrak{gal}_{R}$ where $\Phi_{B}$ is the map defined in the previous section for $M=B$. Since $S$ is bijective, $\Phi_{B}$ is an isomorphism which gives the desired equivalence of bijectivity of $\mathfrak{gal}_{R}$ and $\mathfrak{gal}_{L}$. Thus, the extension $A\subseteq B$ is $\mathcal{H}$-Galois if it is a bialgebroid Galois extension for any of its constituent bialgebroids.

In the case of Galois extension by Hopf algebras, a class of extensions are of particular interest called cleft extensions. Following \cite{bohm}, we will look what cleft extensions are for Hopf algebroids. But before doing so, let us define what is called a \textit{convolution category}. As before, $R$ and $L$ are $k$-algebras. Let $X$ and $Y$ be $k$-modules such that $X$ has an $R$-coring $(X, \Delta_{R}, \epsilon_{R})$ and an $L$-coring $(X,\Delta_{L},\epsilon_{L})$ structures and $Y$ has an $L\otimes_{k}R$-ring structure with multiplications $\mu_{R}:Y\otimes_{R}Y\longrightarrow Y$ and $\mu_{L}:Y\otimes_{L} Y\longrightarrow Y$. Define the convolution category $Conv(X,Y)$ to be the category with two objects labelled $R$ and $L$. For $I,J\in \left\{R,L\right\}$, a morphism $I\longrightarrow J$ is a $J-I$ bimodule map $X\longrightarrow Y$. For $I,J,K\in \left\{R,L\right\}$ and morphisms $J\stackrel{f}{\longrightarrow}I$ and $K\stackrel{g}{\longrightarrow}J$, we define the composition $f\ast g$ to be the following convolution

\[ f\ast g = \mu_{J} \circ (f\tens{J}g) \circ \Delta_{J}. \]

Now, given a Hopf algebroid $\mathcal{H}=(\mathcal{H}_{L},\mathcal{H}_{R},S)$ and a right $\mathcal{H}$-comodule algebra $B$, $B$ only carries an $R$-ring structure. Since the $k$-module $H$ already has an $R$-coring structure coming from $\mathcal{H}_{R}$ and an $L$-coring structure coming from $\mathcal{H}_{L}$, if the $R$-ring structure of $A$ extends to an $L\otimes_{k}R$-ring structure then we can consider the convolution category $Conv(H,B)$. Since there is no reason for the $A$ to carry a compatible $L$-ring structure, we have to add this to the definition of a cleft extension. Explicitly, an extension $A\subseteq B$, where $A=B^{co \ \mathcal{H}}$ is \textit{cleft} if

\begin{enumerate}
\item[(i)] the $R$-ring structure of $B$ extends to an $L\otimes_{k}R$-ring structure, and
\item[(ii)] there is an invertible morphism $R\stackrel{c}{\longrightarrow}L$ in $Conv(H,B)$ which is a right $\mathcal{H}$-comodule map.
\end{enumerate}

Similar to the case of extensions by Hopf algebras, cleft extensions have Galois-normal basis and crossed product characterizations. Let us state it in the following theorem.

\begin{thm}\label{thm2}
Let $\mathcal{H}=(\mathcal{H}_{L},\mathcal{H}_{R},S)$ be a Hopf algebroid with bijective antipode and let $B$ be a right $\mathcal{H}$-comodule algebra with coinvariants $A$. The following conditions are equivalent:
\begin{enumerate}
\item[(i)] $A\subseteq B$ is a cleft extension.
\item[(ii)] $B\cong A\otimes_{L}H$ as left $A$-modules and right $\mathcal{H}$-comodules (\textit{normal basis property}) and $A\subseteq B$ is $\mathcal{H}$-Galois.
\item[(iii)] For some invertible $A$-valued 2-cocycle $\sigma$ on $\mathcal{H}_{L}$, we have $B\cong A\#_{\sigma}\mathcal{H}_{L}$ as left $A$-modules and as right $\mathcal{H}$-comodule algebras.
\end{enumerate}
\end{thm}

\noindent Let us expound on the last characterization of cleft extensions. Consider a left $L$-bialgebroid $\mathcal{B}=(B,s,t,\Delta,\epsilon)$. Let $(N,\mu,\eta)$ be a $\mathcal{B}$-\textit{measured} $L$-ring, i.e one which is equipped with a $k$-module map $B\otimes_{k}N\stackrel{(\cdot)}{\longrightarrow}N$ satisfying

\begin{enumerate}
\item[(i)] $b\cdot 1_{N}=\eta(\epsilon(b))$,
\item[(ii)] $(t(l)b)\cdot n=(b\cdot n)\eta(l)$ and $(s(l)b)\cdot n=\eta(l)(b\cdot n)$,
\item[(iii)] and $b\cdot (nn^{'})=(b_{(1)}\cdot n)(b_{(2)}\cdot n^{'})$,
\end{enumerate}

\noindent for any $b\in B$, $n,n^{'}\in N$ and $l\in L$. Out of these data, we can construct a two-object category $\mathcal{C}(\mathcal{B},N)$ whose objects are conveniently labelled as $\textit{I}$ and $\textit{II}$. Let us describe the morphism in this category. Consider $B\otimes_{k}B$ as an $L$-bimodule by left multiplication of $s$ and $t$ in the first tensor factor. A map $f\in \ _{L}Hom_{L}(B\otimes_{k}B,N)$ is said to be of \textit{type} $(i,j)$ if it satisfies condition $(i)$ on the first list and condition $(j)$ on the second list below.

\begin{center}
  \begin{tabular}{| l | l |}
    \hline
    $1^{st}$ List & $2^{nd}$ List \\ \hline
		& \\
    (I) $f(a\tens{k}t(l)b)=f(at(l)\tens{k}b)$ & (I) $f(a\tens{k}s(l)b)=f(as(l)\tens{k}b)$ \\
		\hline
		& \\
    (II) $f(a\tens{k}t(l)b)=f(a_{(1)}\tens{k}b)(a_{(2)}\cdot\eta(l))$ & (II) $f(a\tens{k}t(l)b)=(a_{(1)}\cdot\eta(l))f(a_{(2)}\tens{k}b)$ \\
    \hline
  \end{tabular}
\end{center}

\noindent where $a,b\in B$ and $l\in L$. For any $i,j\in \left\{\textit{I},\textit{II}\right\}$, a morphism $i\longrightarrow j$ is a map

\[ f\in \ _{L}Hom_{L}(B\otimes_{k}B,N)\]

\noindent of type $(i,j)$. For any $i,j,l\in \left\{I,II\right\}$, the composition of $i\stackrel{f}{\longrightarrow}j$ and $j\stackrel{g}{\longrightarrow}l$ is the following convolution

\[ \left(f\ast g \right)(a\tens{k}b)=f\left(a_{(1)}\tens{k}a_{(1)}\right)g\left(a_{(2)}\tens{k}b_{(2)}\right). \]

The identity morphism $I\longrightarrow I$ is the map $a\otimes_{k}b\longmapsto (ab)\cdot 1_{N}=\eta(\epsilon(ab))$ and the identity morphism $II\longrightarrow II$ is the map $a\otimes_{k}b\longmapsto a\cdot(b\cdot 1_{N})$.

An $N$-\textit{valued} $2$-\textit{cocycle} on $\mathcal{B}$ is a morphism $I\stackrel{\sigma}{\longrightarrow} II$ in the category $\mathcal{C}(\mathcal{B},N)$ satisfying, for any $a,b,c\in B$, the following conditions.

\begin{enumerate}
\item[(i)] $\sigma(1_{B},b)=\eta(\epsilon(b))=\sigma(b,1_{B})$ \ (\textit{normality}),
\item[(ii)] $(a_{(1)}\cdot \sigma(b_{(1)},c_{(1)}))\sigma(a_{(2)},b_{(2)}c_{(2)})=\sigma(a_{(1)},b_{(1)})\sigma(a_{(2)}b_{(2)},c)$ \ (\textit{cocycle condition}).
\end{enumerate}

\noindent If in addition, we have for any $n\in N$ and $a,b\in B$,

\begin{enumerate}
\item[(iii)] $1_{B}\cdot n=n$ \ (\textit{unitality}),
\item[(iv)] $(a_{(1)}\cdot (b_{(1)}\cdot n))\sigma(a_{(2)},b_{(2)})=\sigma(a_{(1)},b_{(1)})(a_{(2)}b_{(2)}\cdot n)$ \ (\textit{associativity}),
\end{enumerate}

\noindent we call the $\mathcal{B}$-measured $L$-ring $N$ a $\sigma$-\textit{twisted} $\mathcal{B}$-\textit{module}.

For such a left $L$-bialgebroid $\mathcal{B}$ and a $\sigma$-twisted $\mathcal{B}$-module $N$, we can construct the \textit{crossed product} $N\#_{\sigma}\mathcal{B}$ as the $k$-algebra whose underlying $k$-module is $N\otimes_{L}B$ where the left $L$-module structure on $B$ is the one via multiplication of $s$. The multiplication in $N\#_{\sigma}\mathcal{B}$ is defined as

\[ (n\#b)(n^{'}\#b^{'})=n(b_{(1)}\cdot n^{'})\sigma(b_{(2)},b^{'}_{(1)})\#b_{(3)}b^{'}_{(2)}, \hspace{.25in} \ \text{for any} \ n\#b,n^{'}\#b^{'}\in N\#_{\sigma}\mathcal{B}. \]

\noindent This multiplication is associative by conditions $(ii)$ and $(iv)$ and unital by conditions $(i)$ and $(iii)$.

Going back to the characterization of cleft extensions by crossed products, the $2$-cocycle $\sigma$ is invertible in the sense that it is invertible as a morphism in the category $\mathcal{C}(\mathcal{H}_{L},A)$.

%%%%%%%%%%%%%%%%%%%%%%%%%%%%%%%%%%%%%%%%%%%%%%%%%%%%%%%%%%%%%%%%%%%%%%%%%%%%%%%%%%%%%%%%

\section{Noncommutative covering spaces}\label{S3.0}

In the classical case, a covering space is a surjective map $Y\longrightarrow X$ with discrete fibers. In formulating the notion of a noncommutative covering space, discreteness plays a serious obstacle. For one, there is no clear way to translate discreteness for algebras. Fortunately, for our purpose we will only be interested with the analogues of finite coverings. In such cases, discreteness is guaranteed once we go back to the classical case. In the subsequent sections, we will give examples.

%%%%%%%%%%%%%%%%%%%%%%%%%%%%%%%%%%%%%%%%%%%%%%%%%%%%%%%%%%%%%%%%%%%%%%%%%%%%%%%%%%%%%%%%

\subsection{Definitions and properties}\label{S3.1}

The noncommutative analogues of principal bundles are Hopf-Galois extensions. Normal covering spaces are principal bundles in which the gauge group has the discrete topology. At present, it is still unclear how to translate discreteness in the language of algebras. However, if we restrict to \textit{finite} normal coverings then the corresponding Hopf-Galois extension has a finite-dimensional Hopf algebra, the dimension being the same as the degree of the covering. Thus, if we restrict our attention to finite coverings the finiteness assumption for $H$ is sufficient.

To justify our notion of noncommutative covering spaces, let us look at what is happening in the classical case from the algebraic point of view. Let $Y\stackrel{p}{\longrightarrow}X$ be a classical Galois covering space with finite deck transformation group $G$. We assume that $X$ has the suitable connectivity properties, see for example \cite{b001}. Denote by $A$ and $B$ the corresponding algebra of continuous functions on $X$ and $Y$, respectively. The surjection $p$ gives an inclusion $A\subseteq B$.

The covering $Y\stackrel{p}{\longrightarrow}X$ gives a groupoid $\mathcal{G}$ whose set of objects is $X$. For $x,y\in X$, we set $Hom_{\mathcal{G}}(x,y)=\emptyset$ if $x\neq y$. Otherwise, an arrow $x\longrightarrow x$ is a bijection $p^{-1}(x)\stackrel{\gamma^{*}}{\longrightarrow}p^{-1}(x)$ induced by lifting a loop $\gamma$ at $x$ to $Y$. The bijection $\gamma^{*}$ only depends on the homotopy class of $\gamma$. Explicitly we have $\mathcal{G}^{(0)}=X$, $\mathcal{G}^{(1)}$ is the set of induced bijections from homotopy classes of loops in $X$, $s=t:\mathcal{G}^{(1)}\longrightarrow \mathcal{G}^{(0)}$ are the source and target maps giving the base point of the loop inducing the bijection in $\mathcal{G}^{(1)}$, $\mathcal{G}^{(2)}$ is the fiber product of $s$ and $t$ i.e. the composable morphisms on $\mathcal{G}$, $\iota:\mathcal{G}^{(0)}\longrightarrow \mathcal{G}^{(1)}$ the map sending $x$ to the identity map on $p^{-1}(x)$, and finally $inv:\mathcal{G}^{(1)}\longrightarrow \mathcal{G}^{(1)}$ the map that associates to $\gamma^{*}$ the bijection $(\gamma^{-1})^{*}$. These structure maps

\begin{equation} \label{eq:1}
\xymatrix{\mathcal{G}^{(2)} \ar[rd]^-{m} & & &\\
& \mathcal{G}^{(1)} \ar@(ur,ul)[rr]^-{s} \ar@(dr,dl)[rr]_-{t} \ar@(l,d)[]_-{inv} & & \mathcal{G}^{(0)} \ar[ll]_-{\iota} \\}
\end{equation}

\noindent make the following diagrams commute

\begin{equation} \label{eq:2}
\xymatrix{\mathcal{G}^{(3)} \ar[rr]^-{m\times id} \ar[dd]_-{id\times m} & & \mathcal{G}^{(2)} \ar[dd]^-{m} \\
& & & & \\
\mathcal{G}^{(2)} \ar[rr]_-{m} & & \mathcal{G}^{(1)} \\}
\xymatrix{\mathcal{G}^{(1)}\times\mathcal{G}^{(1)} \ar[rr]^-{s\times id}  & & \mathcal{G}^{(0)}\times \mathcal{G}^{(1)} \ar[rr]^-{\iota\times id} & & \mathcal{G}^{(2)} \ar[d]^-{m} \\
\mathcal{G}^{(1)} \ar[u]^-{diag} \ar[d]_-{diag} \ar@{=}[rrrr] & & & & \mathcal{G}^{(1)}\\
\mathcal{G}^{(1)}\times\mathcal{G}^{(1)} \ar[rr]_-{id\times s} & & \mathcal{G}^{(1)}\times\mathcal{G}^{(0)} \ar[rr]_-{id\times \iota} & & \mathcal{G}^{(2)} \ar[u]_-{m} \\}
\end{equation}

\begin{equation} \label{eq:3}
\xymatrix{& & \mathcal{G}^{(0)}  & & & & \mathcal{G}^{(0)} & & \\
\mathcal{G}^{(1)} \ar[rru]^-{s} \ar[rrd]_-{s} & & & & \mathcal{G}^{(1)} \ar[llu]_-{t} \ar[rru]^-{t} \ar[rrd]_-{t} & & & & \mathcal{G}^{(1)} \ar[llu]_-{s} \\
& & \mathcal{G}^{(0)} \ar[rru]_-{\iota} & & & & \mathcal{G}^{(0)} \ar[rru]_-{\iota} & & \\}
\end{equation}

\begin{equation} \label{eq:4}
\xymatrix{& & \mathcal{G}^{(1)}\times\mathcal{G}^{(1)} \ar[r(2.8)]^-{inv \times id} & & & & \mathcal{G}^{(2)} \ar[rd(1.6)]^-{m} & & \\
& & & & & & & & \\
\mathcal{G}^{(1)} \ar[rrrr]^-{s} \ar[ru(1.8)]^-{diag} \ar[rd(1.8)]_-{diag}  & & & & \mathcal{G}^{(0)} \ar[rrrr]^-{\iota}   & & & & \mathcal{G}^{(1)} \\
& & & & & & & & \\
& & \mathcal{G}^{(1)}\times\mathcal{G}^{(1)} \ar[r(2.8)]_-{id\times inv} & & & & \mathcal{G}^{(2)} \ar[ru(1.6)]_-{m} & & \\}
\end{equation}

\noindent where, for $n\geqslant 2$, $\mathcal{G}^{(n)}$ denotes the $n$-fold fiber product of $s$ and $t$.

The above data with corresponding compatibility conditions indeed gives us a (topological) groupoid. We will explore a larger groupoid containing the one we constructed here in section \ref{S3.3} The functor $C(-)$ which associates to a topological space $X$ its algebra of continuous complex-valued functions $C(X)$ is a duality (at least for locally compact Hausdorff topological spaces). Applying this functor to the diagram \ref{eq:1} gives us the following diagram of $A$-rings

\begin{equation} \label{eq:5}
\xymatrix{H\otimes_{A}H & & &\\
& H \ar[lu]_-{\Delta} \ar[rr]^-{\epsilon} \ar@(l,d)[]_-{S} & & A \ar@(ul,ur)[ll]_-{s} \ar@(dl,dr)[ll]^-{t} \\}
\end{equation}

\noindent where $H=C(\mathcal{G}^{(1)})$, $\Delta=C(m)$, $\epsilon=C(\iota)$, $S=C(inv)$, and we denote by the same symbol $s$ and $t$ the induced maps of the groupoid's source and target maps.

The diagrams in \ref{eq:2} dualize to the following diagrams

\begin{equation*}
\xymatrix{ H \ar[rr]^-{\Delta} \ar[dd]_-{\Delta} & & H\tens{A}H \ar[dd]^-{\Delta\otimes_{A} id} \\
& & & & \\
H\tens{A}H \ar[rr]_-{id \otimes_{A} \Delta} & & H\tens{A}H\tens{A}H \\}
\xymatrix{H\otimes H \ar[d]_-{\mu^{'}}  & & A\otimes H \ar[ll]_-{s\otimes id} & & H\tens{A}H \ar[ll]_-{\epsilon\times id} \\
H \ar@{=}[rrrr] & & & & H \ar[u]_-{\Delta} \ar[d]^-{\Delta}  \\
H\otimes H \ar[u]^-{\mu^{'}}  & & H\otimes A \ar[ll]^-{id\otimes s}  & & H\tens{A}H \ar[ll]^-{id\otimes \epsilon}  \\}
\end{equation*}

\noindent which express the coassociativity of $\Delta$ and its counitality with respect to $\epsilon$. Diagram \ref{eq:3} dualizes to the following commutative diagram

\begin{equation} \label{eq:6}
\xymatrix{& & A \ar[lld]_-{s} \ar[rrd]^-{t}  & & & & A \ar[lld]_-{t} \ar[rrd]^-{s} & & \\
H & & & & H \ar[lld]_-{\epsilon}  & & & & H.  \ar[lld]_-{\epsilon}  \\
& & A \ar[llu]^-{s}  & & & & A  \ar[llu]^-{t}  & & \\}
\end{equation}

Note that $C(\mathcal{G}^{(1)}\times \mathcal{G}^{(1)})\cong C(\mathcal{G}^{(1)})\otimes C(\mathcal{G}^{(1)})=H\otimes H$. Thus, diagram \ref{eq:4} dualizes to the outer hexagon of the following diagram

\begin{equation} \label{eq:7}
\xymatrix{& & H\otimes_{A}H \ar[r(3.7)]^-{S \otimes id} \ar[drrr]^-{S \otimes_{A} id} & & & & H\otimes H \ar[rd(1.6)]^-{\mu^{'}} \ar@{->>}[dl] & & \\
& & & & & H\otimes_{A}H \ar[drrr]^-{\mu} & & & \\
H \ar[rrrr]^-{\epsilon} \ar[ru(1.8)]^-{\Delta} \ar[rd(1.8)]_-{\Delta}  & & & & A \ar[r(3.5)]^-{s}   & & & & H \\
& & & & & H\otimes_{A}H \ar[urrr]^-{\mu} & & & \\
& & H\otimes_{A}H \ar[r(3.7)]_-{id\otimes S} \ar[urrr]^-{id \otimes_{A} S} & & & & H\otimes H \ar[ru(1.6)]_-{\mu^{'}} \ar@{->>}[ul] & & \\}
\end{equation}

 The diagonal map $diag$ induces a $\mathbb{C}$-algebra structure $H\otimes H\stackrel{\mu^{'}}{\longrightarrow}H$. Since $A\stackrel{s}{\longrightarrow}H$ is a $\mathbb{C}$-algebra map, by Lemma 2.2 of \cite{bohm} there is a unique $A$-ring structure on $H$ with product $H\otimes_{A}H\stackrel{\mu}{\longrightarrow}H$. The inner commutative hexagon of \ref{eq:7} implies that $S$ is the convolution inverse of $id$ is the convolution category $Conv(H,H)$ defined in section \ref{S2.4}. All these diagrams tells us that $H$ is a Hopf algebroid with coinciding left- and right-bialgebroid structures and antipode $S$. Furthermore, $S$ is bijective.

Going back to the covering $Y\stackrel{p}{\longrightarrow}X$ and its associated groupoid $\mathcal{G}$, there is an action of $\mathcal{G}$ on $Y$ defined as follows

\[ \xymatrix{ \mathcal{G}^{(1)} \prescript{}{s}\times_{p} \ Y \ar[rr]^-{\alpha}& & Y } \]
\[ \left(\phi, y\right) \longmapsto \phi(y). \]

Moreover, $Y/\mathcal{G}\cong X$. Also, the covering $Y\stackrel{p}{\longrightarrow}X$ is Galois if and only if the associated action is Galois, i.e. the following map is a bijection.

\[ \xymatrix{ \mathcal{G}^{(1)} \prescript{}{s}\times_{p} \ Y \ar[rr]^-{\alpha}& & Y \prescript{}{p}\times_{p} \ Y } \]
\[ \left(\phi, y\right) \longmapsto \left(\phi(y),y\right) \]

Dually, this gives a coaction $\xymatrix{ B \ar[rr]^-{\rho} & & B\otimes_{A} H }$ whose coinvariants relative to the unit of $H$ is $A$. Furthermore, the associated map

\[ \xymatrix{ B\otimes_{A}B \ar[rr]^-{\mathfrak{gal}} & & B\otimes_{A}H } \]
\[ a\otimes_{A}b \longmapsto (a\otimes_{A}1)\rho(b) \]

\noindent is a linear bijection. In other words, $A\subseteq B$ is an $H$-Galois extension.

Consider a faithful finite-dimensional representation $\pi$ of $\mathcal{G}$. Explicitly, it is a continuous map $\mathcal{G}\stackrel{\pi}{\longrightarrow}GL(E)$ of groupoids where $E\stackrel{q}{\longrightarrow}X$ is a finite-dimensional vector bundle over $X$ and $GL(E)$ is the associated general linear groupoid. $GL(E)$ has objects points of $X$, there is no arrow between different points of $X$ and for $x\in X$, an arrow $x\longrightarrow x$ is a linear automorphism of $E_{x}$. It is clear that the $GL(E)$ acts continuously on $E$. Construct the topological space $W=W(Y,\pi)$ as the space $\left(Y\prescript{}{p}\times_{q} E\right)/\sim$ where $(y,e)\sim(g\cdot y,\pi(g)e)$ for all $y\in Y$, $e\in E$ and $g\in \mathcal{G}^{(1)}$. Here, $(\cdot)$ refers to the action $\alpha$ of the groupoid $\mathcal{G}$ to $Y$. Since the fibers of $p$ are orbits of the $\mathcal{G}$-action on $Y$, there is a well-defined projection $W\stackrel{r}{\longrightarrow}X$ sending $(y,v)\mapsto p(y)$ making $W$ a finite-dimensional vector bundle over $X$. $W=W(Y,\pi)$ is called the associated vector bundle to $Y\stackrel{p}{\longrightarrow}X$ and the representation $\mathcal{G}\stackrel{\pi}{\longrightarrow}GL(E)$.

As before, projection $Y\stackrel{p}{\longrightarrow}X$ gives an algebra inclusion $A\subseteq B$ which makes $B$ into an $A$-module. Also, the global sections $\Gamma(X,W)$ is also a module over $A$ which is finitely generated and projective by the Serre-Swan theorem. Note that by the construction of the associated bundle, $\Gamma(X,W)$ and $B$ are isomorphic as $A$-modules. Thus, $B$ is a finitely-generated projective $A$-module.

Using the arguments above, we present the following definition of a noncommutative covering space.

\begin{defn}\label{D2}
Let $A$ be an algebra over a commutative unital ring $k$. A \textit{(finite, Galois) noncommutative covering} of $A$ is a pair $(B,\mathcal{H})$ where:
\begin{enumerate}
\item[(i)] $\mathcal{H}$ is a finitely generated projective Hopf algebroid with bijective antipode $S$ over $A^{'}$ where $A^{'}$ is a subalgebra of $A$,
\item[(ii)] $A\subseteq B$ is a right $\mathcal{H}$-Galois extension,
\item[(iii)] $B$ is a finitely-generated projective $A$-module via the inclusion $A\subseteq B$.
\end{enumerate}
\noindent If furthermore, $B$ only has $0$ and $1$ as idempotents then the covering $(B,\mathcal{H})$ is said to be \textit{connected}. If $A^{'}=A$, we will call $(B,\mathcal{H})$ a \textit{local} noncommutative covering of $A$. Otherwise, it is called \textit{stratified} with stratification datum $A^{'}\subseteq A$. A local noncommutative covering $(B,\mathcal{H})$ of $A$ is called \textit{uniform} if $A^{'}\cong k$ and $\mathcal{H}$ is a Hopf algebra.
\end{defn}

\begin{rem}
\begin{enumerate}
\item[]

\item[(i)] Since the present work is mainly concerned with noncommutative analogues of (finite, Galois) connected covering spaces, we will simply refer to a (finite, Galois) noncommutative covering as a covering and reserve the name classical covering for classical ones.

\item[(ii)] It is important to note that the Hopf algebroid $\mathcal{H}$ carries several module structure over $A^{'}$ using the sources and targets. However, bijectivity of $S$ implies that finitely-generated projectivity over $A^{'}$ are all equivalent for the module structures induced by multiplication of the sources. Same is true for the modules structures induced by multiplication by the targets. See proposition 4.5 of \cite{bohm}. Now, by definition of the Takeuchi product, these equivalences go between module structures induced from multiplication by a source map and a target map. This makes part $(i)$ of definition \ref{D2} well-defined.

\item[(iii)] In a covering $(B,\mathcal{H})$ of $A$, we call $\mathcal{H}$ the associated \textit{quantum symmetry} or just symmetry for brevity, of the covering. This corresponds to the deck transformation group in the classical set up. Note that for a classical covering space $\xymatrix{ Y \ar@{->>}[r] & X}$ the deck transformation group is completely determined as $G=Aut_{X}(Y)$. In the general case, there might be different quantum symmetries $\mathcal{H}_{1}$ and $\mathcal{H}_{2}$ making an extension $A\subseteq B$ Hopf-Galois and hence $(B,\mathcal{H}_{1})$ and $(B,\mathcal{H}_{2})$ are potentially different coverings. See \cite{m003} for an example of an extension $A\subset B$ which is Galois for different quantum symmetries.

\item[(iv)] The motivation we outlined in this section suggests that in a noncommutative covering $(B,\mathcal{H})$ of $A$, the Hopf algebroid $\mathcal{H}$ is over $A$. However, as we will see in section \ref{S5.2} there are some interesting structures where we need to consider Hopf algebroids over any subalgebra of $A$. In section \ref{S5.3} we will look at the contrast between local and stratified coverings.
\end{enumerate}
\end{rem}

\

The analogues of general finite coverings (possibly non-Galois) are those extensions whose associated Hopf-Galois map is surjective but not necessarily injective. This is justified by the following observation. The deck transformation group of a covering always act freely. But the covering is normal precisely when, aside from being free, the action is transitive. So to get the analogue of general coverings we simply drop the condition that the action is transitive. But transitivity translates to surjectivity of the associated Galois map. The functor $C(-)$, the one that associates to a space $X$ its algebra of functions $C(X)$, is contravariant. Thus, surjectivity of the associated Galois map is equivalent to the injectivity of the associated Hopf-Galois map.

\subsection{Equivalences of coverings}\label{S3.2}

In this section, we will look at two notions of equivalence of coverings. We will focus our attention to local coverings i.e., those coverings $(B,\mathcal{H})$ of $A$ whose quantum symmetry $\mathcal{H}$ is over $A$. We will briefly discuss how these equivalences work when we are dealing with stratified coverings. The first notion, which we call \textit{topological equivalence} is the direct dualization of equivalence of coverings in the classical sense. The second one called \textit{Morita equivalence} is a prominent equivalence in noncommutative geometry.

\begin{defn}\label{D3}
Let $(B,\mathcal{H})$ and $(B^{'},\mathcal{H}^{'})$ be coverings of a noncommutative space $A$. We say that $(B^{'},\mathcal{H}^{'})$ is an \textit{intermediate covering} of $(B,\mathcal{H})$ if there is an intermediate inclusion $A\subseteq B^{'} \subseteq B$ and a monomorphism $\xymatrix{\mathcal{H}^{'} \ar@{^{(}->}[r] & \mathcal{H} }$ of Hopf algebroids such that the restriction of the coaction of $\mathcal{H}$ on $B^{'}\subseteq B$ gives the coaction of $\mathcal{H}^{'}$ on $B^{'}$. Two coverings are \textit{topologically equivalent} if they are intermediate coverings of each other.
\end{defn}

\begin{rem}
In the classical case, a covering $Y\twoheadrightarrow X$ is an intermediate covering of $Z\twoheadrightarrow X$ there is a (continuous) surjection $Z\twoheadrightarrow Y$ and a group epimorphism $Aut_{X}(Z)\twoheadrightarrow Aut_{X}(Y)$.
\end{rem}

\

Let $(B,\mathcal{H})$ and $(B^{'},\mathcal{H}^{'})$ be topologically equivalent coverings of a noncommutative space $A$. Immediately, we see that $B\cong B^{'}$ as $A$-rings. By definition, there are injective maps of Hopf algebroids $\xymatrix{\mathcal{H} \ar@{^{(}->}[r]^-{p} & \mathcal{H}^{'}}$ and $\xymatrix{\mathcal{H}^{'} \ar@{^{(}->}[r]^-{q} & \mathcal{H}}$. Using the map $p$, $\mathcal{H}$ becomes a right $\mathcal{H}^{'}$-comodule algebra via the coactions $\rho_{L}$ and $\rho_{R}$ defined by the composition

\[ \xymatrix{\mathcal{H} \ar[rr]^-{\Delta_{L},\Delta_{R}} & & \mathcal{H}\otimes_{A}\mathcal{H} \ar[r] & \mathcal{H}\otimes_{A}\mathcal{H}^{'} } \]

\noindent whose coassiociativity follows from the commutativity of the following diagram

\[ \xymatrix{ H \ar[rr]^-{\Delta_{L},\Delta_{R}} \ar[dd]_-{\Delta_{R},\Delta_{L}} & & H\otimes_{A}H \ar[rr]^-{id\otimes p} \ar[dd]^-{\Delta_{R},\Delta_{L}}  & & H\otimes_{A}H^{'} \ar[dd]^-{\Delta^{'}_{R},\Delta^{'}_{L}} \\
& & & & \\
H\otimes_{A}H \ar[rr]^-{\left(\Delta_{L},\Delta_{R}\right)\otimes id} \ar[dd]_-{id\otimes p} & & H\otimes_{A}H\otimes_{A}H \ar[rr]^-{id\otimes p\otimes p} \ar[dd]^-{id\otimes id\otimes p} & & H\otimes_{A}H^{'}\otimes_{A}H^{'} \ar[dd]^-{id\otimes id\otimes id} \\
& & & & \\
H\otimes_{A}H^{'} \ar[rr]_-{\left(\Delta_{L},\Delta_{R}\right)\otimes id} & & H\otimes_{A}H\otimes_{A}H^{'} \ar[rr]_-{id\otimes p\otimes id} & & H\otimes_{A}H^{'}\otimes_{A}H^{'}. \\ } \]

\noindent Let us determine the coinvariants of $\mathcal{H}$ under this coaction of $\mathcal{H}^{'}_{R}$. An element $a\in H$ is coinvariant if $\rho_{R}(a)=a\otimes_{A}1$. This means that there exist $h\in H$ such that $p(h)=1$ and $\Delta_{R}(a)=a\otimes_{A}h$. Injectivity of $p$ implies that $h=1$ and hence, $\Delta_{R}(a)=\Delta_{L}(a)=a\otimes_{A}1$. Thus, the coinvariants of $\rho_{L}$ and $\rho_{R}$ coincide with the coinvariants of the regular comodule structure of $\mathcal{H}$ which is $A$ itself.

\noindent Meanwhile, using the map $q$ we can equip $\mathcal{H}$ a structure of a right $\mathcal{H}^{'}$-module via

\[ \xymatrix{H\otimes_{A} H^{'} \ar[r]^-{id\otimes q} & H\otimes_{A}H \ar[r]^-{m} & H  } \]

\noindent which makes $\mathcal{H}$ a right-right $\mathcal{H}^{'}$-Hopf module. The counit of the adjoint equivalence $(-\otimes_{A}H^{'}) \dashv (-)^{co\ \mathcal{H}^{'}_{R}}$ provides an isomorphism

\[ \xymatrix{H^{'}\cong A\otimes_{A}H^{'}\cong (H)^{co \ \mathcal{H}^{'}_{R}}\otimes_{A}H^{'} \ar[rr]^-{counit} & & H} \]

\noindent of right-right $\mathcal{H}^{'}$-Hopf modules. Reversing the roles of $\mathcal{H}$ and $\mathcal{H}^{'}$ in the above computation shows that $H$ and $H^{'}$ are also isomorphic as right-right $\mathcal{H}$-Hopf modules. This is enough to conclude that $\mathcal{H}$ and $\mathcal{H}^{'}$ are isomorphic Hopf algebroids. This gives the following proposition.

\begin{prop}\label{prop1}
Let $(B,\mathcal{H})$ and $(B^{'},\mathcal{H}^{'})$ be coverings of a noncommutative space $A$. Then $(B,\mathcal{H})$ and $(B^{'},\mathcal{H}^{'})$ are topologically equivalent if and only if $B\cong B^{'}$ as $A$-rings and $\mathcal{H}\cong \mathcal{H}^{'}$ as Hopf algebroids.
\end{prop}

\

Before stating the second equivalence, let us recall what a Hopf bimodule is. Given Hopf algebroids $\mathcal{H}$ and $\mathcal{H}^{'}$ over $R$, an $(\mathcal{H},\mathcal{H}^{'})$-\textit{Hopf bimodule} $M$ is an $(\mathcal{H},\mathcal{H}^{'})$-bimodule and an $(\mathcal{H},\mathcal{H}^{'})$-bicomodule such that the bicomodule structure maps $M\stackrel{\rho_{L},\rho_{R}}{\longrightarrow}H\otimes_{R}M$ and $M\stackrel{\rho^{'}_{L},\rho^{'}_{R}}{\longrightarrow}M\otimes_{R}H^{'}$ are $(\mathcal{H},\mathcal{H}^{'})$-bimodule maps.

\begin{defn} \label{D4}
Two coverings $(B,\mathcal{H})$ and $(B^{'},\mathcal{H}^{'})$ of a noncommutative space $A$ are \textit{Morita equivalent} if the following conditions are satisfied:

\begin{enumerate}
\item[(i)] There exist a $(B,B^{'})$-bimodule $\mathcal{X}$ and a $(B^{'},B)$-bimodule $\mathcal{Y}$ such that

\[ \mathcal{X} \tens{B^{'}} \mathcal{Y} \cong B, \hspace{.5in} \mathcal{Y} \tens{B} \mathcal{X} \cong B^{'} \]

\noindent as $B$-bimodules and $B^{'}$-bimodules, respectively.

\item[(ii)] $ B^{'} \tens{B^{'}} \mathcal{Y} \cong B $ and $ B \tens{B} \mathcal{X} \cong B^{'} $ as $(A,B)$-bimodules and $(A,B^{'})$-bimodules, resp.

\item[(iii)] There exist an $(\mathcal{H},\mathcal{H}^{'})$-Hopf bimodule $\mathcal{U}$ and an $(\mathcal{H}^{'},\mathcal{H})$-Hopf bimodule $\mathcal{V}$ such that

\[ \mathcal{U} \tens{\mathcal{H}^{'}} \mathcal{V} \cong \mathcal{H}, \hspace{.5in} \mathcal{V} \tens{\mathcal{H}} \mathcal{U} \cong \mathcal{H}^{'} \]

\noindent as $\mathcal{H}$-Hopf bimodules and $\mathcal{H}^{'}$-Hopf bimodules, respectively.

\item[(iv)] $ \mathcal{H} \tens{\mathcal{H}} \mathcal{U}  \cong \mathcal{H}^{'} $ and $ \mathcal{H}^{'} \tens{\mathcal{H}^{'}} \mathcal{V}  \cong \mathcal{H} $ as $(A^{e},\mathcal{H}^{'})$-Hopf bimodules and $(A^{e},\mathcal{H})$-Hopf bimodules, respectively. Here, $A^{e}$ is the canonical Hopf algebroid over $A$.

\end{enumerate}
\end{defn}

\

\begin{rem}
\begin{enumerate}
\item[]

\item[(a)] Topologically equivalent coverings $(B,\mathcal{H})$ and $(B^{'},\mathcal{H}^{'})$ are Morita equivalent. $B$ and $B^{'}$ provide the bimodules asked in (i) and (ii) while $\mathcal{H}$ and $\mathcal{H}^{'}$ provide the Hopf bimodules required in (iii) and (iv).

\item[(b)] Requiring both the isomorphism in part (ii) of the definition \ref{D4} is redundant as one implies the other one. To see this, assume $B^{'} \tens{B^{'}} \mathcal{Y} \cong B$ as $(A,B)$-bimodules. Then

\[ B \tens{B} \mathcal{X} \cong B^{'} \tens{B^{'}} \mathcal{Y} \tens{B} \mathcal{X} \cong B^{'}\tens{B^{'}}B^{'} \cong B^{'} \]

\noindent as $(A,B^{'})$-bimodules. Similarly, each Hopf bimodule isomorphism in part (iv) implies the other one.

\item[(c)] Morita equivalences of coverings coincide with isomorphisms in a suitable category. Denote by $Cov(A)$ the category whose objects are coverings of a noncommutative space $A$. A morphism $(B,\mathcal{H})\longrightarrow(B^{'},\mathcal{H}^{'})$ is a pair $(\mathcal{X},\mathcal{U})$ of a $(B,B^{'})$-bimodule $\mathcal{X}$ and an $(\mathcal{H},\mathcal{H}^{'})$-Hopf bimodule $\mathcal{U}$ such that

\[ B \tens{B} \mathcal{X} \cong B^{'}, \hspace{.5in} \mathcal{H} \tens{\mathcal{H}} \mathcal{U}  \cong \mathcal{H}^{'} \]

\noindent as $(A,B^{'})$-bimodules and as $(A^{e},\mathcal{H}^{'})$-Hopf bimodules, respectively. The composition rule given by

\[ \xymatrix{ & & (B^{'},\mathcal{H}^{'}) \ar[rrdd]^-{(\mathcal{Y},\ \mathcal{V})} & & \\
& & & & \\
(B,\mathcal{H}) \ar[rruu]^-{(\mathcal{X},\mathcal{U})} \ar[rrrr]_-{\left(\mathcal{X}\tens{B^{'}}\mathcal{Y},\ \mathcal{U}\tens{\mathcal{H}^{'}}\mathcal{V}\right)} & &  & & (B^{''},\mathcal{H}^{''}). } \]

\noindent The identity morphism of the object $(B,\mathcal{H})$ is the pair $(B,\mathcal{H})$ itself. It is now immediate to see that the isomorphisms in $Cov(A)$ are precisely the Morita equivalences. We will call such invertible arrow $(B,\mathcal{H})$ a \textit{Morita equivalence bimodule}.

\item[(d)] Recall that two noncommutative spaces $A$ and $A^{'}$ are \textit{Morita equivalent} if there exist an $(A,A^{'})$-bimodule $\mathcal{P}$ and an $(A^{'},A)$-bimodule $\mathcal{Q}$ such that

\[ \mathcal{P} \tens{A^{'}} \mathcal{Q} \cong A, \hspace{.5in} \mathcal{Q} \tens{A} \mathcal{P} \cong A^{'} \]

\noindent as $A$-bimodules and $A^{'}$-bimodules, respectively. Notice that Morita equivalence of coverings puts together Morita equivalence of the extension algebras (part i) and Hopf-fitted notion of Morita equivalence for the associated symmetries (part iii). Parts ii and iv of the definition take care of the fact that these equivalences are in the \textit{over-category} of noncommutative spaces over $A$.

\item[(e)] In light of remark (d), we will say that two Hopf algebroids $\mathcal{H}^{'}$ and $\mathcal{H}^{'}$ over $A$ are \textit{Morita equivalent} if there exists an $(\mathcal{H},\mathcal{H}^{'})$-Hopf bimodule $\mathcal{U}$ and an $(\mathcal{H}^{'},\mathcal{H})$-Hopf bimodule $\mathcal{V}$ satisfying conditions (iii) and (iv) of definition \ref{D4}.
\end{enumerate}
\end{rem}

\

Let us end this section by stating a result which verifies that coverings are noncommutative geometric invariants.

\begin{prop}\label{prop2}
Let $A$ and $A^{'}$ be Morita equivalent noncommutative spaces. Then $Cov(A)$ and $Cov(A^{'})$ are adjoint equivalent categories.
\end{prop}

\begin{prf}
Consider Morita equivalent noncommutative spaces $A$ and $A^{'}$. Let $(B,\mathcal{H})$ be a covering of $A$. We will construct a covering of $A^{'}$ whose Morita equivalence class is uniquely determined by the Morita equivalence class of $(B,\mathcal{H})$.

By assumption, there is an $(A,A^{'})$-bimodule $\mathcal{P}$ and an $(A^{'},A)$-bimodule $\mathcal{Q}$ such that

\[ \mathcal{P}\tens{A^{'}}\mathcal{Q}\cong A, \hspace{.5in} \mathcal{Q}\tens{A}\mathcal{P} \cong A^{'}. \]

\noindent We claim that $\left(B^{'},\mathcal{H}^{'}\right)=\left( \mathcal{Q}\otimes_{A}B\otimes_{A}\mathcal{P},\mathcal{Q}\otimes_{A}\mathcal{H}\otimes_{A}\mathcal{P} \right)$ is a covering of $A^{'}$. By $\mathcal{Q}\otimes_{A}\mathcal{H}\otimes_{A}\mathcal{P}$ we mean the Hopf algebroid with constituent left- and right-bialgebroids $H_{L}^{'}=\mathcal{Q}\otimes_{A}H_{L}\otimes_{A}\mathcal{P}$ and $H_{R}^{'}=\mathcal{Q}\otimes_{A}H_{R}\otimes_{A}\mathcal{P}$, respectively.

First, let us show that $B^{'}$ is an $A^{'}$-ring. The $A$-bimodule structure maps

\[ \xymatrix{B\otimes_{A}B \ar[r]^-{\mu} & B}, \hspace{.5in} \xymatrix{A \ar[r]^-{\eta} & B  } \]

\noindent of $B$ as an $A$-ring induce the following $A^{'}$-bimodule maps

\[ \xymatrix{B^{'}\otimes_{A^{'}}B^{'}\cong\mathcal{Q}\otimes_{A}B\otimes_{A}B\otimes_{A}\mathcal{P} \ar[rrr]^-{\mathcal{Q}\otimes_{A}\mu\otimes_{A}\mathcal{P}} & & & \mathcal{Q}\otimes_{A}B\otimes_{A}\mathcal{P}\cong B^{'}} \]

\[ \xymatrix{A^{'}\cong\mathcal{Q}\otimes_{A}A\otimes_{A}\mathcal{P} \ar[rrr]^-{\mathcal{Q}\otimes_{A}\eta\otimes_{A}\mathcal{P}} & & & \mathcal{Q}\otimes_{A}B\otimes_{A}\mathcal{P}\cong B^{'} } \]

\noindent which satisfy the associativity and the unitality diagrams. These maps make $B^{'}$ into an $A^{'}$-ring. Note that the above argument is just the application of the functors $\mathcal{Q}\otimes_{A}-$ and $-\otimes_{A}\mathcal{P}$ which are both equivalence by the Morita property. Thus, they preserve diagrams. We will make use of this argument in the rest of the proof.

Now, it is easy to see that $\mathcal{H}^{'}$ is a Hopf algebroid over $A^{'}$ since the maps and diagrams that define the Hopf algebroid structure on $\mathcal{H}$ all live in the category of $A$-bimodules. Applying the functors $\mathcal{Q}\otimes_{A}-$ and $-\otimes_{A}\mathcal{P}$ give the structure maps for $\mathcal{H}^{'}$ which satisfy the relevant diagrams. For the same reason, $B^{'}$ carries an $\mathcal{H}^{'}$-comodule structure via

\[ \xymatrix{ B^{'}\cong \mathcal{Q}\otimes_{A}B\otimes_{A}\mathcal{P} \ar[rrr]^-{\mathcal{Q}\otimes_{A}\rho_{R}\otimes_{A}\mathcal{P}}_-{\mathcal{Q}\otimes_{A}\rho_{L}\otimes_{A}\mathcal{P}} & & & \left(\mathcal{Q}\otimes_{A}B\otimes_{A}\mathcal{P}\right)\tens{A^{'}}\left(\mathcal{Q}\otimes_{A}\mathcal{H}\otimes_{A}\mathcal{P}\right)\cong B^{'}\tens{A^{'}}\mathcal{H}^{'} }. \]

\noindent The $H_{R}^{'}$-coinvariants $\left(B^{'}\right)^{co\ H_{R}^{'}}$ of this comodule structure is the equalizer of $\rho_{R}^{'}$ and $-\otimes_{A^{'}}H_{R}^{'}$, i.e.

\[ \xymatrix{ \left(B^{'}\right)^{co\ H_{R}^{'}} \ar[rr] & & B^{'} \ar@<1ex>[rr]^-{\rho_{L}^{'}} \ar@<-1ex>[rr]_-{\rho_{R}^{'}} & & B^{'}\tens{A^{'}}\mathcal{H}^{'}   }. \]

\noindent This diagram is the image of the equalizer diagram defining $B^{co\ H_{R}}$ after applying $\mathcal{Q}\otimes_{A}-$ and $-\otimes_{A}\mathcal{P}$. Thus, $\left(B^{'}\right)^{co\ H_{R}^{'}}\cong \mathcal{Q}\otimes_{A}B^{co\ H_{R}}\otimes_{A}\mathcal{P}\cong\mathcal{Q}\otimes_{A}A\otimes_{A}\mathcal{P}\cong A^{'}$.

Finally, finitely-generated projectivity of $B^{'}$ and $\mathcal{H}^{'}$ is equivalent to finitely-generated projectivity of $B$ and $\mathcal{H}$. This proves our claim.

Now, any covering of $A^{'}$ Morita equivalent to $(B^{'},\mathcal{H}^{'})$ is of the form

\[\left( B^{'}\tens{B^{'}}\mathcal{X} , \mathcal{H}^{'}\tens{\mathcal{U}}\mathcal{H}^{'} \right) \]

\noindent for some Morita equivalence bimodule $(\mathcal{X},\mathcal{U})$. Again, by $\mathcal{H}^{'}\otimes_{\mathcal{H}^{'}}\mathcal{U}$ we mean the Hopf algebroid whose consituent bialgebroids are the images of that of $\mathcal{H}^{'}$ under the functor $-\otimes_{\mathcal{H}^{'}}\mathcal{U}$. Invertibility of $(\mathcal{X},\mathcal{U})$ implies that there exist a Morita equivalence bimodule $(\mathcal{Y},\ \mathcal{V})$ such that applying the functor $\mathcal{F}=\mathcal{P}\otimes_{A^{'}}\left(\mathcal{Y}\otimes-\right)\otimes_{A^{'}}\mathcal{Q}$ to $B^{'}\otimes_{B^{'}}\mathcal{X}$ and the functor $\mathcal{G}=\mathcal{P}\otimes_{A^{'}}\left(\mathcal{V}\otimes-\right)\otimes_{A^{'}}\mathcal{Q}$ to $\mathcal{V}\otimes_{\mathcal{H}^{'}}\mathcal{H}^{'}$ yields a covering of $A$ Morita equivalent to $(B,\mathcal{H})$. This proves the proposition. $\blacksquare$
\end{prf}

%%%%%%%%%%%%%%%%%%%%%%%%%%%%%%%%%%%%%%%%%%%%%%%%%%%%%%%%%%%%%%%%%%%%%%%%%%%%%%%%%%%%%%%%

\subsection{Composition of coverings}\label{S3.3}

The following commutative diagram of classical covering spaces

\begin{equation} \label{eq:d1}
\xymatrix{ & & Z \ar@{->>}[lld]_-{q} \ar@{->>}[lddd]^-{p} \\
Y \ar@{->>}[rdd]_-{r} & & \\
& & \\
& X & \\ }
\end{equation}

\noindent has three different interpretations which individually has corresponding interpretations in the present set-up. The first one, by viewing $Y\stackrel{r}{\longrightarrow}X$ as an intermediate covering of $Z\stackrel{p}{\longrightarrow}X$, one gets the notion of intermediate covering we defined in section \ref{S3.2}. The second one, by viewing $Z\stackrel{q}{\longrightarrow}Y$ as an arrow from $Z\stackrel{p}{\longrightarrow}X$ to $Y\stackrel{r}{\longrightarrow}X$ in the category of classical coverings of $X$, one is lead to the notion we defined in section \ref{S3.2}. The third one, which is the main subject of this section is the analogue of the fact that $Z\stackrel{p}{\longrightarrow}X$ is the composition of the coverings $Z\stackrel{q}{\longrightarrow}Y$ and $Y\stackrel{r}{\longrightarrow}X$.

Let $G=Aut(Z\stackrel{p}{\longrightarrow}X)$, $H=Aut(Z\stackrel{q}{\longrightarrow}Y)$ and $K=Aut(Y\stackrel{r}{\longrightarrow}X)$ be the automorphism groups of the indicated classical covering maps in the appropriate over-category. Then, we have the following proposition.

\begin{prop}(Exact fitting for classical covering spaces.) \label{prop3}
Using the notation of this section, the commutativity of diagram \ref{eq:d1} implies the exactness of the following sequence

\[ \xymatrix{ 0 \ar[r] & H \ar[r] & G \ar[r] & K \ar[r] & 0. } \]

\noindent Moreover, any extension $G$ of $K$ by $H$ gives a commutative diagram as \ref{eq:d1}.
\end{prop}

\begin{prf} Let us outline a proof of this classical fact. Assume \ref{eq:d1} commutes. Let $\gamma\in H$. Then commutativity of the smaller triangles in the following diagram

\[ \xymatrix{ Z \ar[dd]_-{\gamma} \ar@{->>}[rd]^-{q} \ar@/^/@{->>}[rrrrd]^{p}  & & \\
 & Y \ar@{->>}[rrr]_-{r} & & & X \\
Z \ar@{->>}[ru]_-{q} \ar@/_/@{->>}[rrrru]_{p} \\} \]

\noindent implies that $\gamma\in G$. It is immediate to see that this defines an injection $H\longrightarrow G$. Let us define a map $\chi:G\longrightarrow K$ as follows: for $g\in G$, let $\chi(g):Y\longrightarrow Y$, $y\mapsto qgq^{-1}(y)$. The map $\chi(g)$ is independent of any pre-image of $y$ under $q$. Also, for any $y\in Y$, we have

\[ r\chi(g)(y)=rqg q^{-1}(y)=pg q^{-1}(y) = pq^{-1}(y) =r(y) \]

\noindent which implies that $\chi(g)\in K$. To see that $\chi$ is surjective, for any $\gamma\in K$ let $\gamma^{*}$ be the pullback of $\gamma$ along $q$. Then $\gamma^{*}\in G$ and $\chi(\gamma^{*})=\gamma$. Finally, let us show that $H=ker \ \chi$. Let $g\in G$ such that $\chi(g)=id$. Then we have

\[ \xymatrix{ Z \ar[rr]^-{g} \ar[dd]_-{q} & & Z \ar[dd]^-{q} \\
& & \\
Y \ar[rr]_-{\chi(g)=id} & & Y \\ } \]

\noindent which immediately implies that $g\in H$.

Let $G$ be an extension of $K$ by $H$. Consider the classifying space $BG$ of $G$. By definition, there is a space $EG$ and a surjective map $EG\stackrel{\tilde{p}}{\longrightarrow} BG$ which is a $G$-principal bundle. In other words, $\tilde{p}$ is a classical Galois covering map with $G$ as its deck transformation group. Dividing $EG$ by the restricted action of $H$ gives a diagram

\[ \xymatrix{ & & Z \ar[rrd] \ar@{->>}[lld]_-{q} \ar@{->>}[lddd]|(0.56)\hole^(0.3){p} & & \\
Y \ar[rrd] \ar@{->>}[rdd]_-{r} & & & & EG \ar@{->>}[lld]_-{\tilde{q}} \ar@{->>}[lddd]^-{\tilde{p}} \\
& & EG/H \ar@{->>}[rdd]_-{\tilde{r}} & & \\
& X \ar[rrd] & & & \\
& & & BG & \\ } \]

\noindent of classical covering spaces with $\tilde{q}$ the canonical surjection and $\tilde{r}$ a covering map with $K$ as its deck transformation group. Pulling-back $\tilde{q}$ and $\tilde{r}$ along the classifying map $X\longrightarrow BG$ gives such a commutative diagram as \ref{eq:d1}. This proves the above proposition. $\blacksquare$
\end{prf}

Let us formulate the above proposition in terms of groupoids. To any classical covering $Y\stackrel{p}{\twoheadrightarrow}X$, we can associate a topological groupoid $\mathscr{G}$ as follows. We set $\mathscr{G}^{(0)}=X$, the space of objects. For $x,y\in X$, a morphism $x\rightarrow y$ is a bijection $p^{-1}(x)\longrightarrow p^{-1}(y)$ induced by lifting to $Y$ a continuous path from $x$ to $y$ in $X$. Using this groupoid, we get another topological groupoid $\mathscr{G}^{'}$ by setting $(\mathscr{G}^{'})^{(0)}=\mathscr{G}^{(0)}$ and $(\mathscr{G}^{'})^{(1)}\subseteq \mathscr{G}^{(1)}$ given as

\[ (\mathscr{G}^{'})^{(1)}(x,y)= 
\begin{cases}
    \mathscr{G}^{(1)}(x,x), & \text{if } x=y\\
		& \\
    \emptyset,              & \text{otherwise.}
\end{cases}
\]

\noindent We will call this associated groupoid the \textit{deck transformation groupoid} of $Y\stackrel{p}{\twoheadrightarrow}X$. The covering $Y\stackrel{p}{\twoheadrightarrow}X$ is Galois if and only if the associated groupoid action of $\mathscr{G}^{'}$ on $Y\stackrel{p}{\twoheadrightarrow}X$ is Galois, i.e. the following map is a bijection.

\[ \xymatrix@R=0.2mm{\mathscr{G}^{'} \prescript{}{s}{\times}_{p}Y \ar[rr] && Y \times_{X} Y \\
(g,y) \ar@{|->}[rr] && (gy,y)} \]

A partial converse is true. Let $\mathscr{G}$ be a locally finite, connected groupoid over $X$ where the subspace topology on each hom-set is discrete. Then $Y=\coprod\limits_{x\in X}\mathscr{G}^{(1)}(x,x)\subseteq \mathscr{G}^{(1)}$, equipped with subspace topology, is a principal $G$-bundle with $G=\mathscr{G}^{(1)}(x_{0},x_{0})$ for any fixed $x_{0}\in X$. The bundle map is given by the restriction of the source map on $Y$. Discreteness of the hom-sets imply that this principal bundle is a covering with deck transformation group $G$. This gives us an isomorphism between the category of finite, connected classical coverings of a (pointed) space $X$ and the category of locally finite, connected groupoids over $X$ with discrete hom-sets.

Now, consider locally finite, connected groupoids $\mathscr{G}$ and $\mathscr{K}$ over $X$ with discrete hom-sets. Let $\mathscr{G}\stackrel{\psi}{\longrightarrow}\mathscr{K}$ be a groupoid homomorphism which is identity on objects and surjective on hom-sets. The construction we just illustrated is clearly functorial. Denote by $Z\stackrel{p}{\longrightarrow}X$ and $Y\stackrel{q}{\longrightarrow}X$ the associated covering spaces to $\mathscr{G}$ and $\mathscr{K}$, respectively. The groupoid map $\psi$ then induces a map of classical covering spaces $Z\stackrel{\psi^{*}}{\longrightarrow}Y$. It is easy to see that $\psi^{*}$ is itself a covering. The groupoid $\mathscr{H}$ associated to $\psi^{*}$ is given as $\mathscr{H}^{(0)}=Y$ and $\mathscr{H}^{(1)}(k_{1},k_{2})=\psi^{-1}(k_{2}^{-1}k_{1})$. This gives an exact sequence of groupoids

\[ \xymatrix@R=.5mm{ \mathscr{H} \ar@{^(->}[rr] && \mathscr{G} \ar@{->>}[rr]^-{\psi} && \mathscr{K} \\
&&  && \\
&&  && \\
Y \ar@{->>}[rr] && X \ar[rr]^-{id} && X \\
\mathscr{H}^{(1)} \ar@{^(->}[rr] && \mathscr{G}^{(1)} \ar@{->>}[rr] && \mathscr{K}^{(1)}. } \]

The proposition and the construction above motivate the following definition. Let us diagrammatically write $P\stackrel{\mathcal{S}}{\Longrightarrow}Q$ when $(Q,\mathcal{S})$ is a local covering of $P$.

\begin{defn}\label{D5} Consider inclusions of $k$-algebras $A\subseteq B^{1} \subseteq B^{2}$, Hopf algebroids $\mathcal{H}$ and $\mathcal{H}^{1}$ over $A$, a Hopf algebroid $\mathcal{H}^{2}$ over $B^{1}$ such that $(B^{1},\mathcal{H}^{1})$, $(B^{2},\mathcal{H}^{2})$ and $(B^{2},\mathcal{H})$ are (local) noncommutative coverings of $A$, $B^{1}$ and $A$, respectively. In terms of diagrams, we have

\begin{equation} \label{eq:d2}
\xymatrix{ & & B^{2} \\
B^{1} \ar@{=>}[rru]^-{\mathcal{H}^{2}} & & \\
& & \\
& A \ar@{=>}[uul]^-{\mathcal{H}^{1}} \ar@{=>}[uuur]_-{\mathcal{H}} & \\ }.
\end{equation}

\noindent Let us denote by $\mathfrak{gal}$, $\mathfrak{gal}_{1}$ and $\mathfrak{gal}_{2}$ the respective Galois maps associated to the coactions $B^{2}\stackrel{\rho}{\longrightarrow}B^{2}\otimes_{A}\mathcal{H}$, $B^{1}\stackrel{\rho_{1}}{\longrightarrow}B^{1}\otimes_{A}\mathcal{H}^{1}$ and $B^{2}\stackrel{\rho_{2}}{\longrightarrow}B^{2}\otimes_{B^{1}}\mathcal{H}^{2}$. We say that such a diagram as \ref{eq:d2} \textit{commutes} if the following conditions are satisfied.

\begin{enumerate}
\item[(i)] There is a geometric morphism $\xymatrix{\mathcal{H}^{1}\ar[r]^-{(id,\phi)} & \mathcal{H}}$ of Hopf algebroids such that $\phi$ is injective and the following diagram commutes.

\[ \xymatrix{ B^{1}\otimes_{A}B^{1} \ar[rr]^-{\mathfrak{gal}_{1}} \ar[dd]_-{id\tens{A}id} & & B^{1}\otimes_{A}\mathcal{H}^{1} \ar[dd]^-{id\tens{A}\phi} \\
&& \\
B^{2}\otimes_{A}B^{2} \ar[rr]_-{\mathfrak{gal}} && B^{2}\otimes_{A}\mathcal{H} } \]

\item[(ii)] There is a geometric morphism $\xymatrix{\mathcal{H} \ar[r]^-{(f,\psi)} & \mathcal{H}^{2}}$ of Hopf algebroids such that $f$ is the inclusion $A\subseteq B^{1}$, $\psi$ is surjective and the following diagram commutes.

\[ \xymatrix{ B^{2}\otimes_{A}B^{2} \ar[rr]^-{\mathfrak{gal}} \ar@{->>}[dd] & & B^{2}\otimes_{A}\mathcal{H} \ar@{->>}[dd]^-{id\otimes_{f}\ \psi} \\
&& \\
B^{2}\otimes_{B^{1}}B^{2} \ar[rr]_-{\mathfrak{gal}_{2}} && B^{2}\otimes_{B^{1}}\mathcal{H}^{2} } \]

\end{enumerate}
\end{defn}

\begin{rem}
\begin{enumerate}
\item[]

\item[(1)] Note that we are suppressing a lot of notations here. First, when we denote by $\rho$ the coaction of $\mathcal{H}$ on $B^{2}$ we mean a pair of maps $\rho_{L}$ and $\rho_{R}$ as dscribed in section \ref{S2.3}. Same goes for $\rho^{1}$ and $\rho^{2}$. Correspondingly, by $\mathfrak{gal}$ we mean a pair of maps $\mathfrak{gal}_{L}$ and $\mathfrak{gal}_{R}$ associated to $\rho_{L}$ and $\rho_{R}$, respectively.

\item[(2)] At present writing of this paper, there is no existing Galois connection for Hopf-Galois extensions for Hopf algebras let alone for Hopf algebroids. The two conditions listed above are the minimum requirements one needs to have a noncommutative analogue of porposition \ref{prop3}.

\item[(3)] The above definition is specifically for local coverings. For general stratified coverings, $\mathcal{H}$ is a Hopf algebroid over $A^{'}\subset A$, $\mathcal{H}_{1}$ is a Hopf algebroid over $A_{1}\subseteq A$ and $\mathcal{H}^{2}$ is a Hopf algebroid over $A_{2}\subseteq B^{1}$. For the definition of commutativity of diagram \ref{eq:d2} in this situation, in addition to the existence of $\phi$ and $\psi$ we also assert the existence of $k$-algebra morphisms $f_{1}:A_{1}\longrightarrow A^{'}$ and $f_{2}:A_{2}\longrightarrow A^{'}$. In the appropriate diagrams, we replace $(id,\phi)$ by $(f_{1},\phi)$, $(f,\psi)$ by $(f_{2},\psi)$, $id\tens{A}\phi$ by $id\tens{f_{1}}\phi$ and $id\otimes_{f}\ \psi$ by $id\tens{f_{2}}\psi$.

\end{enumerate}
\end{rem}

If diagram \ref{eq:d2} commutes, we will refer to the local covering $A\stackrel{\mathcal{H}}{\Longrightarrow}B^{2}$ as the \textit{composition} of $A\stackrel{\mathcal{H}^{1}}{\Longrightarrow}B^{1}$ and $B^{1}\stackrel{\mathcal{H}^{2}}{\Longrightarrow}B^{2}$. Note that the commutativity of diagram \ref{eq:d2} depends on $\phi$ and $\psi$. We will call the pair $(\phi,\psi)$ the \textit{commutativity datum} of diagram \ref{eq:d2} of local coverings. The \textit{commutativity datum} of stratified coverings is the quadruple $(f_{1},f_{2},\phi,\psi)$ as described in (3) of the above remarks. The following proposition states the noncommutative analogue of the first part of proposition \ref{prop3} for local coverings. The next proposition is for uniform coverings.

\begin{prop}(Exact fitting for local coverings) \label{prop4}
Let $(B^{2},\mathcal{H})$ and $(B^{1},\mathcal{H}^{1})$ be local coverings of $A$ and let $(B^{2},\mathcal{H}^{2})$ be a local covering of $B^{1}$. Suppose the associated diagram as in \ref{eq:d2} commutes with commutativity datum $(\phi,\psi)$. Then, up to extending scalars, the composite map $\psi\circ \phi$ factors through the source map $B^{1}\stackrel{s}{\longrightarrow}\mathcal{H}^{2}$, i.e. following diagram of $k$-modules commute

\[ \xymatrix{ B^{1}\otimes_{A}\mathcal{H}^{1} \ar@{^(->}[rr]^-{id\otimes_{A}\phi} \ar[dd]_-{id\otimes_{A}\epsilon} & & B^{1}\otimes_{A}\mathcal{H} \ar@{->>}[dd]^-{id\otimes_{f} \psi} \\
& & \\
B^{1}\tens{B^{1}}B^{1} \ar[rr]_-{id\tens{B^{1}} s} & & B^{1}\tens{B^{1}}\ \mathcal{H}^{2} } \]

\noindent where $s$ denotes the pair of source maps $s_{L},s_{R}$ of $\mathcal{H}^{2}$ and $\epsilon$ denotes the pair of counit maps $\epsilon_{L},\epsilon_{R}$ of $\mathcal{H}^{1}$.
\end{prop}

\begin{prf} For $(B^{2},\mathcal{H})$ and $(B^{1},\mathcal{H}^{1})$ local coverings of $A$ and $(B^{2},\mathcal{H}^{2})$ a local covering of $B^{1}$, denote by $\mathfrak{gal},\mathfrak{gal}_{1}$ and $\mathfrak{gal}_{2}$ the associated Galois maps, respectively. Assuming diagram \ref{eq:d2} commutes with commutativity datum $(\phi,\psi)$ gives the following commutative diagram.

\[ \xymatrix{ B^{1}\tens{A}B^{1} \ar[rrrr]^-{\iota\tens{A}\iota} \ar@{->>}[dddd] \ar[ddrr]^-{\mathfrak{gal}_{1}} && && B^{2}\tens{A}B^{2} \ar@{->>}[dddd]|(0.5)\hole \ar[ddrr]^-{\mathfrak{gal}} && \\
&&  &&  &&  \\
&& B^{1}\tens{A}\mathcal{H}^{1} \ar@{->>}[ddd]_(0.33){id\tens{A}\epsilon} \ar[rrrr]^(0.66){id\otimes_{A}\phi} &&  && B^{2}\tens{A}\mathcal{H} \ar@{->>}[dddd]^-{id\otimes_{f}\psi} \\
&&  &&  && \\
B^{1}\tens{B^{1}}B^{1} \ar@{=}[ddrr] \ar[rrrr]|-\hole^(0.33){\iota\tens{B^{1}}\iota} &&  && B^{2}\tens{B^{1}}B^{2} \ar[ddrr]^-{\mathfrak{gal}_{2}} && \\
&& B^{1}\tens{A}A \ar@{=}[d] &&  && \\
&& B^{1}\tens{B^{1}}B^{1} \ar[rrrr]_-{\iota\tens{B^{1}}s} &&  && B^{2}\tens{B^{1}}\mathcal{H}^{2} \\ } \]

\noindent The top and right squares are the commutative diagrams in definition \ref{D5}. The commutativity of the back square, where the arrows going downwards are the canonical surjections, is obvious. To see the commutativity of the left square, take $b,b^{'}\in B^{1}$. Then using the left-Galois map $\mathfrak{gal}^{L}_{1}$ and the left-counit $\epsilon_{L}$ we have

\[ (id\tens{A}\epsilon)\mathfrak{gal}^{L}_{1}(b\tens{A}b^{'})= (id\tens{A}\epsilon)( bb^{'}_{[0]}\tens{A}b^{'}_{[1]}) = bb^{'}_{[0]}\tens{A}\epsilon_{L}(b^{'}_{[1]}) \]
\[= bb^{'}_{[0]}\epsilon_{L}(b^{'}_{[1]})\tens{A}1 = bb^{'}\tens{A}1 =bb^{'}\tens{B^{1}}1=b\tens{B^{1}}b^{'}. \]

\noindent Same computation holds for $\mathfrak{gal}^{R}_{1}$ and $\epsilon_{R}$. The commutativity of the bottom square is due to the fact that the module structure on $\mathcal{H}^{2}$ used to form the tensor product $B^{2}\tens{B^{1}}\mathcal{H}^{2}$ is the one provided by the source maps. Commutativity of the back, right, left, top and bottom squares imply that the front square commutes. By inspection, the front square reduce to the square asserted by the proposition. $\blacksquare$
\end{prf}

\begin{prop}(Exact fitting for uniform coverings) \label{prop5}
Let $(B^{1},H^{1})$ and $(B^{2},H)$ be uniform coverings of $A$ and $(B^{2},H^{2})$ a uniform covering of $B^{1}$. Suppose at least one of $B^{1}$ and $B^{2}$ is faithfully $k$-flat and suppose the associated diagram as in \ref{eq:d2} commutes with commutativity datum $(f_{1},f_{2},\phi,\psi)$. Then $f_{1}$ and $f_{2}$ are both equal to the identity $k$-algebra morphism $k\longrightarrow k$ and the composite map $\psi\circ \phi$ factors through $k$ via the counit $\epsilon_{1}:H^{1}\longrightarrow k$ and the unit $\eta_{2}:k\longrightarrow H^{2}$, i.e. the following diagram commutes.

\[ \xymatrix{ H^{1} \ar[rr]^-{\phi} \ar[dd]_-{\epsilon_{1}} & & H \ar[dd]^-{\psi} \\
& & \\
k \ar[rr]_-{\eta_{2}} & & H^{2} } \]

\end{prop}

\begin{prf}
Following the proof of Proposition \ref{prop4} we have a cube

\[ \xymatrix{ B^{1}\tens{A}B^{1} \ar[rrrr]^-{\iota\tens{A}\iota} \ar@{->>}[dddd] \ar[ddrr]^-{\mathfrak{gal}_{1}} && && B^{2}\tens{A}B^{2} \ar@{->>}[dddd]|(0.5)\hole \ar[ddrr]^-{\mathfrak{gal}} && \\
&&  &&  &&  \\
&& B^{1}\otimes H^{1} \ar@{->>}[dddd]_(0.25){id\otimes\epsilon_{1}} \ar[rrrr]^(0.66){id\otimes\phi} &&  && B^{2}\otimes H \ar@{->>}[dddd]^-{id\otimes\psi} \\
&&  &&  && \\
B^{1}\tens{B^{1}}B^{1} \ar[ddrr]_-{\mathfrak{gal}_{triv}} \ar[rrrr]|-\hole^(0.33){\iota\tens{B^{1}}\iota} &&  && B^{2}\tens{B^{1}}B^{2} \ar[ddrr]^-{\mathfrak{gal}_{2}} && \\
&& &&  && \\
&& B^{1}\otimes k \ar[rrrr]_-{\iota\otimes\eta_{2}} &&  && B^{2}\otimes H^{2} \\ } \]

\noindent with commuting back, right, and top faces. The bottom square commutes by viewing $B^{1}$ as a Hopf-Galois extension of $B^{1}$ with the trivial coaction of the $k$-Hopf algebra $k$. Similar computation as that of the previous proposition implies that the left square commutes as well. Thus, the front square commutes. Finally, the commutative square

\[ \xymatrix{ B^{1}\otimes H^{1} \ar[rrrr]^-{\iota\otimes\phi} \ar[dddd]_-{id\otimes\epsilon_{1}} \ar@{^(->}[rd] & & & & B^{2}\otimes H \ar[dddd]^-{id\otimes\psi} \\
& B^{2}\otimes H^{1} \ar[rrru]_-{id\otimes\phi} \ar@{->>}[dd]^-{id\otimes\epsilon_{1}} & & & \\
& & & & \\
& B^{2}\otimes k \ar[rrrd]^-{\iota\otimes\eta_{2}} & & & \\
B^{1}\otimes k \ar[rrrr]_-{\iota\otimes\eta_{2}} \ar@{^(->}[ru] & & & & B^{2}\otimes H^{2} } \]

\noindent and the faithfully $k$-flatness, say of $B^{2}$, implies the desired result. $\blacksquare$
\end{prf}

\begin{rem}
Note that the commutativity of the diagram in Proposition \ref{prop5} is the naive analogue of exactness of the sequence $H^{1}\stackrel{\phi}{\longrightarrow}H\stackrel{\psi}{\longrightarrow}H^{2}$ of Hopf algebras from the view-point of groups algebras. However, this is not the usual notion of exactness as the zero object in the category of $k$-Hopf algebras is the zero $k$-algebra $\left\{0\right\}$ and not $k$.
\end{rem}

%%%%%%%%%%%%%%%%%%%%%%%%%%%%%%%%%%%%%%%%%%%%%%%%%%%%%%%%%%%%%%%%%%%%%%%%%%%%%%%%%%%%%%%%

\section{Coverings of commutative spaces: Central Case.}\label{S4.0}

We mentioned in the introduction, the formulation of a noncommutative covering space should be guided by the following: (1) they should give, as a special case (when the symmetry is a Hopf algebra), noncommutative principal bundles as currently understood (see for example \cite{m003}), (2) when the algebras involved are commutative then we should be able to get a classical covering spaces i.e., a reconstruction procedure. We will state this reconstruction theorem in this section.
We will only deal with central coverings here, i.e. coverings in which $A$ sits centrally in both $B$ and $\mathcal{H}$. The non-central case is discussed in \cite{canlubo}.

\subsection{Coverings of a point}\label{S4.1}

In this section, we will have a closer look at coverings of a point. In particular, we will see that unlike the classical case, a point has infinitely many connected covers. Also, we will characterize the type of Hopf algebroids $\mathcal{H}$ that can arise in a covering $(B,\mathcal{H})$. In noncommutative geometry, a point is represented by the base ring in consideration. In this section, let us fix the base commutative unital ring $k$.

A priori, a covering of a point is a pair $(B,\mathcal{H})$ where $\mathcal{H}$ is a Hopf algebroid over $k$ and $k\subseteq B$ is a right $\mathcal{H}$-Galois extension. In the literature, $B$ is called a Hopf-Galois object over $k$. Let us give some examples of such coverings.

Given any finitely generated projective Hopf algebra $H$ over $k$, we claim that $(H,H)$ is a covering of a point. Here, we use the regular coaction of $H$ on itself. The left and right-bialgebroid structures of $H$ are both isomorphic to the underlying bialgebra of $H$. All that is left to show is that the Galois map

\[ \xymatrix@R=2mm{ H\otimes H \ar[rr]^-{\mathfrak{gal}} & & H\otimes H \\
a\otimes b  \ar@{|->}[rr] & & ab_{(1)}\otimes b_{(2)} } \]

\noindent is bijective. This is the case, since the map

\[ \xymatrix@R=2mm{ H\otimes H \ar[rr] & & H\otimes H  \\
a\otimes b \ar@{|->}[rr] & & aS(b_{(1)})\otimes b_{(2)}\\ } \]

\noindent is its inverse. In fact more is true, a bialgebra $H$ is a Hopf algebra if and only if it is an $H$-Galois extension of the base ring. This tells us that any connected $k$-Hopf algebra is a connected covering of a point. By a connected Hopf algebra we mean connected as an algebra i.e. one in which the only idempotent elements are $0$ and $1$.

Now, let us look at a more general situation. Let $(B,\mathcal{H})$ be a (finite) covering of $k$. Explicitly, this means that $B$ is a $k$-algebra which is finitely generated and projective as a $k$-module. Also, $\mathcal{H}=(\mathcal{H}_{L},\mathcal{H}_{R},S)$ where $\mathcal{H}_{L}=(H,s_{L},t_{L},\Delta_{L},\epsilon_{L})$ and $\mathcal{H}_{R}=(H,s_{R},t_{R},\Delta_{R},\epsilon_{R})$.

We claim that $\mathcal{H}_{L}$ is a bialgebra. The source and target maps $s_{L}$ and $t_{L}$ define a $k$-algebra map $\eta_{L}=s_{L}\otimes t_{L}:k\longrightarrow H$. The product $\mu_{L}$ on $H$ determined by $\eta_{L}$ is associative and counital with respect to $\eta_{L}$. The coproduct $\Delta_{L}:H\longrightarrow H\otimes H$ is already a $k$-algebra map since the Takeuchi product $H\prescript{}{k}{\times} H$ and $H\otimes H$ coincide. $\Delta_{L}$ is coassociative and counital with respect to $\epsilon_{L}$. All that is left to show is that $\epsilon_{L}:H\longrightarrow k$ is a $k$-algebra map. Part (c) of the definition of a bialgebroid implies that $\epsilon_{L}$ is unital, i.e. $\epsilon_{L}(1)=1$. Applying theorem 5.5 of Schauenburg \cite{sch4} using the identity map $k\longrightarrow k$ and the normalized dual basis of $k$ given by the unit element, we see that $H$ possesses a weak bialgebra structure with coproduct $\Delta_{L}$ and counit $\epsilon_{L}$. This implies that $\epsilon_{L}(xy)=\epsilon_{L}(x1_{[1]})\epsilon_{L}(1_{[2]}y)$ for any $x,y\in H$. But $1\otimes 1=\Delta_{L}(1)=1_{[1]}\otimes 1_{[2]}$. Thus, $\epsilon_{L}$ is a unital $k$-algebra map. This shows that indeed $\mathcal{H}_{L}$ is a bialgebra over $k$.

Now, $\mathcal{H}_{L}$ admits a Galois extension which is $k\subseteq B$ in this case. By a result of Schauenburg \cite{sch1}, the bialgebra $\mathcal{H}_{L}$ is in fact a Hopf algebra, i.e. there is a $k$-module map $S_{L}:H\longrightarrow H$ such that $\mathcal{H}_{L}=(H,\mu_{L},\eta_{L},\Delta_{L},\epsilon_{L},S_{L})$ is a Hopf algebra over $k$. Similar argument shows that there is a $k$-module map $S_{R}:H\longrightarrow H$ making $\mathcal{H}_{R}=(H,\mu_{R},\eta_{R},\Delta_{R},\epsilon_{R},S_{R})$ a Hopf algebra over $k$.

The antipode $S$ of the Hopf algebroid $\mathcal{H}$ provides a coupling map making $\mathcal{H}_{L}$ and $\mathcal{H}_{R}$ coupled Hopf algebras. Thus, we have proved the following proposition.

\begin{prop} \label{prop6}
For $(B,\mathcal{H})$ a covering of a point, $\mathcal{H}$ is a Hopf algebroid coming from coupled Hopf algebras.
\end{prop}

%%%%%%%%%%%%%%%%%%%%%%%%%%%%%%%%%%%%%%%%%%%%%%%%%%%%%%%%%%%%%%%%%%%%%%%%%%%%%%%%%%%%%%%%

\subsection{Commutative coverings of commutative spaces}\label{S4.2}

As we have seen in section \ref{S3.1}, finite Galois (connected) classical covering $Y\stackrel{p}{\longrightarrow}X$ gives a covering $(C(Y),C(\mathcal{G}))$ (in the sense of definition \ref{D2}) where $\mathcal{G}$ is the groupoid we constructed in section \ref{S3.1} and $\rho:C(Y)\longrightarrow C(Y)\otimes_{C(X)} C(\mathcal{G})$ is the induced coaction from the pointwise deck action of $\mathcal{G}$ on $Y$. Conversely, let us show that commutative examples give classical covering spaces. Through out this section, we will restrict our attention to local coverings. We will proceed in two ways, one for commutative $C^{*}$-algebras and the other one for general commutative unital ring $R$.

Let $A$ and $B$ be a commutative unital $C^{*}$-algebras and $\mathcal{H}$ a finitely-generated projective Hopf algebroid over $A$ coacting on $B$ such that $(B,\mathcal{H})$ is a covering of $A$ in the sense of definition \ref{D2}. $A$ and $B$ being commutative implies that $B\otimes_{A}B$ carries an algebra structure by tensorwise product. The Galois maps

\[ \xymatrix@R=2mm{ B\tens{A}B \ar[rr]^-{\mathfrak{gal}_{L}} & & B\tens{A}\mathcal{H}_{L} \\
a\tens{A} b  \ar@{|->}[rr] & & ab_{[1]}\tens{A} b_{[2]} } \hspace{.5in} \xymatrix@R=2mm{ B\tens{A}B \ar[rr]^-{\mathfrak{gal}_{R}} & & B\tens{A}\mathcal{H}_{R} \\
a\tens{A} b  \ar@{|->}[rr] & & ab^{[1]}\tens{A} b^{[2]} } \]

\noindent then become algebra maps. To see this, given $a\otimes_{A}b,a^{'}\otimes_{A}b^{'}\in B\otimes_{A}B$ we have

\begin{eqnarray*} \mathfrak{gal}_{L}\left(\left(a\tens{A}b\right)\left(a^{'}\tens{A}b^{'}\right)\right) &=& \mathfrak{gal}_{L}\left(aa^{'}\tens{A}bb^{'}\right)=aa^{'}b_{[0]}b^{'}_{[0]}\tens{A}b_{[1]}b^{'}_{[1]} \\
&=& \left(ab_{[0]}\tens{A}b_{[1]}\right)\left(a^{'}b^{'}_{[0]}\tens{A}b^{'}_{[1]}\right)\\
&=& \mathfrak{gal}_{L}\left(a\tens{A}b\right)\mathfrak{gal}_{L}\left(a^{'}\tens{A}b^{'}\right).
\end{eqnarray*}

\noindent $\mathfrak{gal}_{L}$ being a linear bijection implies that $\mathcal{H}_{L}$ is a commutative $A$-bialgebroid. Similar computation using $\mathfrak{gal}_{R}$ shows that $\mathcal{H}_{R}$ is a commutative $A$-bialgebroid.

By the Gelfand duality, there are compact Hausdorff spaces $\widehat{A}$ and $\widehat{B}$ such that $A=C(\widehat{A})$ and $B=C(\widehat{B})$. Explicitly, $\widehat{B}$ is the sets of unital homomorphisms $B\stackrel{\varphi}{\longrightarrow}\mathbb{C}$. $B$ being commutative unital Banach algerba forces $\left\|\varphi\right\|=1$. Thus, $\widehat{B}\subset B^{*}$ and we can equip $\widehat{B}$ with the subspace topology it inherits from the weak$-*$ topology on $B^{*}$. Similarly, we can topologize $\widehat{A}$ this way. The inclusion $A\subset B$ induces a projection $\widehat{B}\stackrel{p}{\longrightarrow}\widehat{A}, \varphi\mapsto\varphi|_{A}$. We claim that this is a classical covering space.

First, we need the following lemma generalizing the result in algebraic geometry saying that the category of commutative Hopf algebras is dual to the category of affine group schemes in a particular way.

\begin{lem}\label{L1}
Let $\mathcal{H}=\left(H_{L},H_{R},S\right)$ be a commutative Hopf algebroid (i.e. one whose constituent bialgebroids are commutative) over a commutative algebra $A$ with bijective antipode $S$. Then there is a topological groupoid $\mathcal{G}$ whose algebra of continuous functions is by $\mathcal{H}$.
\end{lem}

\begin{prf}
Applying the $Spec$ functor in the following diagram of commutative $A$-algebras describing the Hopf algebroid $\mathcal{H}$

\[\xymatrix{& & H_{R}\tens{A}H_{R} & & &\\
& & & H_{R} \ar[lu]_-{\Delta_{R}} \ar[rr]^-{\epsilon_{R}} \ar@(l,u)[ldd] & & A \ar@(ul,ur)[ll]_-{s_{R}} \ar@(dl,dr)[ll]^-{t_{R}} \\
& & & & & \\
A \ar@(ur,ul)[rr]^-{s_{L}} \ar@(dr,dl)[rr]_-{t_{L}} & & H_{L} \ar[rd]_-{\Delta_{L}} \ar[ll]_-{\epsilon_{L}} \ar@(r,d)[ruu] \ar@{}[ruu]|-{S} & & & \\
& & & H_{L}\tens{A}H_{L}. & &\\ } \]

\noindent gives topologically enriched small categories $\mathcal{C}_{R}=Spec(H_{R})$ and $\mathcal{C}_{L}=Spec(H_{L})$ over $X=Spec(A)$. To be precise, the underlying space of arrows of these categories come from the commutative $A$-ring structures of $H_{L}$ and $H_{R}$. The categorical compositions and the units come from the $A$-coring structures. We abuse notation by writing $\mathcal{C}_{R}$ (resp. $\mathcal{C}_{L}$) for the space of arrows of the category $\mathcal{C}_{R}$ (resp. $\mathcal{C}_{L}$). Note that $\mathcal{C}_{L}$ and $\mathcal{C}_{R}$ have the same underlying space $C$ as this space is precisely $Spec(H)$ where $H$ is the common underlying $k$-algebra of $H_{L}$ and $H_{R}$.

The antipode $S$ induces a continuous map $C\stackrel{F_{S}}{\longrightarrow}C$. The following diagram of spaces describes the properties of $F_{S}$ in relation with the rest of the categorical structures of $\mathcal{C}_{L}$ and $\mathcal{C}_{R}$.

\begin{equation} \label{eq4.2.1}
\xymatrix{
& & C \times C \ar[rr]^-{F_{S}\times id} & & C\prescript{}{t_{L}}{\times}_{s_{L}} C \ar[rrd]^-{\circ_{L}} & & \\
C \ar[rru]^-{diag} \ar[rrr]^-{s_{R}} & & & X \ar[rrr]^-{\epsilon_{R}} & & & C \\
C \ar[rrd]_-{diag} \ar[rrr]^-{s_{L}} & & & X \ar[rrr]^-{\epsilon_{L}} & & & C \\
& & C \times C \ar[rr]_-{id\times F_{S}} & & C\prescript{}{t_{R}}{\times}_{s_{R}} C \ar[rru]^-{\circ_{R}} & & \\}
\end{equation}

\noindent Here, we denoted by the same notation the maps induced by the source, target and counit maps. As we mentioned above, the counit maps induced the unit maps of the two categories. By part $(1)$ of remark 2, we see that the orientations of elements of $C$ viewed as arrows of $\mathcal{C}_{L}$ are opposite those orientations when viewed as arrows of $\mathcal{C}_{R}$. In particular, this means that the two categories have the same units. Using this fact, we can show that more is true. The two categories are groupoids. Let us show that any $\varphi\in \mathcal{C}_{R}$ is invertible. Using the lower part of diagram \ref{eq4.2.1} implies that for any $f\in H$, we have

\begin{eqnarray*}
f\left(\varphi\circ_{R} F_{S}(\varphi)\right) &=& f^{[1]}(\varphi)f^{[2]}(F_{S}(\varphi))=f^{[1]}(\varphi)S(f^{[2]})(\varphi)\\
&=& f^{[1]}S(f^{[2]})(\varphi)=\left(s_{L}\circ\epsilon_{L}\right)(f)(\varphi)=\epsilon_{L}(f)(s_{L}(\varphi)) = f(id_{s_{L}(\varphi)}).
\end{eqnarray*}

\noindent Thus, $F_{S}(\varphi)$ is the inverse of $\varphi$ in the category $\mathcal{C}_{R}$. The proof for $\mathcal{C}_{L}$ being a groupoid goes the same way.

At this point, we have two groupoids $\mathcal{C}_{L}$ and $\mathcal{C}_{R}$ whose space of units coincide. Recall that the categorical compositions $\circ_{L}$ and $\circ_{R}$ are functorially induced by the coproducts $\Delta_{L}$ and $\Delta_{R}$, respectively. These coproducts commute. Thus, the categorical compositions $\circ_{L}$ and $\circ_{R}$ commute as well. By the groupoid version of Eckmann-Hilton argument, the two compositions are the same. This shows that the groupoids are opposite each other. One can pick either of these groupoids to get the groupoid asserted by the lemma. $\blacksquare$
\end{prf}

\

\begin{rem}
The proof above provides adjoint equivalence between the category of commutative Hopf algebroids and groupoid schemes. This is formally the same as the adjoint equivalence between commutative Hopf algebras and affine group schemes. The only additional ingredient is Grothendieck's relative point of view for schemes. This may lead one to think that $\mathcal{H}$ being a Hopf algebroid over a commutative algebra $R$, $\mathcal{H}$ is simply a Hopf algebra over $R$. This need not be the case, see for example weak Hopf algebras in \ref{S2.2}. Also, proof of the lemma involves a construction inverse to the one we had when we constructed Hopf algebroids from groupoid in \ref{S3.1}. Finally, we note that the constituent bialgebroids of a commutative Hopf algebroid coincide. However, Hopf algebroids exist with one consituent bialgebroid is commutative while the other one is not.
\end{rem}

\

The coaction $B\stackrel{\rho}{\longrightarrow}B\otimes_{A}\mathcal{H}$ defines a groupoid action $\widehat{B}\prescript{}{p}\times_{s}\mathcal{G}\stackrel{\alpha}{\longrightarrow}\widehat{B}$ as follows. Using lemma \ref{L1} we have an isomorphism $B\otimes_{A}B=C(\widehat{B})\otimes_{C(\widehat{A})} C(\mathcal{G})\cong C(\widehat{B} \prescript{}{p}\times_{s}\mathcal{G})$, we can write $C(\widehat{B})\stackrel{\rho}{\longrightarrow}C(\widehat{B}\prescript{}{p}\times_{s} \mathcal{G})$. Define the action $\widehat{B}\times G\stackrel{\alpha}{\longrightarrow}\widehat{B}$ as: for any $\varphi\in\widehat{B}$ and $g\in \mathcal{G}$ such that $p(\varphi)=s(g)$, $\varphi\cdot g\in \widehat{B}$ is defined for any $b\in B$ as $(\varphi \cdot g)(b)=\rho(b)(\varphi,g^{-1})$.Using the identification $C(\widehat{B})\cong B$ we have $(\varphi\cdot g)(b)=\varphi(b_{[0]})b_{[1]}(g^{-1})$. Let us show that indeed, this defines an action. Let $e\in \mathcal{G}$ be an identity arrow of $\mathcal{G}$. Then for any $\varphi\in \widehat{B}$ and $b\in B$ with $p(\varphi)=s(e)$ we have

\[ (\varphi \cdot e)(b)=\varphi(b_{[0]})b_{[1]}(e)=\varphi(b_{[0]})\varepsilon(b_{[1]})=\varphi(b_{[0]}\varepsilon(b_{[1]}))=\varphi(b) \]

\noindent using the definition of the counit $\varepsilon$ of $\mathcal{H}$ and the counit axiom, respectively. Thus, units $e\in G$ act trivially as desired. For the associativity of the action, let $\varphi\in \widehat{B}$ and $g_{1},g_{2}\in \mathcal{G}$ with $p(\varphi)=s(g_{1})=s(g_{2})$. Then for any $b\in B$ we have

\begin{eqnarray*} ((\varphi\cdot g_{1})\cdot g_{2})(b) &=& (\varphi\cdot g_{2})(b_{[0]})b_{[1]}(g_{1}^{-1})=\varphi(b_{[0][0]})b_{[0][1]}(g_{2}^{-1})b_{[1]}(g_{1}^{-1})\\
&=& \varphi(b_{[0]})b_{[1][0]}(g_{2}^{-1})b_{[1][1]}(g_{1}^{-1})=\varphi(b_{[0]})b_{[1]}(g_{2}^{-1}g_{1}^{-1})\\
&=& \varphi(b_{[0]})b_{[1]}((g_{1}g_{2})^{-1})=(\varphi\cdot(g_{1}g_{2}))(b)\\
\end{eqnarray*}

\noindent using the coassociativity of $\rho$ and the definition of the comultiplication on $\mathcal{H}$, respectively.

Let us show that $\widehat{B}/\mathcal{G}\cong\widehat{A}$. Notice that for $g\in \mathcal{G}$ and $\varphi\in\widehat{B}$ with $p(\varphi)=s(g)$, $\varphi \cdot g$ defines the same function on the set of all $b\in B$ for which $b_{[0]}=b$ and $b_{[1]}=1$. Thus, such $b\in B$ satisfies $\rho(b)=b\otimes 1$ which implies that $b\in A$. Thus, classes in $\widehat{B}/\mathcal{G}$ defines an element of $\widehat{A}$. Conversely, any element in $\widehat{A}$ is invariant under the induced action of $\mathcal{G}$. Thus, we have a commutative diagram of $\mathcal{G}$-equivariant continuous maps

\[\xymatrix{
\widehat{B} \ar[rrd(0.8)]_-{proj} \ar[r(4)]^-{p} & & & & \widehat{A} \\
& & \widehat{B}/\mathcal{G} \ar@{=}[urr(0.9)]_-{\cong} & & } \]

\noindent This in particular shows that $\mathcal{G}$ acts by deck transformations on $\widehat{B}\stackrel{p}{\longrightarrow}\widehat{A}$. This means that $\widehat{B}\stackrel{p}{\longrightarrow}\widehat{A}$ is a covering space of degree the order of fiber groups of $\mathcal{G}$.

However, $(B,\mathcal{H})$ being a covering space of $A$ is giving us more. In particular, this tells us that $\widehat{B}\stackrel{p}{\longrightarrow}\widehat{A}$ is in fact a Galois covering. This follows immediately from the fact that $B\otimes_{A}B\stackrel{\mathfrak{gal}}{\longrightarrow}B\otimes_{A} \mathcal{H}$ is bijective. At the level of topological spaces, $\mathfrak{gal}$ induces the corresponding bijective Galois map $\widehat{B}\prescript{}{p}\times_{s} \mathcal{G}\stackrel{\mathfrak{gal}^{'}}{\longrightarrow}\widehat{B}\times_{\widehat{A}}\widehat{B}$, showing that fiberwise, $\mathcal{G}$ acts transitively. Thus, we have shown the following theorem.

\begin{thm}\label{thm3}
Let $A$ be a commutative $C^{*}$-algebra. Let $(B,\mathcal{H})$ be a local covering of $A$ with $B$ a commutative $C^{*}$-algebra and $\mathcal{H}$ a commutative Hopf algebroid with coinciding constituent bialgebroids and unital structure maps. Then there is a classical finite Galois covering $Y\stackrel{p}{\longrightarrow}X$ with finite deck transformation group $G$ such that $A=C(X)$, $B=C(Y)$ and $G$ is the vertex group of the groupoid $\mathcal{G}(\mathbb{C})$ where $\mathcal{G}$ is the groupoid scheme determined by $\mathcal{H}$.
\end{thm}

Now let us look at the case of general commutative rings. Let $k$ be a commutative unital ring and $A$ a commutative algebra over $k$. Let $(B.\mathcal{H})$ be a covering of $A$ with $B$ a commutative algebra and $\mathcal{H}$ a commutative Hopf algebroid with coinciding constituent bialgebroids and bijective antipode.

The inclusion $A\subseteq B$ gives a surjective map $\xymatrix{Spec(B) \ar@{->>}[r]^-{p} & Spec(A)}$. Similar to the case of $C^{*}$-algebras, the coaction $B\stackrel{\rho}{\longrightarrow}B\otimes_{A}\mathcal{H}$ gives an action $Spec(B)\times_{Spec(A)}\mathcal{G}\stackrel{\alpha}{\longrightarrow}Spec(B)$. Since the coinvariants of the coaction $\rho$ is $A$, we have $Spec(B)/\mathcal{G}\cong Spec(A)$. Bijectivity of the Galois map $B\otimes_{A}B\stackrel{\mathfrak{gal}}{\longrightarrow}B\otimes_{A} \mathcal{H}$ translates to bijectivity of the following map.

\[\xymatrix{ Spec(B)\times_{Spec(A)} \mathcal{G} \ar[rr]^-{\mathfrak{gal}^{'}} & & Spec(B)\times_{Spec(A)}Spec(B) }\]

\noindent This tells us that theorem \ref{thm3} is valid in the case of a general commutative rings.

\subsection{Noncommutative coverings of commutative spaces}\label{S4.3}

Let $A$ be a commutative unital $C^{*}$-algebra. Let $(B,\mathcal{H})$ be a local covering of $A$, where $B$ is a unital $C^{*}$-algebra, $A\subseteq Z(B)$, and the image of $A$ under the source and target maps lie in the center of $\mathcal{H}$. We regard this situation as $A$ being \textit{central} in $(B,\mathcal{H})$ or that $(B,\mathcal{H})$ is a \textit{central} covering. We assume that the left and right coactions of $\mathcal{H}$ on $B$ are continuous.

By Gelfand-Naimark duality, $A=C(X)$ where $X$ is a compact Hausdorff space. Specifically, $X$ is the spectrum of $A$, the space of unitary equivalence classes of irreducible $*$-representations of $A$. Since $A$ is commutative, $X$ coincides with the primitive spectrum of $A$, the space of primitive ideals of $A$ with the hull-kernel topology. Since $B$ is a finitely-generated projective module over $C(X)$, the Serre-Swan theorem implies that $B\cong\Gamma(X,E)$ for some finite rank vector bundle $\xymatrix{E \ar@{->>}[r]^-{p} & X}$

Let $x\in X$ and let $B_{x}=\left\{\sigma\in \Gamma(X,E)| \sigma(x)=0\right\}$. Then $B_{x}$ is an ideal of $B$. To see this, given any $\sigma\in B_{x}$, write $\sigma=f\cdot \sigma^{'}$ for some $\sigma^{'}\in B$ and $f\in C(X)$ such that $f(x)=0$. Now, given any $\tau\in B$, we have $\left(\sigma\tau\right)(x)=f(x)\left(\sigma^{'}\tau\right)(x)=0$. Centrality of $A$ in $B$ implies that $B/B_{x}$ is an $\mathbb{C}$-algebra where we identify $\mathbb{C}$ with $A/I_{x}$, $I_{x}=\left\{f\in A|f(x)=0\right\}$.

The evaluation map $ev_{x}:B\longrightarrow E$ at $x$ lifts to a map $e:B/B_{x}\longrightarrow E$. Since $E_{x}$ is the pullback of $x\longrightarrow X \longleftarrow E$, we have a linear map $\varphi$ such that the following diagram commutes

\[\xymatrix{
& & & E \ar[dr] & \\
B/B_{x} \ar@/^/[urrr]^{ev_{x}} \ar@/_/[drrr]
\ar@{-->}[rr]^-{\varphi} & & E_{x} \ar[ur] \ar[dr] & & X .\\
& & & x \ar[ur] & }\]

\noindent In fact, $\varphi$ is an isomorphism from $B/B_{x}$ to $E_{x}$. To see this, note that any element $e\in E_{x}$ can be extended to a section $\sigma\in B$ and any other extension is a section having the same value $e$ at $x$. Thus, they define the same element in $B/B_{x}$. This gives us the following proposition.

\begin{prop} \label{prop7}
Let $A\subseteq B$ be an algebra extension with $A=C(X)$ central in $B$ and $B$ finitely generated and projective as a regular $A$-module. Then $B$ is a bundle of complex algebras over $X$ such that the algebra structure of $B$ is pointwise.
\end{prop}

Let us assume that the images of $C(X)$ under the source and target maps are central. Consider the left bialgebroid structure $\mathcal{H}_{L}$ of $\mathcal{H}$. The $A$-bimodule structure of $\mathcal{H}_{L}$ is finitely-generated and projective in the sense that both the constituent module structures are finitely-generated and projective. In particular, using the same argument we see that as a left $A$-module, $\mathcal{H}_{L}\cong\Gamma(X,H^{L})$ for some finite rank vector bundle $\xymatrix{H^{L} \ar@{->>}[r]^-{q} & X}$. Moreover, each fiber has an algebra structure such that the $A$-ring structure on $\mathcal{H}_{L}$ is isomorphic to the $A$-ring structure one gets by pointwise multiplication in $\Gamma(X,H^{L})$.

By Serre-Swan theorem, the covariant functor $\Gamma(X,-)$ has a left adjoint $Sh$

\[ \left\{ \begin{array}{c}
\text{finitely-generated}\\
\text{projective module} \\
\text{over } C(X) \end{array} \right\}
\begin{array}{c}
\stackrel{Sh}{\xrightarrow{\hspace*{2cm}}}\\
\stackrel{\xleftarrow{\hspace*{2cm}}}{\Gamma(X,-)}\\
\end{array}
\left\{ \begin{array}{c}
\text{finite rank}\\
\text{vector bundle} \\
\text{over } X \end{array} \right\}. \]

\noindent Explicitly, for a $C(X)$-module $M$ the vector bundle $Sh(M)$ is constructed as follows. Let $\mathcal{O}_{X}$ denote the structure sheaf of $X$ and define the presheaf $P(M)$ of $\mathcal{O}_{X}$-modules by

\[ P(M)(U)=M\otimes_{C(X)}\mathcal{O}_{X}(U) \]

\noindent and denote by $Sh(M)$ its sheafification.

Applying the functor $Sh$ to the coproduct $\xymatrix{\mathcal{H}_{L} \ar[r]^-{\Delta_{L}} & \mathcal{H}_{L}\tens{A}\mathcal{H}_{L}}$ gives a map

\[ \xymatrix{H^{L} \ar[rr]^-{Sh(\Delta_{L})} & & H^{L}\otimes H^{L}} \]

\noindent of vector bundles. By definition, the fiber of $H^{L}\otimes H^{L}$ at $x\in X$ is $H^{L}_{x}\otimes H^{L}_{x}$. Thus, there is a linear map $\delta_{L,x}$ making the following diagram commute.

\begin{equation} \label{eq:9}
\xymatrix{
H^{L}_{x} \ar[rrr] \ar[dddd] \ar@{-->}[rdd]^-{\delta_{L,x}}  & & & H^{L} \ar[ddl]|-{Sh(\Delta_{L})} \ar@{->>}[dddd] \\
& & \\
 & H^{L}_{x}\otimes H^{L}_{x} \ar[r] \ar[ldd] & H^{L}\otimes H^{L} \ar@{->>}[ddr] \\
& & \\
x \ar[rrr] & & & X\\}
\end{equation}

Viewing $A$ itself as a finitely-generated projective module over $C(X)$ and applying the functor $Sh$ on the counit map $\xymatrix{\mathcal{H}_{L} \ar[r]^-{\epsilon_{L}} & A}$ gives a map $\xymatrix{H^{L} \ar[rr]^-{Sh(\epsilon_{L})} & & \mathbb{C}_{triv} }$ of vector bundles, where $\mathbb{C}_{triv}$ denotes the trivial line bundle $X\times\mathbb{C}$ over $X$. Since $H^{L}_{x}$ is the pullback of the diagram $x\longrightarrow X \longleftarrow H^{L}$ we see that we get a linear map $\xymatrix{H^{L}_{x} \ar[r]^-{\epsilon_{L,x}} & \mathbb{C} }$.

We claim that $\delta_{L,x}$ is coassociative and counital with respect to $\epsilon_{L,x}$. The back face of the following cube commutes by coassociativity of $\Delta_{L}$ and functoriality of $Sh$

\begin{equation} \label{eq:10}
\xymatrix{
& H^{L} \ar[rr]^-{Sh(\Delta_{L})} \ar'[d][dd]^-{Sh(\Delta_{L})} && \makebox[\widthof{$H^{L}$}][l]{$H^{L}\otimes H^{L}$} \ar[dd]^-{id\otimes Sh(\Delta_{L})} \\
H^{L}_{x} \ar[ru] \ar[rr]^(.65){\delta_{L,x}} \ar[dd]_-{\delta_{L,x}} && H^{L}_{x}\otimes H^{L}_{x} \ar[ru] \ar[dd]^(.25){id\otimes \delta_{L,x}} & \\
& H^{L}\otimes H^{L} \ar'[r]_-{Sh(\Delta_{L})\otimes id}[rr] && \makebox[\widthof{$H^{L}$}][l]{$H^{L}\otimes H^{L}\otimes H^{L}$} \\
H^{L}_{x}\otimes H^{L}_{x} \ar[ru] \ar[rr]_-{\delta_{L,x}\otimes id} && H^{L}_{x}\otimes H^{L}_{x}\otimes H^{L}_{x} \ar[ru]_-{\psi} & }
\end{equation}

\noindent while the lateral faces of diagram \ref{eq:10} commute since they are essentially the upper commuting square of diagram \ref{eq:9}. Commutativity of the five faces and the fact that the map $\psi$ of diagram \ref{eq:10} is injective implies that the front face commutes, i.e. $\delta_{L,x}$ is coassociative. Using the same line of reasoning, we can show counitality of $\delta_{L,x}$ with respect to $\epsilon_{L,x}$ using the leftmost diagram in (10) below.

\begin{equation} \label{eq:11}
\xymatrix{
& H^{L}\otimes H^{L} \ar[ld]_-{id\otimes \epsilon_{L}} \ar'[d][dd]^-{\epsilon_{L}\otimes id} && \makebox[\widthof{$H^{L}$}][l]{$H^{L}_{x}\otimes H^{L}_{x}$} \ar[dd]^-{\epsilon_{L,x}\otimes id} \ar[ll] \ar[ld]_-{id\otimes \epsilon_{L,x}} \\
H^{L}\otimes\mathbb{C}_{triv} \ar@{=}[dd] && H^{L}_{x}\otimes \mathbb{C} \ar@{=}[dd] \ar[ll] & & &  \\
& \mathbb{C}_{triv}\otimes H^{L} \ar@{=}[ld]  && \makebox[\widthof{$H^{L}$}][l]{$\mathbb{C}\otimes H^{L}_{x}$} \ar'[l]|(0.65)\hole[ll] \\
H^{L} \ar[ruuu]|(0.66)\hole^-{Sh(\Delta_{L})}  && H^{L}_{x} \ar@{=}[ru] \ar[ll] \ar[ruuu]^-{\delta_{L,x})} &  \left(H^{L}\right)^{\otimes 2} \ar[rr]^-{Sh(\Delta_{L})\otimes Sh(\Delta_{L})} \ar'[d][dd]_-{m} && \makebox[\widthof{$H^{L}$}][l]{$\left(H^{L}\right)^{\otimes 4}$} \ar[dd]|-{(m\otimes m)\circ \reflectbox{F}} \\
&& \left(H^{L}_{x}\right)^{\otimes 2} \ar[ru] \ar[rr]^(.65){\delta_{L,x}\otimes \delta_{L,x}} \ar[dd]_-{m} && \left(H^{L}_{x}\right)^{\otimes 4} \ar[ru] \ar[dd]^(.25){(m\otimes m)\circ \reflectbox{F}} & \\
&& & H^{L} \ar'[r]_-{Sh(\Delta_{L})}[rr] && \makebox[\widthof{$H^{L}$}][l]{$\left(H^{L}\right)^{\otimes 2}$} \\
&& H^{L}_{x} \ar[ru] \ar[rr]_-{\delta_{L,x}} && \left(H^{L}_{x}\right)^{\otimes 2} \ar[ru] & }
\end{equation}

\noindent whose front, back, top, bottom and left faces are easily seen to commute implying that the right face is commutative as well. In the leftmost diagram, we denoted by $\mathbb{C}\otimes H^{L}$ the tensor product of the trivial line bundle $X\times\mathbb{C}$ and the bundle $H^{L}$.

We also claim that $\delta_{L,x}$ is multiplicative. This follows from the commutativity of the rightmost cube in diagram \ref{eq:11} above.

Thus, each fiber $H^{L}_{x}$ carries a multiplicative coring structure such that the (left) operations on $\mathcal{H}$ are pointwise, i.e $\mathcal{H}_{L}=\Gamma(X,H^{L})$ is an isomorphism, not just of $C(X)$-modules but also of $A$-rings and $A$-corings. This defines a left $\mathbb{C}$-biagebroid structure on $H^{L}$.

Carrying out the same arguments for the right bialgebroid structure $\mathcal{H}_{R}$ of $\mathcal{H}$, we get a finite rank vector bundle $\xymatrix{H^{R} \ar@{->>}[r]^-{r} & X}$ such that $\mathcal{H}_{R}\cong\Gamma(X,H^{R})$ as right $A$-modules. Each fiber $H^{R}_{x}$ of $H^{R}$ carries an algebra structure such that the $A$-ring structure of $\mathcal{H}_{R}$ is isomorphic to the $A$-ring structure of $\Gamma(X,H^{R})$ given by pointwise multiplication. Also, $H^{R}_{x}$ carries a multiplicative coring structure such that the coring structure on $\mathcal{H}_{R}$ is pointwise. Symmetrically, we get a right $\mathbb{C}$-bialgebroid structure on $H^{R}$.

The antipode $S$ defines a $\mathbb{C}$-module map $\Gamma(X,H^{L})\stackrel{S}{\longrightarrow}\Gamma(X,H^{R})$. Part $(c)$ of definition \ref{D1} implies that $S$ induces a fiberwise linear map $H^{L}\stackrel{S^{\wedge}}{\longrightarrow}H^{R}$. We then have the following commutative diagram.

\[ \xymatrix{
& & H^{L}_{x}\otimes H^{L}_{x} \ar[rrrr]^-{S^{\wedge}_{x}\otimes id} & & & & H^{L}_{x}\otimes H^{L}_{x} \ar[rrd]^-{\mu_{L,x}} & & \\
H^{L}_{x} \ar@{-}[d]|-{\simeq} \ar[rru]^-{\delta_{L,x}} \ar[r(4.4)]^-{\epsilon_{R}} & & & & \mathbb{C} \ar[rrrr]^-{s_{R}} & & & & H^{L}_{x} \ar@{-}[d]|-{\simeq} \\
H^{R}_{x}  \ar[rrd]_-{\delta_{R,x}}  \ar[rrrr]_-{\epsilon_{L}} & & & & \mathbb{C} \ar[rrrr]_-{s_{L}} & & & & H^{R}_{x} \\
& & H^{R}_{x}\otimes H^{R}_{x} \ar[rrrr]_-{id\otimes S^{\wedge}_{x}} & & & & H^{R}_{x}\otimes H^{R}_{x} \ar[rru]_-{\mu_{R,x}} & & \\
} \]

\noindent Thus, we have the following result.

\begin{thm} \label{thm4}
A finitely-generated projective Hopf algebroid $\mathcal{H}=\left(\mathcal{H}_{L},\mathcal{H}_{R},S\right)$ over $C(X)$ in which the images of $C(X)$ under the source and target maps are central, is a bundle of $\mathbb{C}$-Hopf algebroids $\mathcal{H}_{x}=\left(H^{L}_{x},H^{R}_{x},S^{\wedge}_{x}\right)$ over $X$.
\end{thm}

Since $(B,\mathcal{H})$ is a covering of $A$, $B$ is an $\mathcal{H}$-Galois extension of $A$ which means that $B$ comes with right coactions $B\stackrel{\rho_{R}}{\longrightarrow}B\tens{A}\mathcal{H}_{R}$ and $B\stackrel{\rho_{L}}{\longrightarrow}B\tens{A}\mathcal{H}_{L}$ by $\mathcal{H}$ whose common coinvariant is $A$. Note that both coactions $\rho_{R}$ and $\rho_{L}$ are $A$-module maps. Thus, applying the functor $Sh$ gives vector bundles maps $ \xymatrix{E \ar[rr]^-{Sh(\rho_{R})} && E\otimes H^{R} }$ and $\xymatrix{E \ar[rr]^-{Sh(\rho_{L})} && E\otimes H^{L} }$. Each of these bundle maps induce coactions $\rho_{L,x}$ and $\rho_{R,x}$ of the fiber Hopf algebroids $\mathcal{H}_{x}$ of $\mathcal{H}$ to the fiber algebras $E_{x}$ of $B$ by the commutativity of the diagrams below for $T=R,L$.

\[ \xymatrix{
& && E \ar[rr]^-{Sh(\rho_{T})} \ar'[d][dd]^-{Sh(\rho_{T})} && \makebox[\widthof{$E$}][l]{$E\otimes H^{T}$} \ar[dd]^-{id\otimes Sh(\Delta_{T})} \\
& E\otimes \mathbb{C}_{triv} \ar@{=}[rru] & E_{x} \ar[ru] \ar[rr]^(.65){\rho_{T,x}} \ar[dd]_-{\rho_{T,x}} && E_{x}\otimes H^{T}_{x} \ar[ru] \ar[dd]^(.25){id\otimes \delta_{T,x}} & \\
E_{x}\otimes\mathbb{C} \ar[ru] \ar@{=}[rru] & && E\otimes H^{T} \ar'[r]_-{Sh(\rho_{T})\otimes id}[rr] \ar[ull]|(0.49)\hole_(0.35){id\otimes\epsilon_{T}}|(0.75)\hole && \makebox[\widthof{$E$}][l]{$E\otimes H^{T}\otimes H^{T}$} \\
&& E_{x}\otimes H^{T}_{x} \ar[ru] \ar[rr]_-{\rho_{T,x}\otimes id} \ar[ull]^-{id\otimes\epsilon_{T,x}} && E_{x}\otimes H^{T}_{x}\otimes H^{T}_{x} \ar[ru] & } \]

The commutativity of the following diagram shows that

\[ \xymatrix{&  E\otimes E \ar[rr]^-{Sh(\rho_{T})\otimes Sh(\rho_{T})} \ar'[d][dd]_-{m} && \makebox[\widthof{$E$}][l]{$E\otimes H^{T}\otimes E \otimes H^{T}$} \ar[dd]|-{(m\otimes m)\circ \reflectbox{F}} \\
E_{x}\otimes E_{x} \ar[ru] \ar[rr]^(.65){\rho_{T,x}\otimes \rho_{T,x}} \ar[dd]_-{m} && \left(E_{x}\otimes H^{T}_{x}\right)^{\otimes 2} \ar[ru] \ar[dd]^(.25){(m\otimes m)\circ \reflectbox{F}} & \\
 & E \ar'[r]_-{Sh(\rho_{T})}[rr] && \makebox[\widthof{$E$}][l]{$E\otimes H^{T}$} \\
 E_{x} \ar[ru] \ar[rr]_-{\rho_{T,x}} && E_{x}\otimes H^{T}_{x} \ar[ru] & }\]

\noindent shows that $\rho_{T,x}$ for $T=L,R$ is multiplicative.

The coinvariants $A$ of the coaction $\rho_{R}$ is the equalizer of $\rho_{R}$ and $id\otimes_{A} 1$. Similarly, $A$ is the equalizer of $\rho_{L}$ and $id\otimes_{A}1$, i.e. we have the following diagrams of $A$-modules.

\[ \xymatrix{A \ar[rr] & & B \ar@<1ex>[rr]^-{\rho_{R}} \ar@<-1ex>[rr]_-{id\tens{A}1} & & B\tens{A}\mathcal{H}_{R}} \]
\[ \xymatrix{A \ar[rr] & & B \ar@<1ex>[rr]^-{\rho_{L}} \ar@<-1ex>[rr]_-{id\tens{A}1} & & B\tens{A}\mathcal{H}_{L}}\]

\noindent Applying the functor $Sh$ to the first diagram gives us the following

\[ \xymatrix{\mathbb{C}_{triv} \ar[rr] & & E \ar@<1ex>[rr]^-{Sh(\rho_{R})} \ar@<-1ex>[rr]_-{id\otimes1} & & E\otimes H^{R}} \]

\noindent which is an equalizer diagram as well since $Sh$ is an equivalence. Thus, the coinvariant of the induced coaction $\rho_{R,x}$ is $\mathbb{C}$. Similarly, $\mathbb{C}$ is the coinvariant of the induced coaction $\rho_{L,x}$.

\noindent Now, let us show that associated Hopf-Galois map $\mathfrak{gal}_{R,x}$ to $\rho_{R,x}$ is a bijective. The $A$-module isomorphism $B\otimes_{A}B\stackrel{\mathfrak{gal_{R}}}{\longrightarrow} B\otimes_{A} \mathcal{H}_{R}$ induces a bundle isomorphism

\[ \xymatrix{E\otimes E \ar[rr]^-{Sh(\mathfrak{gal_{R}})} & & E\otimes H^{R} }\]

\noindent which on fibers give the isomorphism

\[ \xymatrix{E_{x}\otimes E_{x} \ar[r]^-{\mathfrak{gal}_{x}} & E_{x}\otimes H^{R}_{x}.}\]

Similarly, the associated Hopf-Galois maps $\mathfrak{gal}_{L,x}$ to the coactions induced on the fibers by $\rho_{L}$ are all bijective. These give the following result.

\begin{thm} \label{thm5}
Let $(B,\mathcal{H})$ be a local covering of $A=C(X)$ in which $A$ is central. With the notation as above, $(E_{x},\mathcal{H}_{x})$ is a covering of the point $x$.
\end{thm}

\

Using proposition \ref{prop7} and the previous theorem, we get the following corollary.

\

\begin{cor} \label{cor1}
Let $(B,\mathcal{H})$ be a local covering of $C(X)$ in which $C(X)$ is central. Then $\mathcal{H}$ gives two bundles $\xymatrix{H^{1},H^{2} \ar@{->>}[r] & X}$ of coupled Hopf algebras over $X$.
\end{cor}

\

\begin{exa}\label{exa6}
Note that the fiber coverings $(E_{x},\mathcal{H}_{x})$ need not be isomorphic even within a connected component of $X$. As a matter of fact, we already have an example for this in the commutative case. Consider the algebras $E_{t}=\mathbb{C}[x]/(x^{n}-t)$. The underlying vector space of these algebras are all $n$-dimensional and they constitute a vector bundle $\xymatrix{E \ar@{->>}[r]^-{p} & \mathbb{C} }$ over the complex plane where $p^{-1}(t)=E_{t}$. Note that the each fiber carries a natural algebra structure making $E$ an algebra bundle over $\mathbb{C}$ with non-isomorphic fibers. In particular, the fiber algebra $E_{0}$ has a nilpotent element while $E_{1}$ has none. Furthermore, each fiber algebra is spanned by $\left\{1,x,...,x^{n-1}\right\}$. The group $G=\mathbb{Z}/n\mathbb{Z}$ acts on each fiber algebra $E_{t}$ via $(m\cdot x) \mapsto \lambda^{m}x$ extended into an algebra isomorphism where $\lambda$ is a primitive $n^{th}$ root of $1$. This action extends to a Galois action of the group algebra $\mathbb{C}G$ and hence, the function algebra $C(G)$ coacts on $B=\Gamma(\mathbb{C},E)$. This turns $(B,C(\mathbb{C})\otimes C(G))$ into a local covering of $C(\mathbb{C})$.
\end{exa}

\subsection{Coverings with semisimple fibers}\label{S4.4}

In this section, we will continue to look at the case when $A=C(X)$ is central in the local covering $(B,\mathcal{H})$ and whose fibers are semisimple. For simplicity, let us also assume that $X$ is connected. This means that any vector bundle $E$ for which $B=\Gamma(X,E)$ and any vector bundle $H$ for which $\mathcal{H}=\Gamma(X,H)$, the underlying complex algebras of the fiber algebras $E_{x}$ and the fiber Hopf algebroids $(H_{x}^{L},H_{x}^{R},S_{x}^{\vee})$ are semisimple complex algebras. By Wedderburn's theorem, $E_{x}$ is the finite product of matrix algebras i.e.,

\[ E_{x}=M_{n_{1}}(\mathbb{C})\times M_{n_{2}}(\mathbb{C}) \times \cdots \times M_{n_{j}}(\mathbb{C}) \]

\noindent for some positive integers $n_{1},n_{2},...,n_{j}$. This decomposition determines (and is completely determined by) a set of central orthogonal idempotent $\left\{e_{i}\in E_{x}|i=1,2,...,j\right\}$ summing up to 1. Explicitly, $M_{n_{i}}(\mathbb{C})\cong e_{i}E_{x}$ for all $i=1,...,j$. Let us call the (unordered) $j$-tuple $(n_{1},n_{2},...,n_{j})$ the \textit{Wedderburn shape} of the semisimple algebra $E_{x}$. Part of the content of Wedderburn's theorem says that the Wedderburn shape of a semisimple algebra is unique.

\

\begin{exa} \label{exa7}
Let us consider the extreme case when for all $x\in X$, $E_{x}\cong\mathbb{C}^{n}$ for some $n\in \mathbb{N}$ as algebras. In this case, the Hopf algebroid $\mathcal{H}$ is commutative by the bijectivity of the associated Hopf-Galois maps. By assumption, the antipode $S$ is bijective. Assuming the coproduct and the counit are unital maps, lemma \ref{L1} implies that there is a groupoid $\mathcal{G}$ such that $\mathcal{H}\cong C(\mathcal{G})$.

Bijectivity of $\mathfrak{gal}_{x}$ above implies that the underlying $\mathbb{C}$-vector space of $H_{x}$ is finite-dimensional for any $x\in X$. Specifically, each $H_{x}$ is of dimension $n$. Now, given $x\in X$ consider the following diagram in the category of topological spaces

\[ \xymatrix{ \mathcal{G}(x) \ar[rr] \ar[dd] \ar@{}[rrdd]|(0.3){\big\lrcorner} & & Eq(s,t) \ar[dd] \ar[rr] & & \mathcal{G} \ar@<1ex>[dd]^-{t} \ar@<-1ex>[dd]_-{s} \\
& & & & \\
x \ar[rr] & & X \ar@{=}[rr] & & X \\ } \]

\noindent where the left square is a pull-back and the right square is an equalizer diagram. Applying the functor $C(-)$ gives the following diagram

\[ \xymatrix{C(\mathcal{G}(x)) \ar@{}[rrdd]|(0.3){\big\ulcorner} & & C(Eq(s,t)) \ar[ll] & & \mathcal{H} \ar[ll] \\
& & & & \\
\mathbb{C} \ar[uu] & & A \ar@{->>}[ll]_-{ev_{x}} \ar[uu] \ar@{=}[rr] & & A \ar@<1ex>[uu]^-{t} \ar@<-1ex>[uu]_-{s} \\ } \]

\noindent where the left square is a push-out diagram. The right square being a coequalizer implies that the large rectangular diagram (using either $s$ or $t$) is a push-out diagram. The counit of the adjunction $C(-) \dashv Spec$ provides a $\mathbb{C}$-algebra isomorphism $C(\mathcal{G}(x))\cong H_{x}$. This extends to a bialgebroid isomorphism since the coring structure maps of $H_{x}$ and $C(\mathcal{G}(x))$ are morphisms of commutative unital $\mathbb{C}$-algebras. Since $\mathcal{G}(x)$ is a group, $H_{x}$ is then a Hopf algebra. Note that a priori, $G_{x}$ depends on $x\in X$ but connectivity of $X$ implies that the groups $G_{x}$ are all isomorphic, denoted accordingly as $G$. When dualized, the coaction $\rho_{x}:\mathbb{C}^{n}\longrightarrow\mathbb{C}^{n}\otimes H_{x}$ gives an action $\mathbb{C}^{n}\otimes \mathbb{C}G\longrightarrow \mathbb{C}^{n}$. Note that $\mathbb{C}^{n}\otimes \mathbb{C}G\cong \mathbb{C}[Y]\otimes \mathbb{C}G\cong\mathbb{C}[Y\times G]$ where $Y$ is a set consisting of $n$ points and the multiplication in the algebra $\mathbb{C}[Y\times G]$ is pointwise in $Y$ but convolution in $G$. The map $\rho_{x}^{*}$ is completely determines by the map $Y\times G\stackrel{\alpha}{\longrightarrow}Y$ which is an action by the virtue of $\rho_{x}$ being a coaction. The bijectivity of the Hopf-Galois map translate to the bijectivity of the associated map

\[ Y\times G\longrightarrow Y\times Y, \hspace{.25in} (y,g)\mapsto (yg,y) \]

\noindent which means that the action $\alpha$ is free and transitive. Thus, $G\leqslant S_{n}$ is a transitive subgroup with $|G|=n$. $\Box$
\end{exa}

\

Let us consider the general case when the fibers of $E$ are non-commutative algebras. In this case, $E_{x}=M_{n_{1}}(\mathbb{C})\times M_{n_{2}}(\mathbb{C}) \times \cdots \times M_{n_{j}}(\mathbb{C})$ where the Wedderburn shape $(n_{1},n_{2},\cdots,n_{j})$ of $E_{x}$ a priori depends on $x\in X$. Let us describe how these dependence works.

Consider the center $Z(B)$ of $B$. Since $B=\Gamma(X,E)$ equipped with pointwise multiplication, we see that $\sigma\in Z(B)$ precisely when $\sigma(x)\in Z(E_{x})$ for all $x\in X$. The center $Z(B)$ is a $C^{*}$-subalgebra of $B$. In particular, it is a commutative $C^{*}$-algebra and by the Gelfand duality, there is a compact Hausdorff space $Y$ such that $Z(B)=C(Y)$. Note that $A=C(X)$ sits inside $Z(B)=C(Y)$. Thus, there is a continuous surjective map $\xymatrix{Y \ar@{->>}[r]^-{p} & X}$. Consider the following stratification of $X$. Denote by $X^{(n)}=\left\{x\in X| \#(p^{-1}(x))=n\right\}$ where $\#(S)$ denotes the cardinality of the set $S$. Note that $X^{(n)}$'s are generally not connected. Define $X^{(n,i)}$, $i\in I_{n}$ to be the connected components of $X^{(n)}$. Note that the $X^{(n,i)}$'s forms a partition of $X$ and that the $X^{(n,i)}$'s are generally not closed in $X$. We call $\left\{X^{(n,i)}|n\in \mathbb{N}, i\in I_{n}\right\}$ the stratification of $X$ and each $X^{n,i}$ as a stratum. Let us denote by $Y^{(n,i)}=p^{-1}(X^{(n,i)})$. Then $\xymatrix{Y^{(n,i)} \ar@{->>}[r]^-{p} & X^{(n,i)}}$ is a covering space in the classical sense.

Surjectivity of $p$ implies that $X^{(0)}=\emptyset$. We claim that $X^{(n,i)}=\emptyset$ as well for $n\geqslant m$ for some sufficiently large $m$. To see this, note that semisimplicity of $E_{x}$ implies that $Z(E_{x})\subseteq E_{x}$ is complemented. This implies that the dimension of $Z(E_{x})$ is bounded above by the dimension of $E_{x}$. By theorem \ref{thm3}, we see that this dimension is bounded by $dim \ H <\infty$. The center $Z(E_{x})$ of $E_{x}$ is linearly generated by the central orthogonal idempotents $\left\{e_{i}\right\}$ giving the Wedderburn factors. These central orthogonal idempotents can be extended continuously to relative sections $\left\{\sigma_{i}\in\Gamma(X^{(n,j)},Z(E))|x\in X^{(n,j)},\sigma_{i}(x)=e_{i}\right\}$. Since the rank of an idempotent is locally constant, we see that Wedderburn factors are all the same for all $x\in X^{(n,j)}$. Thus, we see that Wedderburn shape of the fibers $E_{x}$ of $E$ only depend on the stratum of $x\in X$.

On the other hand, much can be said about the fiber Hopf algebroids. From section \ref{S4.1} such a Hopf algebroid is a coupled Hopf algebra. There are only finitely many semisimple complex Hopf algebras of a given fixed dimension. Thus, there are only finitely many coupled Hopf algebras of a given dimension. Since the fiber Hopf algebroids have the same dimension, this implies that there are only finitely many posibilities for their structure. Connectivity of $X$ and discreteness of the collection of such coupled Hopf algebras imply that the fiber Hopf algebroids must be isomorphic, say to a fixed one $H_{0}=(H_{x_{0}}^{L},H_{x_{0}}^{R},S_{x_{0}}^{\vee})$.

\begin{prop}\label{prop8}
For any $x,y\in X$, $H_{x}\cong H_{y}$ as coupled Hopf algebras.
\end{prop}

Specializing the notion of an algebraic morphism of Hopf algebroids from section \ref{S2.1}, tells us that a morphism $(H^{L}_{1},H^{R}_{1},S_{1})\stackrel{\phi}{\longrightarrow}(H^{L}_{2},H^{R}_{2},S_{2})$ of coupled Hopf algebras is a linear map $\phi$ which defines Hopf algebra maps $H^{L}_{1}\stackrel{\phi}{\longrightarrow}H^{L}_{2}$ and $H^{R}_{1}\stackrel{\phi}{\longrightarrow}H^{R}_{2}$ intertwining the coupling maps. This makes sense since $H^{L}_{1}$ and $H^{R}_{1}$ have the same underlying algebra. Same goes for $(H^{L}_{2},H^{R}_{2},S_{2})$.

Let $G=Aut(H_{0})$ and let $\phi\in G$. Finite dimensionality of $H^{L}_{x_{0}}$ and $H^{R}_{x_{0}}$ implies that they are Frobenius algebras. Thus, they are equipped with nondegenerate pairings $\left\langle ,\right\rangle_{L}$ and $\left\langle ,\right\rangle_{R}$ making them finite-dimensional Hilbert spaces. The automorphism $\phi$ in particular defines automorphisms of these two Frobenius algebras, i.e. $\phi$ preserves the inner products $\left\langle ,\right\rangle_{L}$ and $\left\langle ,\right\rangle_{R}$. Thus, each $\phi\in G$ is a unitary map with respect to both inner products (actually, since there is a unique Hilbert space up to isomorphism for a particular dimension, the two inner product defines the same Hilbert space structure on $H^{L}_{x_{0}}$) and $H^{R}_{x_{0}}$). Hence, we have the following proposition.

\begin{prop}\label{prop9}
$G\subseteq U(n)$ where $n=dim \ H^{L}_{x_{0}}$.
\end{prop}

The two propositions give a continuous map $\alpha:X\longrightarrow G$, $\alpha(x):H_{x}\stackrel{\simeq}{\longrightarrow}H_{x_{0}}$. By Radford \cite{rad}, the group of automorphisms of a semisimple Hopf algebra over a field of characteristic 0 is finite. Hence, the group of automorphisms of a semisimple coupled Hopf algebra over $\mathbb{C}$ is finite. This implies that $G$ is a finite subgroup of $U(n)$ and thus, finite. Hence, $\alpha$ is a $\breve{C}$ech 1-cocycle since it is locally constant. Therefore, $H\twoheadrightarrow X$ is an algebra bundle, i.e. the local transition maps rather than just being linear maps, are algebra maps. The associated $\breve{C}$ech 1-cocycle is just $\alpha$ followed by the inclusion $G\subseteq GL_{n}(\mathbb{C})$.

\begin{prop}\label{prop10}
$G\subseteq GL_{n}(\mathbb{C})$ is finite and $H\twoheadrightarrow X$ is an algebra bundle.
\end{prop}

\

As we have argued after example \ref{exa7}, the fibers algebras need not be isomorphic. Let us discuss a particular instance when the fiber algebras are all isomorphic. Let $X$ be a compact connected smooth manifold. Let $A=C^{\infty}(X)$ and let $(B,\mathcal{H})$ be a local central covering of $A$. By Serre-Swan, $B=\Gamma^{\infty}(X,E)$ for some finite-rank smooth vector bundle $E\twoheadrightarrow X$. By a \textit{differential connection} $\nabla$ on $E$ we mean a connection $\nabla$ such that for any vector field $\nu$ on $X$ we have

\[ \nabla_{\nu}(\sigma_{1}\sigma_{2})=\sigma_{1}\nabla_{\nu}(\sigma_{2})+\nabla_{\nu}(\sigma_{1})\sigma_{2} \]

\noindent for any sections $\sigma_{1},\sigma_{2}\in B$. We have the following proposition.

\begin{prop}\label{prop11}
If $E$ has a differential connection $\nabla$ then the fiber algebras of $E\twoheadrightarrow X$ are all isomorphic.
\end{prop}

\begin{prf}
Let $x,y\in X$ and let $\gamma:I\longrightarrow X$ be a (piecewise) smooth path in $X$ with $\gamma(0)=x$ and $\gamma(1)=y$. Using the connection $\nabla$, we have a parallel transport map

\[ \Phi(\gamma)^{y}_{x}:E_{x}\longrightarrow E_{y} \]

\noindent which is a linear isomorphism. Thus, all we have to show is that $\Phi(\gamma)^{y}_{x}$ is multiplicative. Given $b_{1},b_{2}\in E_{x}$, there are unique smooth sections $\sigma_{1}$ and $\sigma_{2}$ of $E$ along $\gamma$ such that $\nabla_{\harpoon \gamma}\sigma_{1}=\nabla_{\harpoon \gamma}\sigma_{2}=0$ and $\sigma_{1}(x)=b_{1}$ and $\sigma_{2}(x)=b_{2}$. Here, $\harpoon \gamma$ denotes the smooth tangent vector field of $\gamma$. Note that the product $\sigma_{1}\sigma_{2}$ is the unique smooth section of $\xymatrix{E \ar@{->>}[r] & X}$ along $\gamma$ such that $\left(\sigma_{1}\sigma_{2}\right)(x)=\sigma_{1}(x)\sigma_{2}(x)=b_{1}b_{2}$ and

\[ \nabla_{\harpoon \gamma}(\sigma_{1}\sigma_{2})= \sigma_{1}\nabla_{\harpoon \gamma}(\sigma_{2})+\nabla_{\harpoon \gamma}(\sigma_{1})\sigma_{2} =0. \]

\noindent Thus, by definition of the parallel transport map $\Phi\left(\gamma\right)_{x}^{y}$ we have

\[ \Phi\left(\gamma\right)_{x}^{y}(b_{1}b_{2})=\left(\sigma_{1}\sigma_{2}\right)(y)=\sigma_{1}(y)\sigma_{2}(y)=\Phi\left(\gamma\right)_{x}^{y}(b_{1})\Phi\left(\gamma\right)_{x}^{y}(b_{2}) \]

\noindent which shows that $\Phi\left(\gamma\right)_{x}^{y}$ is multiplicative. $\blacksquare$
\end{prf}

A strong converse of the above proposition, where the isomorphisms among fibers satisfy some coherence conditions, holds. By a coherent collection

\[ \mathscr{P}=\left\{\Phi(\gamma)_{x}^{y}:E_{x}\longrightarrow E_{y}|\forall x,y\in X, \gamma:I\longrightarrow X \ smooth\right\} \]

\noindent of isomorphisms among fibers of $E\twoheadrightarrow X$, we mean a collection satisfying

\begin{enumerate}
\item[(i)] $\Phi(\gamma)^{x}_{x}=id$,
\item[(ii)] $\Phi(\gamma)^{y}_{u}\circ\Phi(\gamma)^{u}_{x}=\Phi(\gamma)^{y}_{x}$,
\item[(iii)] and $\Phi$ depends smoothly on $\gamma$, $y$ and $x$.
\end{enumerate}

\noindent We then have the following proposition.

\begin{prop}\label{prop12}
A coherent collection $\mathscr{P}$ of algebra isomorphisms on $E\twoheadrightarrow X$ gives a differential connection $\nabla$ on $E$.
\end{prop}

\begin{prf}
Using the collection $\mathscr{P}$ we can immediately write an infinitessimal connection $\nabla$ as follows: for any vector $V$ on $X$ we have

\[ \nabla_{V}(\sigma)=\lim\limits_{t\rightarrow 0} \dfrac{\Phi(\gamma)^{x}_{\gamma(t)}\sigma(\gamma(t))-\sigma(x)}{t}=\left.\dfrac{d}{dt}\right|_{t=0}\Phi(\gamma)^{x}_{\gamma(t)}\sigma(\gamma(t)) \]

\noindent for any $\sigma\in B$ and $x=\gamma(0)$. That $\nabla$ is a differential connection follows from the multiplicativity of $\Phi(\gamma)^{y}_{x}$ and the Leibniz property of $\left.\dfrac{d}{dt}\right|_{t=0}$. $\blacksquare$
\end{prf}

\

\begin{exa} \label{exa8}
In this example, we will show that the Wedderburn shape of fibers need not be constant even over a connected base space. Let $G$ be a finite group of \textit{central type}, i.e. $G$ fits in an extension

\[ \xymatrix{1 \ar[r] & Z(\Gamma) \ar[r] & \Gamma \ar[r] & G \ar[r] & 1} \]

\noindent such that $\Gamma$ has an irreducible representation $\Gamma\stackrel{\rho}{\longrightarrow}GL(V)$ of dimension $\sqrt{[\Gamma:Z(\Gamma)]}$.

Now, the group extension above determines a $2$-cocycle $\beta:G\times G\longrightarrow Z(\Gamma)$. Then the composition

\[ \xymatrix{ G\times G \ar[rr]^-{\beta} & & Z(\Gamma) \ar[rr]^-{\rho} \ar@{-->}[rd] & & GL(V) \\
& & & \mathbb{C}^{\times} \ar[ru] & \\ }\]

\noindent determines a $2$-cocycle $\alpha$ such that the associated twisted group algebra $\mathbb{C}^{\alpha}G\cong M_{n}(\mathbb{C})$, where $n=\sqrt{[\Gamma:Z(\Gamma)]}$. The twisted group algerba $\mathbb{C}^{\alpha}G$ is a Hopf algebra with the same coproduct, counit and unit as that of $\mathbb{C}G$ with product given by $g\cdot g^{'}=\alpha(g,g^{'})gg^{'}$ for any $g,g^{'}\in G$. Such a cocycle can be rescaled to get a family of cocycles $\alpha_{t}$ for every $t\in\mathbb{C}$ with $\alpha_{0}=1$ and $\alpha_{t}$ nondegenerate for $t\neq 0$. This means $\mathbb{C}^{\alpha_{t}}G\cong M_{n}(\mathbb{C})$ for $t\neq 0$ while $\mathbb{C}G$ may decompose nontrivially as a direct sum of matrix algebras over $\mathbb{C}$. This gives a bundle of Hopf algebras $E=\coprod_{t\in\mathbb{C}}\mathbb{C}^{\alpha_{t}}G\stackrel{p}{\longrightarrow}\mathbb{C}$. The algebra $B=\Gamma(\mathbb{C},E)$ is then a Hopf-Galois extension of $C(\mathbb{C})$. $\Box$
\end{exa}

\

\subsection{Coverings with cleft fibers}\label{S4.5}

In this section, we are still interested with the case $A=C(X)$ and $(B,\mathcal{H})$ is a local covering in which $A$ is central. As before, $B\cong\Gamma(X,E)$ and $\mathcal{H}\cong\Gamma(X,F)$ where $E$ and $F$ is an algebra bundle and a Hopf algebroid bundle both over $X$, respectively. Moreover, for any $x\in X$, $(E_{x},F_{x})$ is a covering of $\mathbb{C}$. In addition, suppose that $(B,\mathcal{H})$ is a \textit{cleft} covering i.e., $A\subseteq B$ is a cleft extension. Recall from theorem \ref{thm2} that this implies that $B\cong A\otimes_{A}H$ as left $A$-modules and as right $\mathcal{H}$-comodules. These conditions descend to the bundle structures of $E$ and $F$, i.e. $E_{x}\cong \mathbb{C}\otimes F_{x}$ as left $\mathbb{C}$-modules and as right $F_{x}$-comodules. Since $(E_{x},F_{x})$ is a covering of $\mathbb{C}$, again by theorem \ref{thm2} we see that $(E_{x},F_{x})$ is a cleft covering of $\mathbb{C}$. In other words, cleft central coverings of commutative spaces have cleft coverings as fibers.

\begin{thm} \label{thm6}
With the assumption of this section, $A\stackrel{\mathcal{H}}{\Longrightarrow}B$ is a cleft covering implies that the fiber coverings are also cleft.
\end{thm}

%%%%%%%%%%%%%%%%%%%%%%%%%%%%%%%%%%%%%%%%%%%%%%%%%%%%%%%%%%%%%%%%%%%%%%%%%%%%%%%%%%%%%%%%

\section{Coverings of the noncommutative torus}\label{S5.0}

In section \ref{S5.0}, we dealt with the general situation of local coverings $(B,\mathcal{H})$ of a commutative space $A=C(X)$ such that $A$ is noncentral in $(B,\mathcal{H})$. In this section, we will see a particular example of such coverings. Though the algebraic structures are no longer pointwise, they have a nice description for rational and commutative tori as we will see in the following section.

\

\subsection{Commutative and rational noncommutative tori}\label{S5.1}

Let $q\in\mathbb{C}$ be a primitive $n^{th}$ root of unity. Let $B$ be the universal $C^{*}$-algebra generated by unitaries $U$ and $V$ satisfying $UV=qVU$. Let $A$ be the $C^{*}$-subalgebra generated by $U$ and $V^{n}$. Then, $A$ is the universal $C^{*}$-algebra generated by commuting unitaries $U$ and $V^{n}$ and hence $C^{*}$-isomorphic to the continuous functions on the $2$-torus, i.e. $A\cong C(\mathbb{T}^{2})$. As an $A$-module, $B$ is finitely-generated and projective generated by $\left\{1,V,...,V^{n-1}\right\}$. Thus, by the Serre-Swan theorem $B\cong\Gamma(\mathbb{T}^{2},\mathbb{E})$ for some finite-rank vector bundle $\mathbb{E}$ over $\mathbb{T}^{2}$. However, the multiplication in $B$ is not the pointwise multiplication on $\Gamma(\mathbb{T}^{2},\mathbb{E})$ since $A$ is not central in $B$. Let us describe the product in $B$ as an $A$-ring. Since $B$ is free over $A$ via the isomorphism

\[ B \cong \bigoplus\limits_{i=0}^{n-1} A\cdot V^{i}. \]

\noindent Let us index the generating set of $B$ as an $A$-module by $\mathbb{Z}/n$, the group of integers modulo $n$. Elements $f$ and $g$ of $B$ are of the form

\[ f=\sum\limits_{i\in\mathbb{Z}/n} a_{i}V^{i}, \hspace{.5in} g=\sum\limits_{i\in\mathbb{Z}/n} b_{i}V^{i} \]

\noindent for some $a_{i},b_{i}\in A, i=0,...,n-1$.

\begin{figure}[ht!]
\centering
\includegraphics[width=135mm]{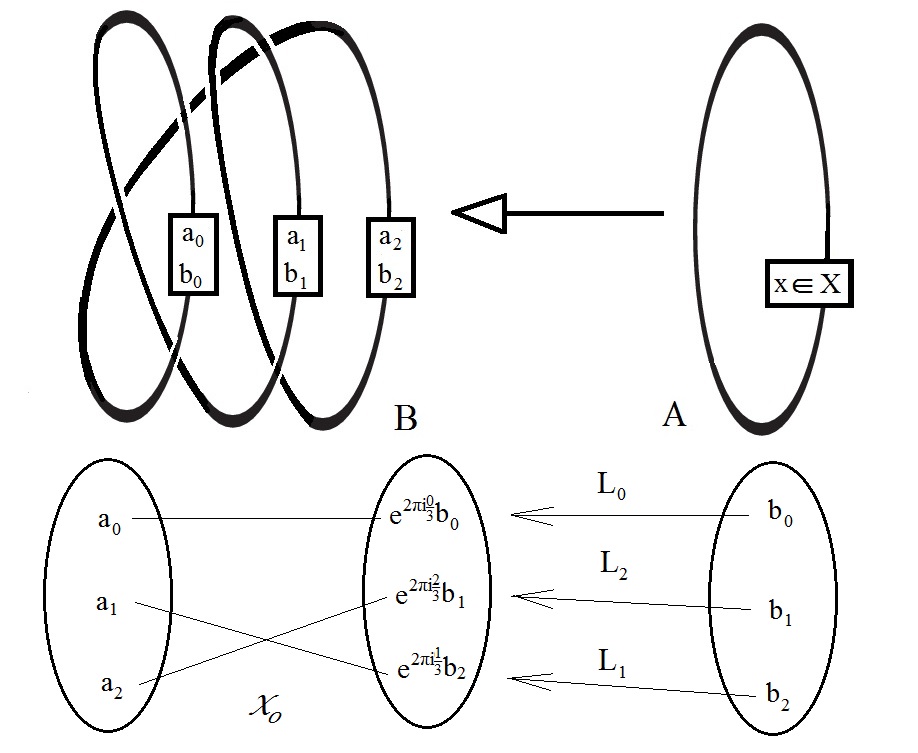}
\caption{Convolution-pointwise product}
\end{figure}

\noindent Then

\[ fg= \sum\limits_{k\in\mathbb{Z}/n} \chi_{k}\left(\alpha,\beta\right)V^{k} \]

\noindent for $\chi_{k}\left(\alpha,\beta\right)\in A, \ k=0,...,n-1$ where $\alpha=(a_{0}, a_{1}, ...,a_{n-1})$ and $\beta=(b_{0}, b_{1}, ...,b_{n-1})$. Let us describe $\chi_{k}$. Denote by $L:A\longrightarrow A$ the diagonal operator defined on linear generators of $A$ by $L(U^{x}V^{ny})=q^{-x}U^{x}V^{ny}$. Consider the group table of $\mathbb{Z}/n$ considered as a matrix, denoted as $\Omega$. Change those entries different from $k\in\mathbb{Z}/n$ to $0$ and change the entries with $k$ to $L^{i-1}$ if that entry is in the $i^{th}$ row. Denote this operator matrix by $\Omega_{k}$. Then

\[ \chi_{k}\left(\alpha,\beta\right) = \alpha \Omega_{k} \beta^{T} = \left(\begin{array}{cccc}
a_{0}, & a_{1}, & ..., & a_{n-1}\\
\end{array}\right)\Omega_{k}\left(\begin{array}{c}
b_{0} \\
b_{1} \\
\vdots \\
b_{n-1}
\end{array}\right)=\sum\limits_{i=0}^{n-1}a_{i}L^{i}(b_{k-i}). \]

\noindent for $k=0,...,n-1$. As an example, for $n=3$ we have

\[  \Omega_{0}=\left(\begin{array}{ccc}
L^{0} &  &  \\
 &  & L^{1} \\
 & L^{2} &  \\
\end{array}\right), \hspace{.3in} \Omega_{1}=\left(\begin{array}{ccc}
 & L^{0} &  \\
L^{1} &  &  \\
 &  & L^{2} \\
\end{array}\right), \hspace{.3in} \Omega_{2}=\left(\begin{array}{ccc}
 &  & L^{0} \\
 & L^{1} &  \\
L^{2} &  &  \\
\end{array}\right) \]

\noindent and so

\[ \chi_{0}=a_{0}L^{0}(b_{0})+a_{1}L^{1}(b_{2})+a_{2}L^{2}(b_{1}) \]
\[ \chi_{1}=a_{0}L^{0}(b_{1})+a_{1}L^{1}(b_{0})+a_{2}L^{2}(b_{2}) \]
\[ \chi_{2}=a_{0}L^{0}(b_{2})+a_{1}L^{1}(b_{1})+a_{2}L^{2}(b_{0}).\]

\

The $A$-ring structure of $B$ is pointwise-convolution as illustrated in figure 2. Denote by $\mathcal{H}=C(G,A)$, where $G=\mathbb{Z}/n$. We claim that $\mathcal{H}$ is a commutative Hopf algebroid. The left- and right-bialgebroid structures of $\mathcal{H}$ are isomorphic, with pointwise product, whose source, target, counit and antipode map is

\[ \xymatrix@R=0.2cm{ A \ar[rr]^-{s,t} & & \mathcal{H},\\
1 \ar@{|->}[rr] & & 1 } \hspace{.5in} \xymatrix@R=0.2cm{ \mathcal{H} \ar[rr]^-{\epsilon} & & A,\\
f \ar@{|->}[rr] & & f(1) } \hspace{.5in} \xymatrix@R=0.2cm{ \mathcal{H} \ar[rr]^-{S} & & \mathcal{H},\\
f \ar@{|->}[r] & Sf, & Sf(x)=f(x^{-1}) }\]

\noindent respectively, and whose coproduct is

\[ \xymatrix@R=0.2cm{ \mathcal{H} \ar[rr]^-{\Delta} & & \mathcal{H}\tens{A}\mathcal{H} \cong C(G\times G,A)\\
f \ar@{|->}[r] & \Delta f, & \Delta f(x,y)=f(xy). } \]

The group $G$ acts on $B$ as follows: $g\cdot U=U$, $g\cdot V=q V$ where $g\in G$ is a generator. This action extends to a module structure over the group algebra $\mathcal{H}^{*}=AG$, the $A$-dual of the Hopf algebroid $\mathcal{H}$. The $\mathcal{H}^{*}$-invariants of $B$ is $A$. Thus, $B$ carries a coaction of $\mathcal{H}$ whose coinvariants is $A$. It is immediate to check that this defines a local covering $(B,\mathcal{H})$ of $A$. $\Box$

\

\begin{rem}
\begin{enumerate}
\item[]

\item[(i)] The covering $(B,\mathcal{H})$ of $A$ above is an example of a covering where $A$ is a commutative space which is not central in $B$. However, the images of $A$ under the source and target map is central in $\mathcal{H}$ as it is a commutative Hopf algebroid. This implies that $\mathcal{H}$ is a bundle of Hopf algebroids (actually, of Hopf algebras) but the coaction is not pointwise.

\item[(ii)] We can generalize the example above as follows. Given integers $n$ and $m$, let $q$ be a primitive $nm^{th}$ root of unity. Let $B$ be the universal $C^{*}$-algebra generated by unitaries $U$ and $V$ satisfying $UV=qVU)$ and let $A$ be the $C^{*}$-subalgebra generated by commuting unitaries $U^{n}$ and $V^{m}$. Thus, $A\cong C(\mathbb{T}^{2})$. Take $\mathcal{H}$ to be the commutative Hopf algebroid $C(G,A)$ over $A$ where $G=\mathbb{Z}/n\times \mathbb{Z}/m$. As a matter of fact, we can construct a coverings of $C(\mathbb{T}^{2})$ for any finite quotient $G$ of $\mathbb{Z}^{2}$. We outline this construction in the next section.

\end{enumerate}
\end{rem}

\

Let $\theta=\frac{n}{m}\in\mathbb{Q}$. The center of the noncommutative torus $\mathbb{T}^{2}_{\theta}$ is the $C^{*}$-subalgebra generated by $U^{m}$ and $V^{m}$. The computation above implies that rational noncommutative tori give local coverings of the commutative torus with commutative quantum symmetries. Thus, we get the following proposition.

\

\begin{prop}\label{prop13}
Let $\theta=\frac{n}{m}$ for coprime integers $n$ and $m$ with $m>0$. Let $\mathbb{T}^{2}_{\theta}$ be the noncommutative torus with parameter $\theta$. Then there is a commutative Hopf algebroid $\mathcal{H}$ such that $(\mathbb{T}^{2}_{\theta},\mathcal{H})$ is a covering of $Z(\mathbb{T}^{2}_{\theta})=C(\mathbb{T}^{2})$.
\end{prop}

\

We have an explicit presentation of $\mathbb{T}^{2}_{\theta}$ as a bundle over $\mathbb{T}^{2}$. Consider the following elements of $\mathbb{T}^{2}_{\theta}\cong\Gamma(\mathbb{T}^{2},M_{m}(\mathbb{C}))$.

\[ U(x,y)=\left(\begin{array}{ccccc}
\exp\left(\frac{2\pi ix}{m}\right) &  &  &  & \\
 & \exp\left({\frac{2\pi i(n+x)}{m}}\right)  & & & \\
 & & \exp\left({\frac{2\pi i(2n+x)}{m}}\right) & \\
 & & & \ddots & \\
 & & & & \exp\left({\frac{2\pi i((m-1)n+x)}{m}}\right)
\end{array}\right), \]

\[ V(x,y)=\left(\begin{array}{ccccc}
 & \exp\left({\frac{2\pi i(n+y)}{m}}\right) & & &  \\
 & & \exp\left(\frac{2\pi iy}{m}\right) & &  \\
 & & & \ddots &  \\
 & & & & \exp\left(\frac{2\pi iy}{m}\right)  \\
\exp\left(\frac{2\pi iy}{m}\right) & & & &  \\
\end{array}\right), \hspace{.25in} x,y\in[0,1]. \]

\noindent They satisfy the canonical commutation relation relation
\[ U(x,y)V(x,y)=e^{2\pi i\theta}V(x,y)U(x,y) \]

\noindent for any $x,y\in[0,1]$. Taking $m^{th}$ powers give the toroidal coordinates

\[ U(x,y)^{m}=e^{2\pi ix}I  \text{ \ \  and  \ \ } V(x,y)^{m}=e^{2\pi iy}I. \]

\

\subsection{Irrational noncommutative tori} \label{S5.2}

The situation of a rational noncommutative torus is closely related to that of the commutative torus as we saw in the previous section. However, the case for an irrational noncommutative torus is far challenging to describe. If we try to mimic the construction of a local covering in section \ref{S5.1}, a natural choice for the quantum symmetry is $\mathbb{T}^{2}_{\theta}\rtimes G$ but this is in general not a Hopf algebroid over $\mathbb{T}^{2}_{\theta}$. The problem is that there are no nice maps $s,t:\mathbb{T}^{2}_{\theta}\longrightarrow\mathbb{T}^{2}_{\theta}\rtimes G$ with commuting images since $\mathbb{T}^{2}_{\theta}$ is centrally simple for $\theta$ irrational. In this section, we will construct stratified coverings of $\mathbb{T}^{2}_{\theta}$ instead.

\begin{exa}\label{exa9}
In the classical case, any finite covering of the 2-torus $\mathbb{T}^{2}$ is again a 2-torus. Such covering spaces are of the form

\[ \xymatrix@R=0.1cm{\mathbb{T}^{2} \ar@{->>}[rr]^-{p} & & \mathbb{T}^{2} \\
(z_{1},z_{2}) \ar@{|->}[rr] & & (z_{1}^{n},z_{2}^{m}) } \]

\noindent and whose associated deck transformation group is $G=\mathbb{Z}/n\times\mathbb{Z}/m$. In the noncommutative set up, there is no reason for a covering of a noncommutative torus torus to be a noncommutative torus as well. This is easily seen with comparison with the noncommutative point having more than one connected covering space. Let us look at coverings of an irrational noncommutative torus $\mathbb{T}^{2}_{\theta}$ which are themselves noncommutative tori.

Let $0<\theta\in\mathbb{R}$ be an irrational number. Let $\mathbb{T}^{2}_{\theta}$ be the universal $C^{*}$-algebra generated by unitaries $U$ and $V$ satisfying $UV=e^{2\pi i\theta}VU$. It is well known that $\mathbb{T}^{2}_{\theta}$ is simple. The $K$-theory groups of $\mathbb{T}^{2}_{\theta}$ are $K_{0}(\mathbb{T}^{2}_{\theta})\cong K_{1}(\mathbb{T}^{2}_{\theta})\cong\mathbb{Z}^{2}$. More precisely, $K_{0}(\mathbb{T}^{2}_{\theta})\cong \mathbb{Z}+\theta\mathbb{Z}$ as an ordered group. For the purpose of what follows, we will say that two irrational numbers $\theta$ and $\eta$ are of the \textit{same type} if $\theta=n+m\eta$ for some integers $n,m$.

Consider an injective unital $C^{*}$-morphism $\mathbb{T}^{2}_{\theta}\stackrel{j}{\longrightarrow}\mathbb{T}^{2}_{\eta}$. There is an induced map $\mathbb{Z}+\theta\mathbb{Z}\stackrel{j_{*}}{\longrightarrow} \mathbb{Z}+\phi\mathbb{Z}$ in $K_{0}$, a map of ordered groups. Without loss of generality, we may assume $0<\theta<1$. Let $j_{*}(\theta)=n+m\eta$ for some integers $n,m$. By unitality of $j$, we have $j_{*}(1)=1$. We claim that $\theta$ and $\eta$ are of the same type. Suppose otherwise. In particular, this implies that $n+m\eta\neq\theta$. Without loss of generality, assume $n+m\eta>\theta$. Then, there is an integer $N$ such that $N\theta<M<N(n+m\eta)$ for some integer $M$. Thus, $N\theta<M$ and $M<N(n+m\eta)$. This implies that $N\theta<M$ and $\phi(M)<\phi(N\theta)$, which contradicts the fact that $\phi$ is order-preserving.

Using theorem 3.2.6 and proposition 3.2.7 in \cite{b010}, any injective $*$-homomorphism $\mathbb{T}^{2}_{\theta}\stackrel{\phi}{\longrightarrow}\mathbb{T}^{2}_{\eta}$ is approximately unitarily equivalent to an injective $*$-map $\mathbb{T}^{2}_{\theta}\stackrel{\alpha}{\longrightarrow}\mathbb{T}^{2}_{\eta}$ with $K_{1}\alpha:\mathbb{Z}^{2}\longrightarrow\mathbb{Z}^{2}$, $(x,y)\mapsto (n_{1}x+m_{1}y,n_{2}x+m_{2}y)$. In particular,

\[ \mathbb{T}^{2}_{\theta}\stackrel{\alpha}{\longrightarrow}\mathbb{T}^{2}_{\eta}, \hspace{.25in} U\mapsto P^{n_{1}}Q^{m_{1}}, \hspace{.25in} V\mapsto P^{n_{2}}Q^{m_{2}} \]

\noindent does the job. Here, $P$ and $Q$ are the unitary generators of $\mathbb{T}^{2}_{\eta}$. Since $UV=e^{2\pi i \theta}VU$, we have $\alpha(U)\alpha(V)=e^{2\pi i \theta}\alpha(V)\alpha(U)$. This implies that $e^{2\pi i(\theta-(n_{1}m_{2}-n_{2}m_{1})\eta)}=1$, and hence $\theta-(n_{1}m_{2}-n_{2}m_{1})\eta\in\mathbb{Z}$. This verifies that $\theta$ and $\eta$ are of the same type and at the same time gives the multiplier $N=n_{1}m_{2}-n_{2}m_{1}$. Let $G=\mathbb{Z}^{2}/\left\langle \left(n_{1},m_{1}\right),\left(n_{2},m_{2}\right)\right\rangle$, a group of order $N$. Let $\mathcal{H}=C\left(G\right)$, the Hopf algebra dual to $\mathbb{C}G$.

\begin{figure}[ht!]
\centering
\includegraphics[width=145mm]{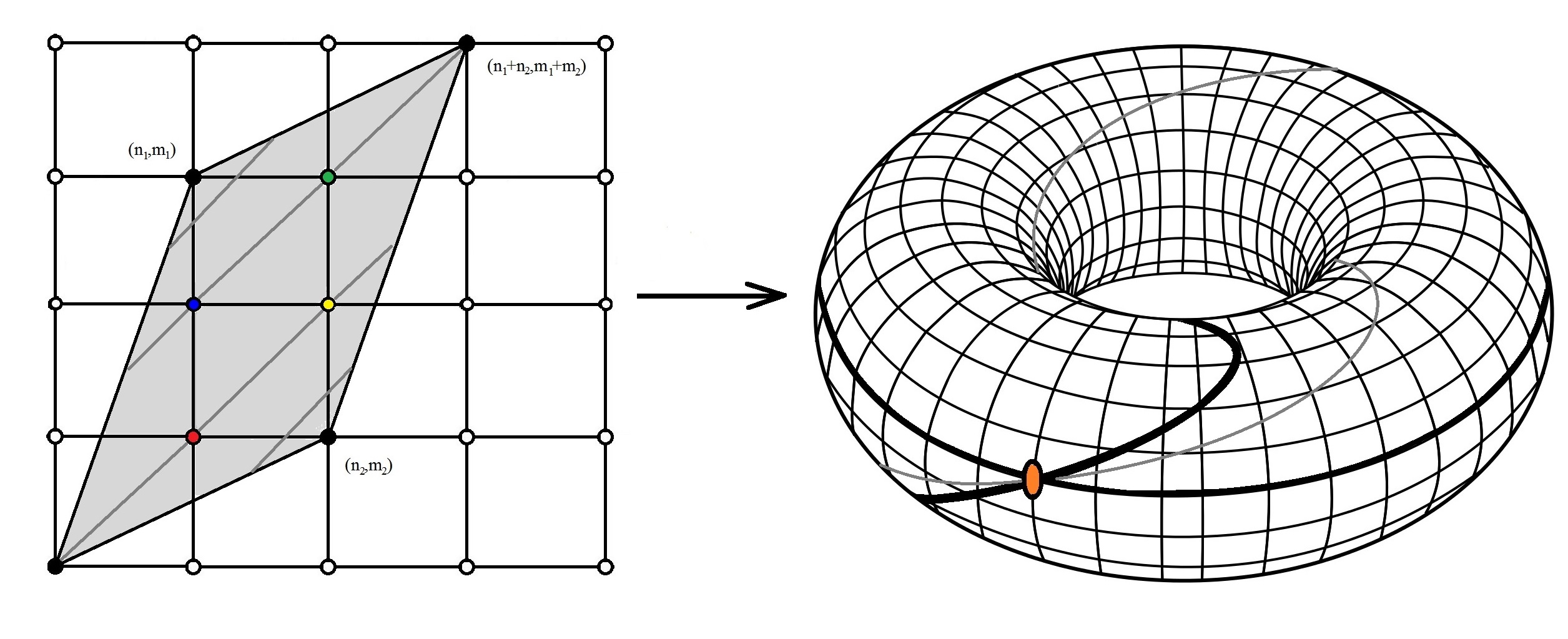}
\caption{Action of $G$ on $\mathbb{T}^{2}_{\eta}$}
\end{figure}

\

\noindent Let us show that $G$ acts on $\mathbb{T}^{2}_{\eta}$ with invariants $\mathbb{T}_{\theta}^{2}$ and hence, $\mathcal{H}$ coacts on $\mathbb{T}^{2}_{\eta}$ with coinvariants $\mathbb{T}^{2}_{\theta}$. Consider a fundamental domain for $G$. One can for example take the integral region in $\mathbb{Z}^{2}$ inside the parallelogram with vertices $(0,0),(n_{1},m_{1}),(n_{2},m_{2})$ and $(n_{1}+n_{2},m_{1}+m_{2})$ including $(0,0)$. This fundamental region can be identified with the Pontryagin dual $\widehat{G}$ of $G$. As an $\mathbb{T}^{2}_{\theta}$-module, $\mathbb{T}^{2}_{\eta}$ is freely generated by elements of the form $P^{n}Q^{m}$ where $(n,m)\in \widehat{G}$. Consider the canonical pairing $\left\langle ,\right\rangle:G\times\widehat{G}\longrightarrow\mathbb{S}^{1}$. Then $G$ acts on $\mathbb{T}^{2}_{\eta}$ by algebra isomorphisms defined for all $(n,m)\in \mathbb{Z}^{2}$ by

\[ (i,j)\cdot P^{n}Q^{m} = \left\langle (i,j),(n,m)\right\rangle P^{n}Q^{m}, \hspace{.25in} \text{for } (i,j)\in G. \]

\noindent Note that an element of $\mathbb{T}^{2}_{\eta}$ is invariant with this action precisely when $(n,m)$ is in the integral span of $(n_{1},m_{1})$ and $(n_{2},m_{2})$. This shows that the space of invariants is $\mathbb{T}^{2}_{\theta}$. This proves our claim. To show that the extension $\mathbb{T}^{2}_{\theta}\subseteq\mathbb{T}^{2}_{\eta}$ is $\mathcal{H}$-Galois, we have to check that the following linear map is an isomorphism.

\[ \xymatrix{\mathbb{T}^{2}_{\eta}\tens{\mathbb{T}^{2}_{\theta}}\mathbb{T}^{2}_{\eta} \ar[rr] & & \mathbb{T}^{2}_{\eta} \otimes \mathbb{C}G } \]

But this is immediate from the fact that $G$ acts freely and transitively on the $\mathbb{T}^{2}_{\theta}$-module generators of $\mathbb{T}^{2}_{\eta}$. This gives us a stratified covering $(\mathbb{T}^{2}_{\eta},\mathcal{H})$ of $\mathbb{T}^{2}_{\theta}$ with stratification $\mathbb{C}\subseteq \mathbb{T}^{2}_{\theta}$. $\Box$
\end{exa}

\begin{exa}\label{exa10}
Let us construct another stratified covering of $\mathbb{T}^{2}_{\theta}$. Let $n\in\mathbb{N}$ and let

\[ B=\mathbb{T}^{2}_{\theta/n}=C^{*}\left\langle U,V|U^{*}U=UU^{*}=1=V^{*}V=VV^{*},UV=e^\frac{2\pi i\theta}{n}VU\right\rangle \]

\noindent and let $A$ be the $C^{*}$-subalgebra of $B$ generated by $U$ and $V^{n}$. Note that $A\cong\mathbb{T}^{2}_{\theta}$. Let $A^{'}=C^{*}\left\langle U\right\rangle\subseteq A$. Note that $A^{'}\cong C(S^{1})$. Consider the Hopf algebra $\mathcal{H}=C(G,A^{'})$ where $G=\left\{1,\zeta,\zeta^{2},...,\zeta^{n-1}\right\}$, the group of $n^{th}$ roots of unity. $G$ acts on $\mathbb{T}^{2}_{\theta/n}$ as follows: $\zeta\cdot U=U$  and $\zeta\cdot V=\zeta V$. This action extends to an action of the Hopf algebra $A^{'}G$ with invariants $A$. Thus, $\mathcal{H}$ coacts on $\mathbb{T}^{2}_{\theta/n}$ with coinvariants $\mathbb{T}^{2}_{\theta}$. Using similar argument as the previous example, $A\subseteq B$ is an $\mathcal{H}$-Galois extension. This gives us a stratified covering of $\mathcal{T}^{2}_{\theta}$ with stratification $A^{'}\cong C(S^{1})$. $\Box$
\end{exa}

\subsection{Local, stratified, and uniform coverings}\label{S5.3}

Let us describe the contrast between local and stratified coverings. We aim to give a geometric intuition behind such stratifications and we will be less precise in doing so. First, note that local coverings can be regarded as a stratified coverings whose stratification is trivial (i.e., stratification by points). However, it will be useful to use local as we shall see soon.

In sections \ref{S5.1} and \ref{S5.2} we have constructed coverings of noncommutative tori with stratifications $A^{'}=A$, $A^{'}=C(S^{1})$ and $A^{'}=\mathbb{C}$. Pretending $A$ has points, these stratifications correspond to geometric stratifications illustrated in figure 4.

\begin{figure}[ht!]
\centering
\includegraphics[width=145mm]{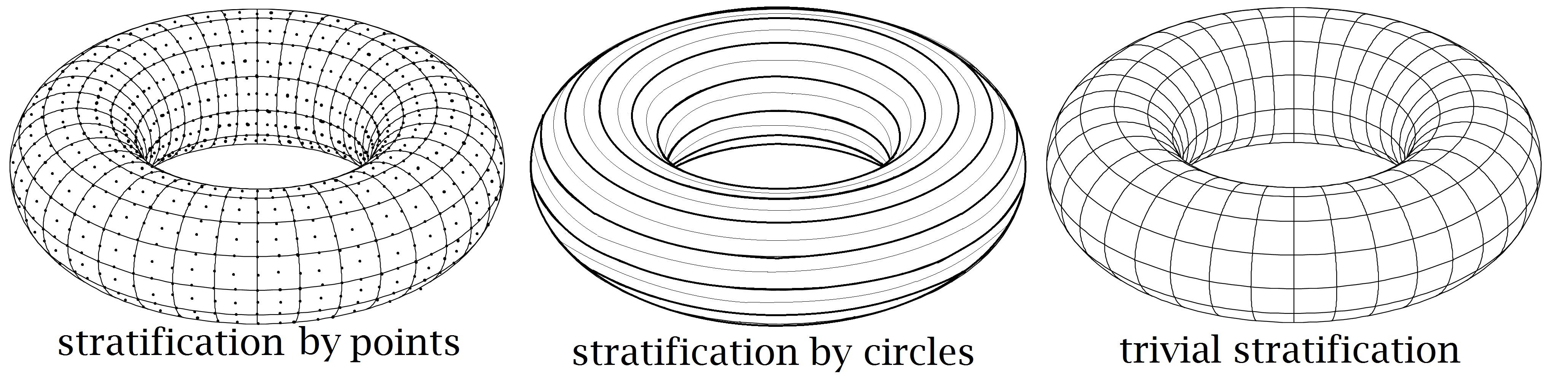}
\caption{Geometric stratifications associated with $A^{'}=A$, $A^{'}=C(S^{1})$ and $A^{'}=\mathbb{C}$.}
\end{figure}

A covering $(B,\mathcal{H})$ of $A$ with stratification $A^{'}\subseteq A$, by definition, has its quantum symmetry defined over $A^{'}$. By the duality between noncommutative spaces and algebras, the inclusion $A^{'}\subseteq A$ induces a surjection $\xymatrix{\widehat{A} \ar@{->>}[r] & \widehat{A^{'}}}$. This suggests that the quantum symmetry varies among the leaves of the stratification defined by $\xymatrix{\widehat{A} \ar@{->>}[r] & \widehat{A^{'}}}$ but remain constant within the leaves. As a concrete illustration, let us consider coverings of the (commutative) torus $\mathbb{T}^{2}$ with stratifications $A^{'}=C(\mathbb{T}^{2})$, $A^{'}=C(S^{1})$ and $A^{'}=\mathbb{C}$. The covering with stratification $A^{'}=C(\mathbb{T}^{2})$ has its quantum symmetry a Hopf algebroid $\mathcal{H}$ defined over the commutative algebra $C(\mathbb{T}^{2})$. If $C(\mathbb{T}^{2})$ is central in $\mathcal{H}$ then $\mathcal{H}$ is a bundle of complex Hopf algebroids over $\mathbb{T}^{2}$. These fiber Hopf algebroids need not be isomorphic. This suggest that the quantum symmetry can vary over $A^{'}=C(\mathbb{T}^{2})$. For the second case, $A^{'}=C(S^{1})$ using the same argument and assumptions imply that $\mathcal{H}$ is a bundle of complex Hopf algebroids over $S^{1}$ whose fibers may be nonisomorphic. These fibers Hopf algebroid varries among the fibers of $\xymatrix{\mathbb{T}^{2} \ar@{->>}[r]^-{p} & S^{1}}$ which defines the stratification. If $C(S^{1})$ is the largest subalgebra of $A=C(\mathbb{T}^{2})$ for which $\mathcal{H}$ is defined over then by the Galois condition, $\mathcal{H}$ must be constant along each fibers of $p$. The third case suggest that we have the same quantum symmetry $\mathcal{H}$ over each point of $\mathbb{T}^{2}$.

Meanwhile, uniform coverings are a special case of stratified coverings. Aside from $A^{'}=k$ we also require that $\mathcal{H}$ is a Hopf algebra. This in particular requires that the bialgebroid structures to coincide.

%%%%%%%%%%%%%%%%%%%%%%%%%%%%%%%%%%%%%%%%%%%%%%%%%%%%%%%%%%%%%%%%%%%%%%%%%%%%%%%%%%%%%%%%

\end{document}